\documentclass[draft=false]{scrartcl}%

\usepackage{amsfonts}
\usepackage{amsmath}
\usepackage{amsthm}
\usepackage{cite}
\usepackage[T1]{fontenc}
\usepackage{mathtools}
\usepackage{microtype}
\usepackage[automark]{scrpage2}%

\usepackage{142arxiv}

\theoremstyle{definition}

\theoremstyle{remark}

\mathtoolsset{mathic=true}

\numberwithin{equation}{section}
\title{On Matrix-Valued Stieltjes Functions with an Emphasis on Particular Subclasses}
\author{Bernd Fritzsche \and Bernd Kirstein \and Conrad M\"adler}
\date{}
\dedication{Dedicated to our friend and colleague Albrecht B\"ottcher on the occasion of his 60th birthday}

\begin{document}
\maketitle

\begin{abstract}
 The paper deals with particular classes of \tqqa{matrix-valued} functions which are holomorphic in \(\Cs\), where \(\lep\) is an arbitrary real number. These classes are generalizations of classes of holomorphic complex-valued functions studied by Kats and Krein~\zita{KK74} and by Krein and Nudelman~\zita{MR0458081}. The functions are closely related to truncated matricial Stieltjes problems on the interval \(\rhl\). Characterizations of these classes via integral representations are presented. Particular emphasis is placed on the discussion of the Moore-Penrose inverse of these matrix-valued functions.
\end{abstract}

\section{Introduction}\label{S0948}
In their papers~\zitas{MR3014199,MR2988005,1151}, the authors developed a simultaneous approach to the even and odd truncated matricial Hamburger moment problems. This approach was based on three cornerstones. One of them, namely the paper~\zita{MR2988005} is devoted to several functiontheoretical aspects concerning special subclasses of matrix-valued Herglotz-Nevanlinna functions. Now we are going to work out a similar simultaneous approach to the even and odd truncated matricial Stieltjes moment problems. Our approach is again subdivided into three steps. This paper is concerned with the first step which is aimed at a closer analysis of several classes of holomorphic matrix-valued functions in the complex plane which turn out to be closely related to Stieltjes type matrix moment problems. In the scalar case, the corresponding classes were carefully studied by Kats/Krein~\zita{KK74} and Krein/Nudelman~\zitaa{MR0458081}{Appendix}. What concerns the treatment of several matricial and operatorial generalizations, we refer the reader to the monograph Arlinskii/Belyi/Tsekanovskii~\zita{MR2828331} and the references therein.

In order to give a precise formulation of the matricial moment problem standing in the background of our investigations, we introduce some notation. Throughout this paper, let \(p\)\index{p@\(p\)} and \(q\)\index{q@\(q\)} be positive integers. Let \(\C\)\index{c@\(\C\)} and \(\R\)\index{r@\(\R\)} be the set of all complex numbers and the set of all real numbers, respectively. Furthermore, let \(\NO\)\index{n@\(\NO\)} and \(\N\)\index{n@\(\N\)} be the set of all \tnn{} integers and the set of all positive integers, respectively. Further, for every choice of \(\alpha,\beta\in\R\cup\set{-\infty,+\infty}\), let \(\mn{\alpha}{\beta}\)\index{z@\(\mn{\alpha}{\beta}\)} be the set of all integers \(k\) for which \(\alpha\leq k\leq\beta\) holds. If \(\cX\) is a \tne{} set, then let \(\cX^\pxq\)\index{\(\cX^\pxq\)} be the set of all \tpqa{matrices} each entry of which belongs to \(\cX\), and \(\cX^p\)\index{\(\cX^p\)} is short for \(\cX^\xx{p}{1}\). The notations \(\CHq\)\index{c@\(\CHq\)}, \(\Cggq\)\index{c@\(\Cggq\)}, and \(\Cgq\)\index{c@\(\Cgq\)} stand for the subsets of \(\Cqq\) which are formed by the sets of \tH{}, \tnnH{}, and \tpH{} matrices, respectively. If \((\Omega,\gA)\) is a measurable space, then each countably additive mapping whose domain is \(\gA\) and whose values belong to \(\Cggq\) is called a \tnnH{} \tqqa{measure} on \((\Omega,\gA)\). Let \(\BsaR\)\index{b@\(\BsaR\)} (resp.\ \(\BsaC\)\index{b@\(\BsaC\)}) be the \(\sigma\)\nobreakdash-algebra of all Borel subsets of \(\R\) (resp.\ \(\C\)). For each \(\Omega\in\BsaR\setminus\set{\emptyset}\), let \(\BsaO\defeq\BsaR\cap\Omega\)\index{b@\(\BsaO\)}, and let \(\Mggqa{\Omega}\)\index{m@\(\Mggqa{\Omega}\)} be the set of all \tnnH{} \tqqa{measure}s on \((\Omega,\BsaO)\). Furthermore, for each \(\Omega\in\BsaR\setminus\set{\emptyset}\) and every \(\kappa\in\NOinf \), let \(\Mgguqa{\kappa}{\Omega}\)\index{m@\(\Mgguqa{\kappa}{\Omega}\)} be the set of all \(\sigma\in\Mggqa{\Omega}\) such that the integral \(\mom{j}{\sigma}\defeq\int_\Omega t^j\sigma(\dif t)\)\index{s@\(\mom{j}{\sigma}\)} exists for all \(j\in\mn{0}{\kappa}\). In this paper, for an arbitrarily fixed \(\alpha\in\R\), we study classes of \tqqa{matrix-valued} functions which are holomorphic in \(\Cs\). These classes turn out to be closely related via \tStieltjestransform{} with the following truncated matricial Stieltjes type moment problem:
\begin{description}
 \item[\mproblem{\iraa{\alpha}}{m}{=}] Let \(\alpha\in\R\), let \(m\in\NO\), and let \(\seq{\su{j}}{j}{0}{m}\) be a sequence of complex \tqqa{matrices}. Describe the set \(\Mggqaag{\iraa{\alpha}}{\seq{\su{j}}{j}{0}{m}}\)\index{m@\(\Mggqaag{\iraa{\alpha}}{\seq{\su{j}}{j}{0}{m}}\)} of all \(\sigma\in\Mgguqa{m}{\iraa{\alpha}}\) for which \(\suo{j}{\sigma}=\su{j}\) is fulfilled for each \(j\in\mn{0}{m}\).\index{m@\mproblem{\iraa{\alpha}}{m}{=}}
\end{description}
In a forthcoming paper, we will indicate how classes of holomorphic functions studied in this paper can be used to parametrize the set \(\Mggqaag{\iraa{\alpha}}{\seq{\su{j}}{j}{0}{m}}\).

This paper is organized as follows. In \rsec{S1237}, we introduce several classes of holomorphic \tqqa{matrix-valued} functions. A particular important role will be played by the class \(\SFqlep\) of \tlSF{s} of order \(q\), which was considered in the special case \(q=1\) and \(\lep=0\) by I.~S.~Kats and M.~G.~Krein in~\zita{KK74} (see \rdefn{D1242}). In \rsec{S1508}, we derive several integral representations for functions belonging to \(\SFqlep\) (see \rthmss{T1542}{368}). Furthermore we analyse the structure of ranges and null spaces of the values of functions belonging to \(\SFqlep\) (see \rprop{L374}). In \rsec{S1038}, we state characterizations of the membership of a function to the class \(\SFqlep\). In \rsec{S1242}, we investigate the subclass \(\SOFqlep\) (see \rnota{D1351} below) of \(\SFqlep\). It is shown in \rthm{T1434} that this class is formed exactly by those \tqqa{matrix-valued} functions defined in \(\Cs\) which can be written as \tStieltjestransform{} of a \tnnH{} \tqqa{measure} defined on the Borelian \(\sigma\)\nobreakdash-algebra of the interval \(\rhl\). In \rsec{S1044}, we investigate the Moore-Penrose inverses of the functions belonging to \(\SFqlep\). In particular, we show that the Moore-Penrose inverse \(F^\mpi\) of a function \(F\in\SFqlep\) is holomorphic in \(\Cs\) and that the function \(G\colon\Cs\to\Cqq\) defined by \(G(z)\defeq-(z-\lep)^\inv\ek{F(z)}^\mpi\) belongs to \(\SFqlep\) as well (see \rthm{L375}). The second main theme of \rsec{S1044} is concerned with the investigation of the class \(\SiFqlep\) (see \rnota{D0924}), which was considered in the special case \(q=1\) and \(\lep=0\) by I.~S.~Kats and M.~G.~Krein in~\zita{KK74}. The main result on the class \(\SiFqlep\) is \rthm{T1305} which contains an integral representation which is even new for the case \(q=1\) and \(\lep=0\). The application of \rthm{T1305} enables us to obtain much information about ranges, null spaces and Moore-Penrose inverses of the functions belonging to the class \(\SiFqlep\) (see \rprop{P1622} and \rthm{P0931}). In the remaining sections of this paper, we carry out corresponding investigations for dual classes of \tqqa{matrix-valued} functions which are related to an interval \(\blhl\). These classes occur in the treatment of a matrix moment problem~\mproblem{\blhl}{m}{=}, which is analogous to Problem~\mproblem{\rhl}{m}{=} formulated above. In \rapp{A1608}, we summarize some facts from the integration theory with respect to \tnnH{} \tqqa{measures}.

\section{On several classes of holomorphic matrix-valued functions}\label{S1237}

In this section, we introduce those classes of holomorphic \tqqa{matrix-valued} functions which form the central objects of this paper. For each \(A\in\Cqq\), let \(\re A\defeq\frac{1}{2}(A+A^\ad)\)\index{r@\(\re A\)} and \(\im A\defeq\frac{1}{2\iu}(A-A^\ad)\)\index{i@\(\im A\)} be the real part of \(A\) and the imaginary part of \(A\), respectively. Let \(\ohe\defeq\setaa{z\in\C}{\im z>0}\)\index{p@\(\ohe\)} and \(\uhe\defeq\setaa{z\in\C}{\im z<0}\)\index{p@\(\uhe\)} 
be the open upper half plane and open lower half plane of \(\C\), respectively. The first two dual classes of holomorphic matrix-valued functions, which are particularly important for this paper, are the following.
\bdefnl{D1242}
 Let \(\lep\in\R\) and let \(F\colon\Cs\to\Cqq\). Then \(F\) is called a \emph{\tlSFq{}} if \(F\) satisfies the following three conditions:
 \bAeqi{0}
  \il{D1242.I} \(F\) is holomorphic in \(\Cs\).
  \il{D1242.II} For all \(w\in\ohe\), the matrix \(\im\ek{F(w)}\) is \tnnH{}.
  \il{D1242.III} For all \(w\in\lhl\), the matrix \(F(w)\) is \tnnH{}.
 \eAeqi
 We denote by \(\SFqlep\)\index{s@\(\SFqlep\)} the set of all \tlSF{s} of order \(q\).
\edefn

\bexal{E1256}
 Let \(\lep\in\R\) and let \(A,B\in\Cggq\). Let \(F\colon\Cs\to\Cqq\) be defined by \(F(z)\defeq A+\frac{1}{\lep-z}B\). Since \(\im F(z)=\frac{\im z}{\abs{\lep-z}^2}B\) holds true for each \(z\in\Cs\), we have \(F\in\SFqlep\).
\eexa

\bdefnl{D1253}
 Let \(\rep\in\R\) and let \(G\colon\Ct\to\Cqq\). Then \(G\) is called a \emph{\trTFq{}} if \(G\) fulfills the following three conditions:
 \bAeqi{0}
  \il{D1253.I} \(G\) is holomorphic in \(\Ct\).
  \il{D1253.II} For all \(w\in\ohe\), the matrix \(\im\ek{G(w)}\) is \tnnH{}.
  \il{D1253.III} For all \(w\in\brhl\), the matrix \(-G(w)\) is \tnnH{}.
 \eAeqi
 We denote by \(\TFqrep\)\index{s@\(\TFqrep\)} the set of all \trTF{s} of order \(q\).
\edefn

\bexal{E0913}
 Let \(\rep\in\R\) and let \(A,B\in\Cggq\). Let \(G\colon\Ct\to\Cqq\) be defined by \(G(z)\defeq-A+\frac{1}{\rep-z}B\). Since \(\im G(z)=\frac{\im z}{\abs{\rep-z}^2}B\) holds true for each \(z\in\Ct\), we have \(G\in\TFqrep\).
\eexa

The particular functions  belonging to the class \(\SFqlep\) resp.\ \(\TFqrep\), which were introduced in \rexa{E1256} (resp.\ \rexa{E0913}), were called by V.~E.~Katsnelson~\zita{MR2805421} \emph{special} functions belonging to \(\SFqlep\) (resp.\ \(\TFqrep\)), whereas all remaining functions contained in \(\SFqlep\) (resp.\ \(\TFqrep\)) were called \emph{generic} functions belonging to \(\SFqlep\) (resp.\ \(\TFqrep\)). In the case \(q=1\) and \(\lep=0\) the general theory of multiplicative representations of functions belonging to \(\SFuu{1}{0}\) was treated in detail by Aronszajn/Donoghue~\zita{MR0168769}.

It should be mentioned that Yu.~M.~Dyukarev and V.~E.~Katsnelson studied in~\zitas{MR686076,MR645305,MR752057} an interpolation problem for functions belonging to the class \(\TFuu{q}{0}\). Their approach was based on V.~P.~Potapov's method of fundamental matrix inequalities. V.~E.~Katsnelson~\zita{MR2805421} used the class \(\TFuu{1}{0}\) to construct Hurwitz stable entire functions.

First more or less obvious interrelations between the above introduced classes of functions can be described by the following result. For each \(\lep\in\R\), let the mapping \(T_\lep\colon\R\to\R\) be defined by \(T_\lep(x)\defeq x+\lep\)\index{t@\(T_\lep(x)\)}. If \(\cX\), \(\cY\), and \(\cZ\) are \tne{} sets with \(\cZ\subseteq\cX\) and if \(f\colon\cX\to\cY\) is a mapping, then we will use \(\Rstr_\cZ f\)\index{r@\(\Rstr_\cZ f\)} to denote the restriction of \(f\) onto \(\cZ\).

\breml{366}
 \benui
  \item If \(\lep\in\R\) and if \(F\colon\Cs \to \Cqq \), then \(F\) belongs to \(\SFqlep \) if and only if the function  \(F \circ \Rstr_{\C\setminus [0,+\infty)} T_\lep\) belongs to the class \(\SFuu{q}{0}\).
  \item If \(\lep\in\R\) and if \(F\colon\C\setminus [0,+\infty)\to \Cqq \), then \(F\in\SFuu{q}{0}\) if and only if \(F \circ \Rstr_{\Cs } T_{-\lep} \in\SFqlep \).
 \eenui
\erem

\breml{R1341}
 \benui
  \item If \(\rep\in\R\) and if \(G\colon\Ct\to \Cqq \), then \(G\) belongs to \(\TFqrep \) if and only if the function  \(G\circ \Rstr_{\C\setminus(-\infty,0]} T_\rep\) belongs to the class \(\TFuu{q}{0}\).
  \item If \(\rep\in\R\) and if \(G\colon\C\setminus(-\infty,0]\to \Cqq \), then \(F\in\TFuu{q}{0}\) if and only if \(F \circ \Rstr_{\Ct} T_{-\rep} \in\TFqrep \).
 \eenui
\erem

Now we introduce particular subclasses of \(\SFqlep\) and \(\TFqrep\), which will turn out to be important in studying matricial versions of the Stieltjes moment problem. If \(A\in\Cpq\), then we denote by \(\normE{A}\defeq\sqrt{\tr\rk{A^\ad A}}\)\index{e@\(\normE{A}\)} the Euclidean norm of \(A\). Special attention will be put in the sequel also to the following subclasses of the classes of holomorphic matrix-valued functions introduced in \rdefnss{D1242}{D1253}.

\bnotal{D1351}
 Let \(\lep\in\R\). Then let \(\SOFqlep\)\index{s@\(\SOFqlep\)} be the class of all \(S\) which belong to \(\SFqlep\) and which satisfy
 \beql{Nr.YS}
  \sup_{y\in[1,+\infty)}y\normE*{S(\iu y)}
  <+\infty.
 \eeq
 Furthermore, let \(\TOFuu{q}{\lep}\)\index{s@\(\TOFuu{q}{\lep}\)} be the class of all \(S\in\TFuu{q}{\lep}\) which satisfy \eqref{Nr.YS}.
\enota

\breml{R1621}
 Let \(\lep\in\R\). If \(S\in\SOFqlep\) or if \(S\in\TOFuu{q}{\lep}\), then \(\lim_{y\to+\infty}S(\iu y)=\Oqq\).
\erem

\rrem{R1621} leads us to the following classes, which will play an important role in the framework of a truncated matricial Stieltjes moment problems.

\bnotal{D1026}
 Let \(\lep\in\R\). Then by \(\SdFqlep\)\index{s@\(\SdFqlep\)} (resp.\ \(\TdFuu{q}{\lep}\)\index{s@\(\TdFqrep\)}) we denote the set of all \(F\in\SFqlep\) (resp.\ \(\TdFuu{q}{\lep}\)) which satisfy \(\lim_{y\to+\infty}F(\iu y)=\Oqq\).
\enota

\breml{R1031}
 If \(\lep\in\R\), then \rrem{R1621} shows that \(\SOFqlep\subseteq\SdFqlep\) and \(\TOFuu{q}{\lep}\subseteq\TdFuu{q}{\lep}\).
\erem

\breml{R1036}
 Let \(\lep\in\R\) and let
 \[
  \mathcal{S}
  \in\set{\SFqlep,\SOFqlep,\SdFqlep,\TFuu{q}{\lep},\TOFuu{q}{\lep},\TdFuu{q}{\lep}}.
 \]
 Then \(F\in\mathcal{S}\) if and only if \(F^\tra\in\mathcal{S}\).
\erem

\breml{R1352}
 \benui
  \il{R1352.a} Let \(\lep\in\R\) and let \(F\colon\Cs \to\Cqq\) be a matrix-valued function. Then it is readily checked that \(F\) belongs to one of the classes \(\SFqlep\), \(\SOFqlep \), and \(\SdFqlep\) if and only if for each \(u\in\Cq\) the function \(u^\ad Fu\) belongs to the corresponding classes \(\SFuu{1}{\lep}\), \(\SOFuu{1}{\lep}\), and \(\SdFuu{1}{\lep}\), respectively.
  \il{R1352.b} Let \(\rep\in\R\) and let \(G\colon\Ct\to\Cqq\) be a matrix-valued function. Then it is readily checked that \(G\) belongs to one of the classes \(\TFqrep\), \(\TOFqrep \), and \(\TdFqrep \) if and only if for each \(u\in\Cq\) the function \(u^\ad Gu\) belongs to the corresponding classes \(\TFuu{1}{\rep}\), \(\TOFuu{1}{\rep}\), and \(\TdFuu{1}{\rep}\), respectively.
 \eenui
\erem

In order to get integral representations and other useful information about the one-sided Stieltjes functions of order \(q\), we exploit the fact that these classes of functions can be embedded via restriction to the upper half plane \(\ohe\) to the well-studied Herglotz-Nevanlinna class \(\RFq\) of all Herglotz-Nevanlinna functions in \(\ohe\). A matrix-valued function \(F\colon\ohe\to\Cqq\) is called a \emph{\tqqa{Herglotz-Nevanlinna} function in \(\ohe\)} if \(F\) is holomorphic in \(\ohe\) and satisfies the condition \(\im\ek{F(w)}\in\Cggq\) for all \(w\in\ohe\). For a comprehensive study on the class \(\RFq\), we refer the reader to the paper~\zita{MR1784638} by F.~Gesztesy and E.~R.~Tsekanovskii and to the paper~\zita{MR2988005}. In particular, in~\zita{MR2988005} one can find a detailed discussion of the holomorphicity properties of the Moore-Penrose pseudoinverse of matrix-valued Herglotz-Nevanlinna functions. Before we recall the well-known characterization of the class \(\NFq\) via Nevanlinna parametrization, we observe that, for each \(\nu\in\MggqR\) and each \(z\in\C\setminus\R\), the function \(\hu{z}\colon\R\to\C\) defined by \(\hua{z}{t}\defeq(1+tz)/(t-z)\) belongs to \(\LaaaC{\R}{\BAR}{\nu}\).

\bthml{T1554}
 \benui
  \il{T1554.a} Let \(F\in\NFq\). Then there are unique matrices \(\Npa\in\CHq\) and \(\Npb\in\Cggq\) and a unique \tnnH{} measure \(\Npnu\in\MggqR\) such that
  \begin{align}\label{B1}
   F(z)
   &=\Npa+z\Npb+\int_\R\frac{1+tz}{t-z}\Npnu(\dif t)&\text{for each }z&\in\ohe.
  \end{align}
  \item If \(\Npa\in\CHq\), if \(\Npb\in\Cggq\), and if \(\Npnu\in\MggqR\), then \(F\colon\ohe\to\Cqq\) defined by \eqref{B1} belongs to the class \(\NFq\).
 \eenui
\ethm

For each \(F\in\NFq\), the unique triple \((\Npa,\Npb,\Npnu)\in\CHq\times\Cggq\times\MggqR\) for which the representation \eqref{B1} holds true is called the \emph{Nevanlinna parametrization of \(F\)} and we will also write \((\alphaF,\betaF,\nuF)\) for \((\Npa,\Npb,\Npnu)\).

\breml{R1049}
 Let \(F\in\RFq\) with Nevanlinna parametrization \((\alphaF,\betaF,\nuF)\). Then it is immediately seen that \(F^\tra\) belongs to \(\RFq\) and that \((\alphau{F^\tra},\betau{F^\tra},\nuu{F^\tra})=(\alphaF^\tra,\betaF^\tra,\nuF^\tra)\).
\erem

\breml{R1602}
 Let \(F\in\RFq\) with Nevanlinna parametrization \((\Npa,\Npb,\Npnu)\). In view of \rthm{T1554}, we have then \(-\ko z\in\ohe\) and
 \[
  -\ek*{F(-\ko z)}^\ad
  =-\Npa+z\Npb+\int_\R\frac{1-tz}{-t-z}\Npnua{\dif t}
  =-\Npa+z\Npb+\int_\R\frac{1+tz}{t-z}\theta(\dif t)
 \]
 for all \(z\in\ohe\), where \(\theta\) is the image measure of \(\Npnu\) under the reflection \(t\mapsto-t\) on \(\R\). Because of \(-\Npa\in\CHq\) and \(\Npb\in\Cggq\), \rthm{T1554} yields then that \(G\colon\ohe\to\Cqq\) defined by \(G(z)\defeq-\ek{F(-\ko z)}^\ad\) belongs to \(\RFq\) and that the Nevanlinna parametrization of \(G\) is given by \((-\Npa,\Npb,\theta)\).
\erem

From \rdefnss{D1242}{D1253} we see immediately that \(\setaa{\Rstr_{\ohe}F}{F\in\SFqlep}\subseteq\RFq\) and \(\setaa{\Rstr_{\ohe}G}{G\in\TFqrep}\subseteq\RFq\). Now we analyse the Nevanlinna parametrizations of the restrictions on \(\ohe\) of the members of these classes of functions.

\bpropl{HCT}
 Let \(\lep\in\R\) and let \(F\in\SFqlep \). Then the Nevanlinna parametrization \((\Npa,\Npb,\Npnu)\) of \(\Rstr_{\ohe } F\) fulfills
 \begin{align}\label{Nr.Neu1}
  \NpnuA{(-\infty, \alpha)}&= \Oqq,&
  \Npb&= \Oqq,&
  &\text{and}&
  \Npnu&\in\MgguqR{1}. 
 \end{align}
 In particular, for each \(z\in\Cs \), then 
 \beql{HCT.B1}
  F(z)
  =\Npa + \int_{\rhl} \frac{1+tz}{t-z} \Npnua{\dif t}.
 \eeq
\eprop
\bproof
 From~\zitaa{MR2222521}{\cprop{8.3}} and its proof we obtain \(\Npnua{\lhl}= \Oqq\) and
 \beql{HCT.1}
  F(z)
  =\Npa+z\Npb+\int_{\rhl}\frac{1+tz}{t-z}\nu(\dif t)
 \eeq
 for all \(z\in\Cs\). Let \(u\in\Cq\). Then \(u^\ad\Npnu u\in\Mggoa{1}{\R}\) and \eqref{HCT.1} yields
 \beql{HCT.2}
  u^\ad F(x)u
  =u^\ad\Npa u+xu^\ad\Npb u+\int_{\rhl}\frac{1+tx}{t-x}\rk{u^\ad\nu u}(\dif t)
 \eeq
 for all \(x\in\lhl\). Let \(\lep_1\defeq\min\set{\lep-1,-1}\) and \(\lep_2\defeq\max\set{\lep+1,1}\). Since \(F(x)\in\Cggq\) for all \(x\in(-\infty,\lep_1)\), and furthermore \((1+tx)/(t-x)<0\) for all \(x\in(-\infty,\lep_1)\) and all \(t\in(\lep_2,+\infty)\), we  conclude from \eqref{HCT.2} then
 \[
  -xu^\ad\Npb u
  \leq u^\ad\Npa u+\int_{[\lep,\lep_2]}\frac{1+tx}{t-x}\rk{u^\ad\nu u}(\dif t)
 \]
 for all \(x\in(-\infty,\lep_1)\). One can easily check that there exists a constant \(L_\lep\in\R\) depending only on \(\lep\) such that \(\abs{(1+tx)/(t-x)}\leq L_\lep\) for all \(x\in(-\infty,\lep_1)\) and all \(t\in[\lep,\lep_2]\). For all \(x\in(-\infty,\lep_1)\), hence
 \beql{HCT.3}
  \abs*{\int_{[\lep,\lep_2]}\frac{1+tx}{t-x}\rk{u^\ad\nu u}(\dif t)}
  \leq L_\lep\cdot\rk{u^\ad\nu u}\rk*{[\lep,\lep_2]}
  <+\infty.
 \eeq
 Setting \(K\defeq u^\ad\Npa u+L_\lep\cdot\rk{u^\ad\nu u}\rk{[\lep,\lep_2]}\), we get then \(-xu^\ad\Npb u\leq K<+\infty\). In view of \(\Npb\in\Cggq\), we have \(u^\ad\Npb u\geq0\), where \(u^\ad\Npb u>0\) is impossible, since \(-xu^\ad\Npb u\) would then tend to \(+\infty\) as \(x\) tends to \(-\infty\). Thus, \(u^\ad\Npb u=0\). Since \(u\in\Cq\) was arbitrarily chosen, we get \(B=\Oqq\) and, in view of \eqref{HCT.1}, thus \eqref{HCT.B1} holds true for all \(z\in\Cs\). Taking additionally into account \(F(x)\in\Cggq\) for all \(x\in(-\infty,\lep_1)\), we conclude from \eqref{HCT.2} and \eqref{HCT.3} furthermore \(\int_{[\lep_2,+\infty)}\ek{-\rk{1+tx}/\rk{t-x}}\rk{u^\ad\nu u}(\dif t)\leq K<+\infty\) for all \(x\in(-\infty,\lep_1)\). Now we consider an arbitrary sequence \(\seq{x_n}{n}{1}{\infty}\) from \((-\infty,\lep_1)\) with \(\lim_{n\to\infty}x_n=-\infty\). Then \(-(1+tx_n)/(t-x_n)>0\) for all \(n\in\N\) and all \(t\in[\lep_2,+\infty)\) and, furthermore, \(\abs{t}=\liminf_{n\to\infty}\ek{-(1+tx_n)/(t-x_n)}\) for all \(t\in[\lep_2,+\infty)\). The application of Fatou's lemma then yields
 \[
  \int_{[\lep_2,+\infty)}\abs{t}\rk{u^\ad\nu u}(\dif t)
  \leq\liminf_{n\to\infty}\int_{[\lep_2,+\infty)}\rk*{-\frac{1+tx}{t-x}}\rk{u^\ad\nu u}(\dif t)
  \leq K<+\infty.
 \]
 Since \(\rk{u^\ad\nu u}\rk{\lhl}=0\) and since \(\int_{[\lep,\lep_2]}\abs{t}\rk{u^\ad\nu u}(\dif t)\) is finite, we conclude then \(\int_\R\abs{t}\rk{u^\ad\nu u}(\dif t)<+\infty\). Because \(u\in\Cq\) was arbitrarily chosen, we get \(\Npnu\in\MgguqR{1}\).
\eproof

\section{Integral representations for the class $\SFqlep$}\label{S1508}
The main goal of this section is to derive some integral representations for \tlSF{s} of order \(q\).

\bthml{T1542}
 Let \(\lep\in\R\) and let \(F\colon\Cs\to\Cqq\). Then:
 \benui
  \il{T1542.a} Suppose \(F\in\SFqlep \). Denote by \((\Npa,\Npb,\Npnu)\) the Nevanlinna parametrization of \(\Rstr_{\ohe } F\). Then \(\tilde\Npnu\defeq\Rstr_{\Bsarhl}\Npnu\) belongs to \(\Mgguqrhl{1}\) and there is a unique pair \((C,\eta)\) belonging to \(\Cggq\times\Mgguqrhl{1}\) such that
  \beql{T1542.B1}
   F(z)
   =C+\int_{\rhl}\frac{1+t^2}{t-z}\eta(\dif t)
  \eeq
  for all \(z\in\Cs\), namely \(C=\Npa-\int_{\rhl}t\tilde\Npnu(\dif t)\) and \(\eta=\tilde\Npnu\). Furthermore, \(C=\lim_{r\to+\infty}F(\lep+r\ec^{\iu\phi})\) for all \(\phi\in(\pi/2,3\pi/2)\).
  \il{T1542.b} Let \(C\in\Cggq\) and let \(\eta\in\Mgguqrhl{1}\) be such that \eqref{T1542.B1} holds true for all \(z\in\Cs\). Then \(F\) belongs to \(\SFqlep\).
 \eenui
\ethm
\bproof
 \eqref{T1542.a} From \rprop{HCT} we conclude \(\tilde\Npnu\in\Mgguqrhl{1}\). Let \(u\in\Cq\). Then \(\int_{\rhl}t\rk{u^\ad\tilde\Npnu u}\rk{\dif t}\) is finite. Obviously, for all \(z\in\C\) and all \(t\in\R\) with \(t\neq z\), we get
 \beql{T1542.0}
  \frac{1+tz}{t-z}+t
  =\frac{1+t^2}{t-z}.
 \eeq
 Thus, in view of \rprop{HCT}, we obtain
 \beql{T1542.2}
  \int_{\rhl}\frac{1+t^2}{t-z}\rk{u^\ad\tilde\Npnu u}\rk{\dif t}
  <+\infty
 \eeq
 and
 \beql{T1542.1}
  u^\ad F(z) u
  =u^\ad\Npa u+\int_{\rhl}\frac{1+t^2}{t-z}\rk{u^\ad\tilde\Npnu u}\rk{\dif t}-\int_{\rhl}t\rk{u^\ad\tilde\Npnu u}\rk{\dif t}
 \eeq
 for all \(z\in\Cs\). Since \(u\in\Cq\) was arbitrarily chosen, it follows \eqref{T1542.B1} for all \(z\in\Cs\) with \(C\defeq\Npa-\int_{\rhl}t\tilde\Npnu(\dif t)\) and \(\eta\defeq\tilde\Npnu\). Let \(\phi\in(\pi/2,3\pi/2)\). Then \(\cos\phi<0\). To show \(C=\lim_{r\to+\infty}F(\lep+r\ec^{\iu\phi})\), we consider an arbitrary sequence \(\seq{r_n}{n}{1}{\infty}\) from \(\R\) with \(r_n\geq1/\abs{\cos\phi}\) for all \(n\in\N\) and \(\lim_{n\to\infty}r_n=+\infty\). We have then \(\lim_{n\to\infty}(1+t^2)/(t-\lep-r_n\ec^{\iu\phi})=0\) for all \(t\in\rhl\). For all \(n\in\N\) and all \(t\in\rhl\), we get furthermore
 \[
  \abs{t-\lep-r_n\ec^{\iu\phi}}
  \geq t-\lep-r_n\cos\phi
  =t-\lep+r_n\abs{\cos\phi}
  \geq t-\lep+1
  \geq1
 \]
  and hence \(\abs{(1+t^2)/(t-\lep-r_n\ec^{\iu\phi})}\leq(1+t^2)/(t-\lep+1)\). Since, because of \eqref{T1542.2}, the integral \(\int_{\rhl}(1+t^2)/(t-\lep+1)\rk{u^\ad\tilde\Npnu u}\rk{\dif t}\) is finite, the application of Lebesgue's dominated convergence theorem yields \(\lim_{n\to\infty}\int_{\rhl}(1+t^2)/(t-\lep-r_n\ec^{\iu\phi})\rk{u^\ad\tilde\Npnu u}\rk{\dif t}=0\). From \eqref{T1542.1} we conclude then \(u^\ad Cu=\lim_{n\to\infty}u^\ad F(\lep+r_n\ec^{\iu\phi}) u\). Since \(u\in\Cq\) was arbitrarily chosen, we obtain \(C=\lim_{n\to\infty}F(\lep+r_n\ec^{\iu\phi})\). Taking into account \(F(x)\in\Cggq\) for all \(x\in(-\infty,\lep)\), we get with \(\phi=\pi\) in particular \(C\in\Cggq\).
 
 Now let \(C\in\Cggq\) and \(\eta\in\Mgguqrhl{1}\) be such that \eqref{T1542.B1} holds true for all \(z\in\Cs\). Then \(\chi\colon\BsaR\to\Cggq\) defined by \(\chi(M)\defeq\eta(M\cap\rhl)\) belongs to \(\MgguqR{1}\) and the matrix \(C+\int_\R t\chi(\dif t)\) is \tH{}. Using \eqref{T1542.0}, we infer from \eqref{T1542.B1} that the integral \(\int_\R(1+tz)/(t-z)\chi(\dif t)\) exists and that
 \[
  F(z)
  =C+\int_\R t\chi(\dif t)+z\cdot\Oqq+\int_\R\frac{1+tz}{t-z}\chi(\dif t)
 \]
 is fulfilled for all \(z\in\ohe\). \rthmp{T1554}{T1554.a} yields then \(C+\int_\R t\chi(\dif t)=\Npa\) and \(\chi=\Npnu\). Hence \(\eta=\tilde\Npnu\) and \(C=\Npa-\int_{\rhl}t\tilde\Npnu(\dif t)\).
 
 \eqref{T1542.b} Let \(C\in\Cggq\) and \(\eta\in\Mgguqrhl{1}\) be such that \eqref{T1542.B1} holds true for all \(z\in\Cs\). Using a result on holomorphic dependence of an integral on a complex parameter (see, \eg{}~\zitaa{MR2257838}{Ch.~IV, \S5, \cSatz{5.8}}), we conclude then that \(F\) is a matrix-valued function which is holomorphic in \(\Cs\). Furthermore,
 \[
  \im F(z)
  =\int_{\rhl}\im\rk*{\frac{1+t^2}{t-z}}\eta(\dif t)
  =\int_{\rhl}\frac{(1+t^2)\im z}{\abs{t-z}^2}\eta(\dif t)
  \in\Cggq
 \]
 for all \(z\in\ohe\) and
 \[
  F(x)
  =C+\int_{\rhl}\frac{1+t^2}{t-x}\eta(\dif t)
  \in\Cggq
 \]
 for all \(x\in\lhl\). Thus, \(F\) belongs to \(\SFqlep\).
\eproof

In the following, if \(\lep\in\R\) and \(F\in \SFqlep \) are given, then we will write \((\Cu{F},\etau{F})\)\index{c@\(\CF\)}\index{e@\(\etaF\)} for the unique pair \((C,\eta)\) belonging to \(\Cggq\times\Mgguqrhl{1}\) which fulfills \eqref{T1542.B1} for all \(z\in\Cs \). Furthermore, if \(A\) and \(B\) are complex \tqqa{matrices}, then we write \(A\lleq B\)\index{$A\lleq B$} or \(B\lgeq A\)\index{$B\lgeq A$} to indicate that the matrices \(A\) and \(B\) are \tH{} and that \(B-A\) is \tnnH{}.

\breml{R0833}
 Let \(\lep\in\R\) and let \(F\in\SFqlep\). For all \(x_1,x_2\in\lhl\) with \(x_1\leq x_2\), then \(\Oqq\leq F(x_1)\leq F(x_2)\), by virtue of \rthmp{T1542}{T1542.a}.
\erem

\breml{R1618}
 Let \(\lep\in\R\) and let \(z\in\Cs \). Then, for each \(\mu\in\Mggqrhl \), in view of the equation \((1+t-\lep)/(t-z)=1+(1+z-\lep)/(t-z)\), which holds for each \(t\in\rhl \), and \rlemp{63A}{63A.a}, one can easily see that the function \(h_{\lep,z}\colon\rhl \to\C\)\index{h@\(h_{\lep,z}\)} defined by \(h_{\lep,z}(t)\defeq(1+t-\lep)/(t-z)\) belongs to \(\Loaaaa{1}{\rhl }{\Bsarhl }{\mu}{\C}\).
\erem

\bleml{NN}
 Let \(\lep\in\R\) and let \(F\colon\Cs\to\Cqq\) be a continuous matrix-valued function such that \(\rk{\im F}(\ohe)\subseteq\Cggq\) and \(\rk{-\im F}(\uhe)\subseteq\Cggq\). Then \(F(x)=\re F(x)\) and \(\im F(x)=\Oqq\) for each \(x\in\lhl\).
\elem
\bproof
 Let \(x\in\lhl\). Then \(\seq{\im F(x+\iu/n)}{n}{1}{\infty}\) and \(\seq{-\im F(x-\iu/n)}{n}{1}{\infty}\) are sequences of \tnnH{} complex \tqqa{matrices} which converge to the \tnnH{} complex matrices \(\im F(x)\) and \(-\im F(x)\), respectively. Consequently, \(\im F(x)=\Oqq\), which implies \(F(x)=\re F(x)\).
\eproof

\bleml{N6N}
 Let \(\lep\in\R\), let \(\gamma\in\Cggq\), and let \(\mu\in\Mggqrhl\). Let \(F\colon\Cs\to\Cqq\) be defined by
 \beql{Nr.FM}
  F(z)
  = \gamma +\int_{\rhl} \frac{1+ t-\alpha}{t-z} \mu (\dif t).
 \eeq
 Then
 \beql{368.2}
  \re F(z)
  =\gamma+\int_{\rhl}\frac{1+t-\alpha}{\abs{t-z}^2}(t- \re z)\mu(\dif t)
 \eeq
 and
 \beql{368.1}
  \im F(z)
  =(\im z)\int_{\rhl} \frac{1+t-\alpha}{\abs{t-z}^2} \mu (\dif t)
 \eeq
 hold true for each \(z\in\Cs\). In particular,
 \begin{align*}
  \rk{\re F}(\lhelep)&\subseteq\Cggq,&
  \rk{\im F}(\ohe)&\subseteq\Cggq,&
  &\text{and}&
  \rk{-\im F}(\uhe)&\subseteq\Cggq.
 \end{align*}
 Furthermore, \(F(x)=\re F(x)\) and \(F(x)\in\Cggq\) for each \(x\in\lhl\).
\elem
\bproof
 For each \(z\in\Cs\) and each \(t\in\rhl\), we have
 \beql{IK}
  \re\rk*{\frac{1+t-\lep}{t-z}}
  =\frac{1+t-\lep}{\abs{t-z}^2}(t-\re z)
 \eeq
 and \(\im\ek{\rk{1+t-\lep}/\rk{t-z}}=\rk{\im z}\rk{1+t-\lep}/\abs{t-z}^2\).
 Taking into account \(\gamma\in\Cggq\), thus \eqref{368.2} and \eqref{368.1} follow for each \(z\in\Cs\). For each \(z\in\lhelep\) and each \(t\in\rhl\), the right-hand side of \eqref{IK} belongs to \([0,+\infty)\). Thus, \(\gamma\in\Cggq\) and \eqref{368.2} show that \(\re F(z)\) belongs to \(\Cggq\) for each \(z\in\lhelep\). Since \((1+t-\alpha)/\abs{t-z}^2\in[0,+\infty)\) for every choice of \(z\in\Cs\) and \(t\in\rhl\), from \eqref{368.1} we see that \(\im F(w)\in\Cggq\) for each \(w\in\ohe\) and \(-\im F(v)\in\Cggq\) for each \(v\in\uhe\) are fulfilled. Applying \rlem{NN} completes the proof.
\eproof

Now we give a further integral representation of the matrix-valued functions which belong to the class \(\SFqlep\). In the special case that \(q=1\) and \(\lep=0\) hold, one can find this result in~\zitaa{MR0458081}{Appendix}.
\bthml{368}
 Let \(\lep \in\R\) and let \(F\colon\Cs  \to \Cqq \). Then:
 \benui
 \il{368.a} If \(F \in \SFqlep \), then there are a unique matrix \(\gamma \in \Cggq\) and a unique \tnnH{} measure \(\mu \in\Mggqrhl \) such that \eqref{Nr.FM} holds true for each \(z\in\Cs \). Furthermore, \(\gamma=\CF\).
 \il{368.b} If there are a matrix \(\gamma\in\Cggq\) and a \tnnH{} measure \(\mu \in\Mggqrhl \) such that \(F\) can be represented via \eqref{Nr.FM} for each \(z\in \Cs \), then \(F\) belongs to the class \(\SFqlep \).
 \eenui
\ethm
\bproof
 Denote \(f\defeq  F\circ \Rstr_{\C\setminus [0,+\infty)} T_\lep \).

 \eqref{368.a} Let \(F \in\SFqlep \). According to \rrem{366}, the function \(f\) belongs to \(\SFuu{q}{0}\). In view of \rprop{HCT}, then \(\Rstr_{\ohe } f\) belongs to \(\RFq\), the Nevanlinna parametrization \((\Npa,\Npb,\Npnu)\) of \(\Rstr_{\ohe } f\) fulfills \eqref{Nr.Neu1}, and, for each \(w\in\C\setminus [0,+\infty)\), we have \(f (w)= \Npa + \int_{[0,+\infty)} \rk{1+xw}/\rk{x-w}\Npnua{\dif x}\).
 Because of \eqref{Nr.Neu1}, the integral \(\int_{[0,+\infty)}x\Npnua{\dif x}\) exists and the mapping \(\hat\mu\colon\Bsa{[0,+\infty)}\to\Cggq\) given by \(\hat\mu(B)\defeq\int_B(1+x^2)/(1+x)\Npnua{\dif x}\) is well defined and belongs to \(\Mggqa{[0,+\infty)}\). Setting \(\gamma\defeq\Npa-\int_{[0,+\infty)}x\Npnua{\dif x}\), for each \(w\in\C\setminus[0,+\infty)\), we get
 \[
  \begin{split}
   f(w)
   &=\Npa + \int_{[0,+\infty)}\rk*{\frac{1+x^2}{x-w}-x}\Npnua{\dif x}\\
   &=\Npa-\int_{[0,+\infty)}x\Npnua{\dif x}+\int_{[0,+\infty)}\rk*{\frac{1+x}{x-w}\cdot\frac{1+x^2}{1+x}}\Npnua{\dif x}\\
   &=\gamma+\int_{[0,+\infty)}\frac{1+x}{x-w}\hat\mu (\dif x).
  \end{split}
 \]
 Obviously, \(\mu\defeq(\Rstr_{[0,+\infty)}T_\lep )(\hat\mu)\) belongs to \(\Mggqrhl \). For each \(z\in\Cs \), \rprop{P1100} yields
 \[
  \begin{split}
   F (z)
   &= f(z-\lep )
   = \gamma + \int_{[0,+\infty)} \frac{1+x+\lep -\lep }{x + \lep -z}\hat\mu (\dif x)\\
   & = \gamma +  \int_{[0,+\infty)} (\huu{\lep}{z} \circ T_\lep  )(x)\hat\mu (\dif x)
   = \gamma +  \int_{T_\lep  ([0,+\infty))} \huu{\lep}{z} \dif\mu
  \end{split}
 \]
 and, hence, \eqref{Nr.FM} holds for each \(z\in\Cs \).

 Now we assume that \(\gamma\) is an arbitrary complex \tqqa{matrix} and that \(\mu\) is an arbitrary \tnnH{} measure belonging to \(\Mggqrhl \) such that \eqref{Nr.FM} holds for each \(z\in\Cs \). Observe that \(\lim_{n\to+\infty} h_{\lep ,\lep -1-n} (t) = 0\) is valid for each \(t\in \rhl \). Moreover, for every choice of \(n\in\NO\) and \(t\in \rhl \), one can easily check that the estimation \(\abs{h_{\lep ,\lep -1-n} (t)}\leq1\) holds. Consequently, a matrix generalization of Lebesgue's dominated convergence theorem (see \rprop{P1102}) yields
 \[
  \lim_{n\to+\infty} \int_{\rhl} \frac{1+t-\lep }{t-(\lep -1-n)} \mu (\dif t) 
  = 0.
 \]
 From \eqref{Nr.FM} we obtain then
 \[
  \gamma
  = \gamma + \lim_{n\to+\infty} \int_{\rhl} \frac{1 + t-\lep }{t-(\lep -1-n)} \mu (\dif t) 
  = \lim_{n\to+\infty} F(\lep -1-n).
 \]
 From \rthmp{T1542}{T1542.a} we conclude then \(\gamma=\CF\) and thus \(\gamma\in\Cggq\) follows. The mapping \(\tilde T_{-\lep }\colon\rhl \to[0,+\infty)\) defined by \(\tilde T_{-\lep }(t)\defeq t-\lep \) is obviously bijective and \(\Bsa{\rhl }\)\nobreakdash-\(\Bsa{[0,+\infty)}\)\nobreakdash-measurable. Further, the mapping \(\rho \colon\BsaR  \to\Cggq\) given by \(\rho (B)  \defeq\int_{B\cap[0,+\infty)}\rk{1+x}/\rk{1+x^2}\ek{\tilde T_{-\lep } (\mu)}(\dif x)\)
 is well defined, belongs to \(\MggqR \), and satisfies \(\rho ((-\infty, 0)) = \Oqq\). Furthermore, the integral \(\int_{[0,+\infty)}x\rho(\dif x)\) exists. For each \(w\in\C\setminus[0,+\infty)\), using \eqref{Nr.FM}, the relation \(\tilde T_{-\lep } (\rhl ) = [0,+\infty)\) and \rprop{P1100} provide us
 \beql{Nr.GLB}
  \begin{split}
   f(w)
   &= F(w+\lep )
   =\gamma +\int_{\rhl} \frac{1+ t-\lep }{t-(w+\lep )} \mu (\dif t)\\
   & =\gamma + \int_{\rhl} \frac{1+\tilde T_{-\lep }(t)}{\tilde T_{-\lep } (t) -w} \mu (\dif t)
   =\gamma + \int_{[0,+\infty)} \frac{1+x}{x-w}\ek*{\tilde T_{-\lep } (\mu)}(\dif x)\\
   &=\gamma + \int_{[0,+\infty)}\rk*{\frac{1+x^2}{x-w}\cdot\frac{1+x}{1+x^2}}\ek*{\tilde T_{-\lep } (\mu)}(\dif x)\\
   &=\gamma + \int_{[0,+\infty)} \frac{1+x^2}{x-w}\rho(\dif x)
   = \gamma + \int_{[0,+\infty)}  \rk*{ x +\frac{1+xw}{x-w}} \rho (\dif x)\\
   &= \gamma + \int_{[0,+\infty)}x\rho(\dif x)+\int_{[0,+\infty)}\frac{1+xw}{x-w} \rho (\dif x)
 \end{split}
 \eeq
 and, consequently,
 \[
  \Rstr_{\ohe } f(w)
  = \gamma + \int_{[0,+\infty)}x\rho(\dif x)+w\cdot\Oqq+\int_\R\frac{1+xw}{x-w} \rho (\dif x)
 \]
 for each \(w\in\ohe\). Since \(\gamma \) is \tnnH{} we see that \(A^\diamond\defeq\gamma +\int_{[0,+\infty)}x\rho(\dif x)\) belongs to \(\Cggq\). Thus, \({(A^\diamond, \Oqq,\rho)}\) coincides with the Nevanlinna parametrization \((\Npa,\Npb,\Npnu)\) of the Herglotz-Nevanlinna function \(\Rstr_{\ohe } f\). In particular, \(\rho\) is exactly the (unique) Nevanlinna measure \(\Npnu\) of \(\Rstr_{\ohe } f\). For each \(B\in\Bsarhl\), we have \(\tilde T_{-\lep } (B)\in \Bsa{[0,+\infty)}\) and hence
 \[
  \begin{split}
   \mu (B)
   &=\ek*{\tilde T_{-\lep } (\mu)}\rk*{\tilde T_{-\lep } (B)}
   =\int_{\tilde T_{-\lep } (B)}\rk*{\frac{1+x^2}{1+x}\cdot\frac{1+x}{1+x^2}}\ek*{\tilde T_{-\lep } (\mu)}(\dif x)\\
   &=\int_{\tilde T_{-\lep } (B)}\frac{1+x^2}{1+x}\rho(\dif x)
   =\int_{\tilde T_{-\lep } (B)}\frac{1+x^2}{1+x}\Npnua{\dif x}.
  \end{split}
 \]
 In particular, \(\mu\) is uniquely determined.

 \eqref{368.b} Let \(\gamma\in\Cggq\) and \(\mu\in\Mggqrhl \) be such that \(F\) can be represented via \eqref{Nr.FM} for each \(z\in\Cs \). Then the mapping \(\hat{\rho}\colon\Bsa{[0,+\infty)}\to\Cggq\) given by \(\hat\rho(B) \defeq\int_B(1+x)/(1+x^2)[\tilde T_{-\lep } (\mu)](\dif x)\) is well defined and belongs to \(\Mggqa{[0,+\infty)}\). Furthermore, \(f\) satisfies \eqref{Nr.GLB} with \(\hat\rho\) instead of \(\rho\) for each \(w\in\C\setminus [0,+\infty)\). Hence, using a result on holomorphic dependence of an integral on a complex parameter (see, \eg{}~\zitaa{MR2257838}{Ch.~IV, \S5, \cSatz{5.8}}), we conclude that \(f\) is a matrix-valued function which is holomorphic in \(\C\setminus [0,+\infty)\). Because of \(F (z) = f(z-\lep )\) for each \(z\in\Cs \), we obtain then that \(F\) is holomorphic in \(\Cs \). From \rlem{N6N} we get \(\im F(w)\in\Cggq\) for each \(w\in\ohe\) and \(F(x)\in\Cggq\) for each \(x\in\lhl \). Thus, \(F\) belongs to \(\SFqlep\).
\end{proof}

\breml{JB}
 In the following, if \(F\in \SFqlep \) is given, then we will write \((\gammau{F},\muu{F})\)\index{g@\(\gammau{F}\)}\index{m@\(\muu{F}\)} for the unique pair \((\gamma,\mu)\in\Cggq\times\Mggqrhl \) which realizes the integral representation \eqref{Nr.FM} for each \(z\in\Cs \).
\erem

\bexal{E0754}
 Let \(\lep\in\R\) and let \(A,B\in\Cggq\). Let \(F\in\Cs\to\Cqq\) be the function from \(\SFqlep\) which is defined in \rexa{E1256}. Then \(\gammaF=A\) and \(\muF=\Kron{\lep}B\) where \(\Kron{\lep}\) denotes the Dirac measure on \((\rhl,\Bsarhl)\) with unit mass at \(\lep\).
\eexa

\bexal{E0945}
 Let \(\lep\in\R\), let \(\gamma\in\Cggq\), and let \(F\colon\Cs\to\Cqq\) be defined for each \(z\in\Cs \) by \(F(z)\defeq\gamma\). In view of \rthm{368}, then \(F\in\SFqlep\), \(\gammau{F}=\gamma\), and \(\muu{F}\) is the zero measure belonging to \(\Mggqrhl \).
\eexa

Now we state some observations on the arithmetic of the class \(\SFqlep\).
\breml{R1018}
 If \(F\in\SFqlep\) then \(F^\tra\in\SFqlep\) and \((\gammau{F^\tra},\muu{F^\tra})=(\gammaF^\tra,\muF^\tra)\).
\erem

\breml{R370}
 Let \(\lep \in\R\), let \(n\in\N\), and let \((q_k)_{k=1}^n\) be a sequence of positive integers. For each \(k\in\mn{1}{n}\), let \(F_k\in\SFuu{q_k}{\lep }\). Then \(F \defeq\diag\matauuo{F_k}{k}{1}{n}\) belongs to \(\SFuu{\sum_{k=1}^n q_k}{\lep }\) and \((\gammau{F},\muu{F})=\rk{\diag\matauuo{\gammau{F_k}}{k}{1}{n},\diag\matauuo{\muu{F_k}}{k}{1}{n}}\). Moreover, if \(A_k\in\Coo{q_k}{q}\) for each \(k\in\mn{1}{n}\), then \(G \defeq  \sum_{k=1}^n A_k^\ad F_k A_k\) belongs to \(\SFqlep \) and \((\gammau{G},\muu{G})=\rk{\sum_{k=1}^n A_k^\ad \gammau{F_k} A_k, \sum_{k=1}^n A_k^\ad \muu{F_k} A_k}\).
\erem

\breml{R1612}
 Let \(\lep \in\R\) and let \(F\in \SFqlep \). For each matrix \(A\in\Cqq \) for which the matrix \(\gammau{F} + A\) is \tnnH{}, from \rthm{368} one can easily see that the function \(G\defeq F+A\) belongs to \(\SFqlep \) and that \((\gammau{G},\muu{G})=(\gammau{F}+ A, \muu{F})\).
\erem

\bpropl{R367Ad}
 Let \(\lep\in\R\) and let \(F\in\SFqlep\). Then \(\lim_{y\to+\infty}F(\iu y)=\gammau{F}\).
\eprop
\bproof
 Let \(\seq{y_n}{n}{1}{\infty}\) be a sequence from \([1,+\infty)\) such that \(\lim_{n\to+\infty}y_n=+\infty\). Obviously, in view of \rrem{R1618}, we have \(\lim_{n\to+\infty}\huua{\lep }{\iu y_n}{t}=0\) for each \(t\in\rhl \). Furthermore, for each \(t\in\rhl \), we get \(\abs{\huua{\lep }{\iu y_n}{t}}\leq3+\abs{\lep }\). By virtue of \rprop{P1102}, we obtain then \(\lim_{n\to+\infty}\int_{\rhl} \huu{\lep }{\iu y_n}\dif\muu{F}=\Oqq\). Application of the integral representation stated in \rthmp{368}{368.a} completes the proof.
\eproof

\bcorl{C1024}
 Let \(\lep\in\R\). Then \(\SdFqlep=\setaa{F\in\SFqlep}{\gammaF=\Oqq}\).
\ecor
\bproof
 Combine \rprop{R367Ad} and \rrem{JB}.
\eproof

If \(\cX\) is a \tne{} subset of \(\Cq\), then we will use \(\cX^\oc\)\index{\(\cX^\oc\)} to denote the orthogonal space of \(\cX\). For each \(A\in\Cpq\), let \(\nul{A}\)\index{n@\(\nul{A}\)} be the null space of \(A\) and let \(\ran{A}\)\index{r@\(\ran{A}\)} be the column space of \(A\).

Recall that a complex \tqqa{matrix} \(A\) is called an \emph{EP~matrix} if \(\ran{A^\ad}=\ran{A}\). The class \(\CEPq\)\index{c@\(\CEPq\)} of these complex \tqqa{matrices} was introduced by Schwerdtfeger~\zita{MR0038923}. For a comprehensive treatment of the class \(\CEPq\) against to the background of this paper, we refer the reader to~\zitaa{MR3014198}{\capp{A}}. If \(G\in\RFq\) then it was proved in~\zitaa{MR3014198}{\clem{9.1}} that, for each \(w\in\ohe\), the matrix \(G(w)\) belongs to \(\CEPq\). \rPart{L374.a} of the following proposition shows that an analogous result is true for functions belonging to the class \(\SFqlep\). Furthermore, the following proposition contains also extensions to the class \(\SFqlep\) of former results (see~\zitaa{MR3014198}{\cthm{9.4}},~\zitaa{MR2988005}{\cprop{3.7}}) concerning the class \(\RFq\).

\bpropl{L374}
 Let \(\lep \in\R\) and let \(F\in\SFqlep \). Then:
 \benui
  \il{L374.aa} Let \(z\in\Cs\). Then \(\ko z\in\Cs\) and \(F^\ad(z)=F(\ko z)\).
  \il{L374.a} For each \(z\in \Cs \), the equations
  \begin{align}\label{Nr.NM}
   \Nul{F(z)}&=\nul{\gammau{F}}\cap\Nul{\muuA{F}{\rhl }}\\
  \shortintertext{and}\label{L374.1}
   \Ran{F(z)}&=\ran{\gammau{F}}+\Ran{\muuA{F}{\rhl }}
  \end{align}
  are valid. In particular, \(\nul{F(z)}=\nul{F^\ad(z)}=\ek{\ran{F(z)}}^\oc\) and \(\ran{F(z)}=\ran{F^\ad(z)}=\ek{\nul{F(z)}}^\oc\) for each \(z\in\Cs  \).
  \il{L374.b} Let \(r\in\NO\). Then the following statements are equivalent:
  \baeqii{0}
   \item For each \(z\in\Cs \), the equation \(\rank  F(z) = r\) holds.
   \item There is some \(z_0\in\Cs \) such that \(\rank  F(z_0) = r\).
   \item \(\dim[ \ran{\gammau{F}} + \ran{\muu{F} ( \rhl  }]= r\).
  \eaeqii
  \il{L374.d} If \(\det\gammaF\neq0\) or \(\det\ek{\muFa{\rhl}}\neq0\) then \(\det\ek{F(w)}\neq0\) for all \(w\in\Cs\).
 \eenui
\eprop
\bproof
 In view of \rrem{JB}, we have
 \begin{align}\label{Nr.BM1}
  \gammau{F}&\in\Cggq&
 &\text{and}&
  \muu{F}&\in\MggqRhl.
 \end{align}
 Taking into account \rrem{R1618}, for all \(z\in\Cs\), we get then that \(\huu{\lep}{z}\) belongs to  \(\Loaaaa{1}{\rhl }{\Bsarhl }{\muu{F}}{\C}\).
 
 \eqref{L374.aa} Let \(z\in\Cs \). From \rrem{JB} and \eqref{Nr.BM1} we obtain
 \[
  F^\ad(z)
  = \gammau{F}^\ad + \ek*{\int_{\rhl} \frac{1+ t-\lep }{t-z} \muua{F}{\dif t}}^\ad
  = \gammau{F} + \int_{\rhl} \frac{1+ t-\lep }{t- \ko z } \muua{F}{\dif t}
  = F(\ko z ).
 \]
 
 \eqref{L374.a} Let \(z\in\Cs \). For each \(t\in \rhl \), we get then
 \begin{align}\label{Nr.BM7-1}
  \re  \huua{\lep}{z}{t}&= \frac{(t- \re  z)(1 + t-\lep )}{\abs{t-z}^2}\\
 \shortintertext{and}\label{Nr.BM7-2}
 \im  \huua{\lep}{z}{t}&= \frac{(\im z)(1 + t-\lep )}{\abs{t-z}^2}.
 \end{align}
 Since \(\huu{\lep}{z}\) belongs to  \(\Loaaaa{1}{\rhl }{\Bsarhl }{\muu{F}}{\C}\), from \rlemp{62B}{62B.c} we see that
 \beql{Nr.BM4C}
  \Nul{\muuA{F}{\rhl}}
  \subseteq\Nul{\int_{\rhl} \huu{\lep}{z} \dif\muu{F}}
 \eeq
 holds true. Now we consider an arbitrary \(u\in\nul{F(z)}\). In view of the definition of the pair \((\gammau{F},\muu{F})\) (see \rrem{JB}), we have
 \beql{Nr.BM9}
  u^\ad \gammau{F} u + \int_{\rhl} \huu{\lep}{z}\dif(u^\ad \muu{F} u).
  = u^\ad F (z) u
  = 0.
 \eeq
 Consequently, because of \eqref{Nr.BM9}, \eqref{Nr.BM1}, and \eqref{Nr.BM7-2}, then
 \beql{Nr.BM20}
  0 
  = u^\ad \gammau{F} u + \int_{\rhl} \re \huu{\lep}{z}\dif(u^\ad \muu{F} u)
  \geq  \int_{\rhl} \re  \huu{\lep}{z}\dif(u^\ad \muu{F} u)
 \eeq
 and
 \beql{Nr.BM12}
  0
  = \im  \ek*{\int_{\rhl} \huu{\lep}{z}\dif(u^\ad \muu{F} u)}
  = \int_{\rhl}  \frac{(\im  z)(1+ t-\lep )}{\abs{t-z}^2} (u^\ad\muu{F} u) (\dif t)
 \eeq
 follow. In the case \(\im  z\neq 0\), from \eqref{Nr.BM12} and \eqref{Nr.BM1} we get \(\rk{u^\ad \muu{F} u}(\rhl )=0\). If \(z\in \lhl \), then from \eqref{Nr.BM7-1} we see that \(\re  \huua{\lep}{z}{t} \in (0,+\infty)\) holds for each \(t\in \rhl \), and, by virtue of \eqref{Nr.BM20} and \eqref{Nr.BM1}, we obtain \(\rk{u^\ad \muu{F} u} (\rhl ) = 0\). Thus \(\rk{u^\ad \muu{F} u} (\rhl ) = 0\) is proved in each case, which, in view of \eqref{Nr.BM1}, implies \(u\in\nul{\muua{F}{\rhl }}\). Taking into account a standard argument of the integration theory of \tnnH{} measures and \eqref{Nr.BM4C}, we conclude that
 \[
  \int_{\rhl} \huu{\lep}{z}\dif(u^\ad\muu{F} u) 
  = u^\ad \rk*{ \int_{\rhl} \huu{\lep}{z}\dif\muu{F}} u
  =u^\ad\cdot\Ouu{q}{1}
  = 0.
 \]
 Consequently, from \eqref{Nr.BM9} we infer \(u^\ad \gammau{F} u = 0\). Thus, \eqref{Nr.BM1} shows that \(u\) belongs to \(\nul{\gammau{F}}\). Hence,
 \beql{Nr.BM28}
  \Nul{F(z)}
  \subseteq \nul{\gammau{F}} \cap \Nul{\muuA{F}{\rhl}}
 \eeq
 is valid. Now we are going to check that 
 \beql{Nr.BM33}
  \nul{\gammau{F}} \cap \Nul{\muuA{F}{\rhl}}
  \subseteq \Nul{F(z)}
 \eeq
 holds. For this reason, we consider an arbitrary \(u\in\nul{\gammau{F}}\cap \nul{\muu{F}(\rhl )}\). From \eqref{Nr.BM4C} we get then  \(F(z) u = \gammau{F} u + \rk{ \int_{\rhl} \huu{\lep}{z} \dif\muu{F}} u =\Ouu{q}{1}\)
 and therefore \(u\in\nul{F(z)}\). Hence \eqref{Nr.BM33} is verified. From \eqref{Nr.BM28} and \eqref{Nr.BM33} then 
 \eqref{Nr.NM} follows. Keeping in mind~\eqref{L374.aa}, \eqref{Nr.NM} for \(\ko z\) instead of \(z\), and \eqref{Nr.BM1}, standard arguments of functional analysis yield then
 \[
  \begin{split}
   \Ran{F (z)}
   & =\ek*{\Nul{F(z)^\ad}}^\oc 
   =\ek*{\Nul{F (\ko z )}}^\oc\\
   &= \ek*{\nul{\gammau{F}} \cap \Nul{\muuA{F}{\rhl}}}^\oc
   = \ek*{\ran{\gammau{F}}^\oc  \cap \Ran{\muuA{F}{\rhl}}^\oc  }^\oc\\
   &= \rk*{\ek*{\Span\rk*{\ran{\gammau{F}}\cup \Ran{\muuA{F}{\rhl}}}}^\oc  }^\oc
   = \ran{\gammau{F}}+\Ran{\muuA{F}{\rhl}}.
  \end{split}
 \]
 Thus, \eqref{L374.1} is proved. Using \eqref{Nr.NM} for \(z\) and for \(\ko z\) instead of \(z\), from~\eqref{L374.aa} we obtain \(\nul{F(z)}=\nul{F^\ad(z)}=\ek{\ran{F(z)}}^\oc\). Similarly, \(\ran{F(z)}=\ran{F^\ad(z)}=\ek{\nul{F(z)}}^\oc\) follows from \eqref{L374.1} and~\eqref{L374.aa}.
 
 \eqref{L374.b}--\eqref{L374.d} These are immediate consequences of~\eqref{L374.a}.
\end{proof}

\rprop{L374} yields a generalization of a result due to Kats and Krein~\zitaa{KK74}{\ccor{5.1}}:
\bcorl{C1350}
 Let \(\lep \in\R\), let \(F\in\SFqlep \), and let \(z_0\in\Cs\). Then \(F(z_0)=\Oqq\) if and only if \(F(z)=\Oqq\) for all \(z\in\Cs\).
\ecor
\bproof
 This is an immediate consequence of \rpropp{L374}{L374.b}.
\eproof

\bcorl{C1127}
 Let \(\lep \in\R\), let \(F\in\SFqlep \), and let \(\lambda\in\R\) be such that the matrix \(\gammaF-\lambda\Iq\) is \tnnH{}. Then
 \begin{align}\label{C1127.B1}
  \Ran{F(z)-\lambda\Iq}&=\Ran{F(w)-\lambda\Iq}&
  &\text{and}&
  \Nul{F(z)-\lambda\Iq}&=\Nul{F(w)-\lambda\Iq}
 \end{align}
 for all \(z,w\in\Cs\). In particular, if \(\lambda\leq0\), then \(\lambda\) is an eigenvalue of the matrix \(F(z_0)\) for some \(z_0\in\Cs\) if and only if \(\lambda\) is an eigenvalue of the matrix \(F(z)\) for all \(z\in\Cs\). In this case, the eigenspaces \(\nul{F(z)-\lambda\Iq}\) are independent of \(z\in\Cs\).
\ecor
\bproof
 In view of \rrem{JB} and \rthm{368}, we conclude that the function \(G\colon\Cs\to\Cqq\) defined by \(G(z)\defeq F(z)-\lambda\Iq\) belongs to \(\SFqlep\). The application of \rpropp{L374}{L374.a} to the function \(G\) yields then \eqref{C1127.B1}. Since the matrix \(\gammaF\) is \tnnH{}, we have \(\gammaF-\lambda\Iq\in\Cggq\) if \(\lambda\leq0\). Thus, the remaining assertions are an immediate consequence of \eqref{C1127.B1}.
\eproof

At the end of this section we add a useful technical result.
\bleml{L1716}
 Let \(\lep \in\R\), let \(A\in\Cpq\), and let \(F\in\SFqlep\). Then the following statements are equivalent:
 \baeqi{0}
  \il{L1716.i} \(\nul{A}\subseteq\nul{F(z)}\) for all \(z\in\Cs  \).
  \il{L1716.i1} There is a \(z_0\in\C\setminus\rhl \) such that \(\nul{A}\subseteq\nul{F(z_0)}\).
  \il{L1716.ii} \(\nul{A}\subseteq\nul{\gammau{F}}\cap\nul{\muua{F}{\rhl }}\).
  \il{L1716.iii} \(FA^\mpi A=F\).
  \il{L1716.iv} \(\ek{\nul{A}}^\oc\supseteq\ran{F(z)}\) for all \(z\in\Cs  \).
  \il{L1716.iv1} There is a \(z_0\in\C\setminus\rhl \) such that \(\ek{\nul{A}}^\oc\supseteq\ran{F(z_0)}\).
  \il{L1716.v} \(\ek{\nul{A}}^\oc\supseteq\ran{\gammau{F}}+\ran{\muua{F}{\rhl }}\).
  \il{L1716.vi} \(A^\mpi AF=F\).
 \eaeqi
\elem
\bproof
 \impl{L1716.i}{L1716.i1} and \impl{L1716.iv}{L1716.iv1}: These implications hold true obviously.
 
 \aequ{L1716.i}{L1716.ii} and \impl{L1716.i1}{L1716.ii}: Use equation \eqref{Nr.NM} in \rpropp{L374}{L374.a}.
 
 \baeq{L1716.i}{L1716.iii}
  This equivalence follows from a well-known result for the Moore-Penrose inverse of complex matrices.
 \eaeq
 
 \baeq{L1716.i}{L1716.iv}
  Because of \rpropp{L374}{L374.a}, we have \(\nul{F(z)}=\ran{F(z)}^\oc\) for all \(z\in\Cs\). Hence,~\rstat{L1716.i} and~\rstat{L1716.iv} are equivalent.
 \eaeq
 
 \aequ{L1716.iv}{L1716.v} and \impl{L1716.iv1}{L1716.v}: Use equation~\eqref{L374.1} in \rpropp{L374}{L374.a}.
 
 \baeq{L1716.iv}{L1716.vi}
  Use \(\nul{A}^\oc=\ran{A^\ad}\) and \(A^\mpi A\ran{A^\ad}=\ran{A^\ad}\).
 \eaeq
\eproof

Now we apply the preceding results to the subclass \(\SdFqlep\) of \(\SFqlep\) (see \rnota{D1026}).

\bexal{E0949}
 Let \(\lep\in\R\) and let \(F\colon\Cs\to\Cqq\) be defined by \(F(z)\defeq\Oqq\). In view of \rexa{E0945} and \rcor{C1024}, one can easily see then that \(F\) belongs to \(\SdFqlep\) and that \(\muu{F}\) is the zero measure belonging to \(\Mggqrhl \).
\eexa

\breml{R1717}
 Let \(\lep\in\R\), let \(n\in\N\), and let \((q_k)_{k=1}^n\) be a sequence of positive integers. For each \(k\in\mn{1}{n}\), let \(F_k\in\SdFuu{q_k}{\alpha}\) and let \(A_k\in\Coo{q_k}{q}\). In view of \rcor{C1024} and \rrem{R370}, then:
 \benui
  \il{R1717.a} The function \(G\defeq\sum_{k=1}^n A_k^\ad F_k A_k\) belongs to \(\SdFqlep \) and \(\muu{G}=\sum_{k=1}^n A_k^\ad\muu{F_k} A_k\).
  \item The function \(F\defeq\diag\matauuo{F_k}{k}{1}{n}\) belongs to \(\SdFuu{\sum_{k=1}^n q_k}{\alpha}\) and \(\muu{F}=\diag\matauuo{\muu{F_k}}{k}{1}{n}\).
 \eenui
\erem

\brem
 Let \(\lep\in\R\) and let \(F\in\SdFqlep \). In view of \rcor{C1024} and \rpropp{L374}{L374.a}, then \(\nul{F(z)}=\nul{\muua{F}{\rhl }}\) and \(\ran{F(z)}=\ran{\muua{F}{\rhl }}\) for all \(z\in\Cs  \).
\erem

\section{Characterizations of the class $\SFqlep$}\label{S1038}
In this section, we give several characterizations of the class \(\SFqlep\).

\bleml{CHR}
 Let \(\lep\in\R\), let \(F\colon\Cs\to\Cqq\) be holomorphic, and let \(F^\mulz\colon\Cs\to\Cqq\) be defined by
 \beql{NDF}
  F^\mulz(z)
  \defeq(z-\lep)F(z).
 \eeq
 Suppose that \(\Rstr_\ohe F\) and \(\Rstr_\ohe F^\mulz\) both belong to \(\RFq\). Then \(\rk{\re F}(\lhl)\subseteq\Cggq\).
\elem
\bproof
 We consider an arbitrary \(x\in\lhl\). For each \(n\in\N\), we have then
 \beql{RI}
  \re F\rk*{x+\frac{\iu}{n}}
  =n\im F^\mulz\rk*{x+\frac{\iu}{n}}+n(\lep-x)\im F\rk*{x+\frac{\iu}{n}}.
 \eeq
 For each \(n\in\N\), \(\Rstr_\ohe F^\mulz\in\RFq\) implies \(n\im F^\mulz\rk{x+\iu/n}\in\Cggq\), whereas \(\Rstr_\ohe F\in\RFq\) yields \(n(\lep-x)\im F\rk{x+\iu/n}\in\Cggq\). Thus, \eqref{RI} provides us \(\re F\rk{x+\iu/n}\in\Cggq\) for each \(n\in\N\). Since \(F\) is continuous, we get then \(\re F(x)=\lim_{n\to\infty}\re F\rk{x+\iu/n}\). In particular, the matrix \(\re F(x)\) is \tnnH{}.
\eproof

In order to give further characterizations of the class \(\SFqlep\), we state the following technical result. The proof of which uses an idea which is originated in~\zitaa{MR0458081}{\cthm{A.5}}.
\bleml{NID}
 Let \(\lep\in\R\) and let \(F\in\SFqlep\). Then \(F^\mulz\colon\Cs\to\Cqq\) defined by \eqref{NDF} is holomorphic and fulfills
 \beql{M121N}
  \im F^\mulz(z)
  =(\im z)\ek*{\gammau{F}+\int_{\rhl}\frac{(1+t-\lep)(t-\lep)}{\abs{t-z}^2}\muua{F}{\dif t}}
 \eeq
 for each \(z\in\Cs\). Furthermore, \(\frac{1}{\im z}\im F^\mulz(z)\in\Cggq\) for each \(z\in\C\setminus\R\) and \(F^\mulz(x)\in\CHq\) for each \(x\in\lhl\).
\elem
\bproof
 Since \(F\) is holomorphic, the matrix-valued function \(F^\mulz\) is holomorphic as well. In view of \rrem{JB}, using a well-known result on integrals with respect to \tnnH{} measures, we have
 \beql{Nr.FMN}
  \ek{F(z)}^\ad
  =\gammau{F}+\int_{\rhl}\frac{1+t-\lep}{t-\ko z}\muua{F}{\dif t}
 \eeq
 for each \(z\in\Cs\). Thus, from \rrem{JB} and \eqref{Nr.FMN} we get
 \beql{Nr.FMM}
  \begin{split}
   2\iu\im F^\mulz(z)
   &=(z-\lep)F(z)-(\ko z-\lep)\ek*{F(z)}^\ad\\
   &=(z-\ko z)\gammau{F}+\int_{\rhl}(1+t-\lep)\rk*{\frac{z-\lep}{t-z}-\frac{\ko z-\lep}{t-\ko z}}\muua{F}{\dif t}
  \end{split}
 \eeq
 for each \(z\in\Cs\). Since \((z-\lep)/(t-z)-(\ko z-\lep)/(t-\ko z)=2\iu(\im z)(t-\lep)/\abs{t-z}^2\) holds for every choice of \(z\in\Cs\) and \(t\in\rhl\), from \eqref{Nr.FMM} it follows \eqref{M121N} for each \(z\in\Cs\). Since \((1+t-\lep)(t-\lep)/\abs{t-z}^2\in[0,+\infty)\) holds true for each \(z\in\Cs\) and each \(t\in\rhl\), from \(\gammau{F}\in\Cggq\) and \eqref{M121N} we get \(\frac{1}{\im z}\im F^\mulz(z)\in\Cggq\) for each \(z\in\C\setminus\R\). In view of \rlem{NN}, then \(F^\mulz(x)=\re F^\mulz(x)\) and hence \(F^\mulz(x)\in\CHq\) for each \(x\in\lhl\).
\eproof

\bpropl{M1.47N}
 Let \(\lep\in\R\), let \(F\colon\Cs\to\Cqq\) be a matrix-valued function, and let \(F^\mulz\colon\Cs\to\Cqq\) be defined by \eqref{NDF}. Then \(F\) belongs to \(\SFqlep\) if and only if the following two conditions hold true:
 \bAeqi{0}
  \il{M1.47N.I} \(F\) is holomorphic in \(\Cs\).
  \il{M1.47N.II} The matrix-valued functions \(\Rstr_\ohe F\) and \(\Rstr_\ohe F^\mulz\) both belong to \(\RFq\).
 \eAeqi
\eprop
\bproof
 First suppose that \(F\) belongs to \(\SFqlep\). Then~\rstat{M1.47N.I} and \(\Rstr_\ohe F\in\RFq\) follow from the definition of the class \(\SFqlep\). Furthermore, \rlem{NID} provides us \(\Rstr_\ohe F^\mulz\in\RFq\).
 
 Conversely, now suppose that~\rstat{M1.47N.I} and~\rstat{M1.47N.II} hold true. Because of the definition of the classes \(\RFq\) and \(\SFqlep\), it then remains to prove that \(F(\lhl)\subseteq\Cggq\). We consider an arbitrary \(x\in\lhl\). First we show that \(\im F(x)=\Oqq\). Because of~\rstat{M1.47N.II}, for each \(n\in\N\), the matrices \(\im F(x+\iu/n)\) and \(\im F^\mulz(x+\iu/n)\) are \tnnH{}. Thus, the matrices \(\im F(x)\) and \(\im F^\mulz(x)\) are (as limits of the sequences \(\seq{\im F(x+\iu/n)}{n}{1}{\infty}\) and \(\seq{\im F^\mulz(x+\iu/n)}{n}{1}{\infty}\), respectively) \tnnH{} as well. Since \eqref{NDF} implies \(\im F^\mulz(x)=(x-\lep)\im F(x)\), we get then \(-\im F(x)=\frac{1}{\lep-x}\im F^\mulz(x)\in\Cggq\), which together with \(\im F(x)\in\Cggq\) yields \(\im F(x)=\Oqq\). Hence, \(\re F(x)=F(x)\). Because of~\rstat{M1.47N.I},~\rstat{M1.47N.II}, and \rlem{CHR}, we have \(\re F(x)\in\Cggq\). Thus, \(F(\lhl)\subseteq\Cggq\).
\eproof

Let \(\lhelep\defeq\setaa{z\in\C}{\re z\in(-\infty,\alpha)}\)\index{c@\(\lhelep\)}.

\bpropl{214-1}
 Let \(\lep \in\R\) and let \(F\colon\Cs \to\Cqq\) be a matrix-valued function. Then \(F\) belongs to \(\SFqlep\) if and only if the following four conditions are fulfilled:
 \bAeqi{0}
  \il{214-1.I} \(F\) is holomorphic in \(\Cs \).
  \il{214-1.II} For each \(z\in\ohe\), the matrix \(\im F(z)\) is \tnnH{}.
  \il{214-1.III} For each \(z\in\uhe\), the matrix \(-\im F(z)\) is \tnnH{}.
  \il{214-1.IV} For each \(z\in\lhelep \), the matrix \(\re F(z)\) is \tnnH{}.
 \eAeqi
\eprop
\bproof
 First suppose that \(F\in\SFqlep\). By definition of the class \(\SFqlep\), conditions~\rstat{214-1.I} and~\rstat{214-1.II} are fulfilled. From \rthm{368} and \rlem{N6N} we obtain~\rstat{214-1.III} and~\rstat{214-1.IV}.
 
 Conversely,~\rstat{214-1.I}--\rstat{214-1.III} and \rlem{NN} imply \(\im F(x)=\Oqq\) for all \(x\in\lhl \). In view of~\rstat{214-1.IV}, we have then \(F(\lhl )\subseteq\Cggq\). Together with~\rstat{214-1.I} and~\rstat{214-1.II}, this yields \(F\in\SFqlep\).
\eproof

\section{The class $\SOFqlep$}\label{S1242}

In this section, we prove an important integral representation for functions which belong to the class \(\SOFqlep\). It can be considered as modified integral representation of the functions belonging to the class \(\ROFq\defeq\setaa{F\in\RFq}{\sup_{y\in[1,+\infty)}y\norm{F(\iu y)}<+\infty}\) (see~\zitaa{MR2222521}{\cthm{8.7}}). Observe that if \(\lep\in\R\) and if \(z\in\Cs \), then in view of \rlemp{63A}{63A.a}, for each \(\sigma\in\Mggqrhl \), the integral \(\int_{\rhl}1/(t-z)\sigma(\dif t)\) exists.

\bthml{T1434}
 Let \(\lep\in\R\) and let \(F\colon\Cs  \to \Cqq \). Then:
 \benui
  \il{T1434.a} If \(F\in \SOFqlep  \), then there is a unique \tnnH{} measure \(\sigma\in\Mggqrhl \) such that
  \beql{T1434.1}
   F(z)
   =\int_{\rhl} \frac{1}{t-z} \sigma(\dif t)
  \eeq
  for each \(z\in\Cs \).
  \il{T1434.b} If there is a \tnnH{} measure \(\sigma\in\Mggqrhl \) such that \(F\) can be represented via \eqref{T1434.1} for each \(z\in \Cs \), then \(F\) belongs to the class \(\SOFqlep  \).
 \eenui
\ethm
\bproof
 We modify ideas of proofs of integral representations of similar classes of holomorphic functions (see~\zitaa{MR0458081}{Appendix}).

 \eqref{T1434.a} First suppose \(F\in\SOFqlep \). Then \(F\in\SFqlep\) and the function \(F\defeq\Rstr_\ohe F\) belongs to the class \(\ROFq\). From a matricial version of a well-known integral representation of functions belonging to \(\ROFu{1}\) (see, \eg{}~\zitaa{MR2222521}{\cthm{8.7}}) we know that there is a unique \(\mu\in\MggqR\) such that
 \beql{T1434.2}
  F(w)
  =\int_\R\frac{1}{t-w}\mu(\dif t)
 \eeq
 for all \(w\in\ohe\), namely the so-called spectral measure of \(F\), \ie{}, for each \(B\in\BsaR\), we have \(\mu(B)=\int_B(1+t^2)\Npnua{\dif t}\), where \(\Npnu\) denotes the Nevanlinna measure of \(F\). \rprop{214-1} shows that \(F\) is holomorphic in \(\Cs\) and that \(\im F(z)\in\Cggq\) for all \(z\in\ohe\) and \(-\im F(z)\in\Cggq\) for all \(z\in\uhe\). Hence, for each \(t\in\lhl \), we have \(F^\ad(t)=F(t)\). Applying the Stieltjes-Perron inversion formula (see, \eg{}~\zitaa{MR2222521}{\cthm{8.2}}), one can verify that \(\Npnua{\lhl}=\Oqq\). Hence \(\mu(\lhl )=\Oqq\). Consequently, formula \eqref{T1434.2} shows that \eqref{T1434.1} holds for each \(z\in\ohe\), where \(\sigma\defeq\Rstr_\Bsarhl\mu\). Since \(\rhl \) is a closed interval the function \(G\colon\Cs \to\Cqq\) defined by \(G(z)\defeq\int_{\rhl}1/(t-z)\sigma(\dif t)\) is holomorphic (see, \eg{}~\zitaa{MR2257838}{Ch.~IV, \S5, \cSatz{5.8}}). Because of \(F(w)=G(w)\) for each \(w\in\ohe\), we have \(F=G\). 
 If \(\sigma\) is an arbitrary measure belonging to \(\Mggqrhl \) such that \eqref{T1434.1} holds for each \(z\in\Cs \), then using standard arguments of measure theory and the uniqueness of the \tnnH{} \tqqa{measure} \(\Npnu\) in the integral representation \eqref{B1}, one gets necessary \(\sigma=\Rstr_\Bsarhl\mu\).
 
 \eqref{T1434.b} Let \(\sigma\in\Mggqrhl \) such that \eqref{T1434.1} holds for each \(z\in\Cs \). Then \(F\) is holomorphic (see, \eg{}~\zitaa{MR2257838}{Ch.~IV, \S5, \cSatz{5.8}}) and, for each \(z\in\C\setminus\R\), we have
 \[
  \frac{1}{\im  z} \im F(z)
  =  \int_{\rhl} \frac{1}{\im  z} \im  \frac{1}{t-z} \sigma(\dif t)
  = \int_{\rhl} \frac{1}{\abs{t-z}^2} \sigma(\dif t)
 \in\Cggq.
 \]
 and, for each \(z\) belonging to \(\lhelep \), moreover
 \[
  \re  F(z)
  =\int_{\rhl} \re  \frac{1}{t-z} \sigma(\dif t)
  =\int_{\rhl} \frac{t- \re  z}{\abs{t-z}^2}\sigma(\dif t)
  \in\Cggq.
 \]
 Thus, \(F\in\SFqlep\). From the definition of \(F\) and~\zitaa{MR2222521}{\cthm{8.7(b)}} we see that \(\Rstr_\ohe F\in\ROFq\), where \(\ROFq\) is the class of all \(H\in\RFq\) satisfying \(\sup_{y\in[1,+\infty)}y\normE{H(\iu y)}<+\infty\). Thus, \eqref{Nr.YS} is satisfied. Hence, \(F\in\SOFqlep \) holds.
\eproof

If \(\sigma\) is a measure belonging to \(\Mggqrhl \), then we will call the matrix-valued function \(F\colon\Cs \to\Cqq\) which is, for each \(z\in\Cs \), given by \eqref{T1434.1} the \emph{\tasto{\lep}{\(\sigma\)}}. If \(F\in\SOFqlep \), then the unique measure \(\sigma\) which belongs to \(\Mggqrhl \) and which fulfills \eqref{T1434.1} for each \(z\in\Cs \) is said to be the \emph{\tasmo{\lep}{\(F\)}} and will be denoted by \(\stm{F}\)\index{s@\(\stm{F}\)}.

Note that, in view of \rthm{T1434}, the matricial Stieltjes moment problem~\mproblem{\iraa{\alpha}}{m}{=} can be obviously reformulated in the language of \tast{\lep}{s} of \tnnH{} measures. We omit the details.

\breml{R0906}
 Let \(\lep\in\R\), let \(n\in\N\), and let \((q_k)_{k=1}^n\) be a sequence of positive integers. For each \(k\in\mn{1}{n}\), let \(F_k\in\SOFuu{q_k}{\alpha}\), and let \(\stm{F_k}\) be the \tasmo{\alpha}{\(F_k\)}. Then \(F\defeq\diag\matauuo{F_k}{k}{1}{n}\) belongs to \(\SOFuu{\sum_{k=1}^n q_k}{\alpha}\) and \(\diag\matauuo{\stm{F_k}}{k}{1}{n}\) is the \tasmo{\alpha}{\(F\)}. Moreover, if \(A_k\in\Coo{q_k}{q}\) for each \(k\in\mn{1}{n}\), then \(G \defeq  \sum_{k=1}^n A_k^\ad F_k A_k\) belongs to \(\SOFqlep  \) and \(\sum_{k=1}^n A_k^\ad \stm{F_k} A_k\) is the \tasmo{\alpha}{\(G\)}.
\erem

\bpropl{F394}
 Let \(\lep\in\R\), let \(F\in\SOFqlep \), and let \(\stm{F}\) be the \tasmo{\alpha}{\(F\)}. For each \(z\in\Cs \), then
 \begin{align}\label{F394.1}
  \Nul{F(z)}& =\Nul{\stmA{F}{\rhl}}&
 &\text{and}&
  \Ran{F(z)}&=\Ran{\stmA{F}{\rhl}}.
 \end{align}
 Furthermore,
 \beql{N14}
  \stmA{F}{\rhl}
  =-\iu\lim_{y\to+\infty}yF(\iu y).
 \eeq
 In particular, \(\rank F(z)=\rank\stma{F}{\rhl}\) holds true for each \(z\in\Cs\).
\eprop
\bproof
 Let \(z\in\Cs \). From \rthmp{T1434}{T1434.a} and \rlemp{63A}{63A.b} we obtain the second equation in \eqref{F394.1}. The first one is an immediate consequence of the second one and the equation \(F(z)=F^\ad(\ko z)\), which can be seen from \eqref{T1434.1}. Because of \rthmp{T1434}{T1434.a} and \rlemp{63A}{63A.c}, the equation \eqref{N14} holds true.
\eproof

\section{Moore-Penrose inverses of functions belonging to the class $\SFqlep$}\label{S1044}
We start with some further notation. If \(\cZ\) is a \tne{} subset of \(\C\) and if a matrix-valued function \(F\colon\cZ\to\Cpq\) is given, then let \(F^\mpi\colon\cZ\to\Cqp\) be defined by \(F^\mpi(z)\defeq\ek{F(z)}^\mpi\)\index{\(F^\mpi(z)\)}, where \(\ek{F(z)}^\mpi\) stands for the Moore-Penrose inverse of the matrix \(F(z)\). In~\zita{MR2988005} (see also~\zita{MR3014198}), we investigated the Moore-Penrose inverse of an arbitrary function \(F\) belonging to the class \(\RFq\). In particular, it turned out that \(-F^\mpi\) belongs to \(\RFq\) (see~\zitaa{MR3014198}{\cthm{9.4}}). The close relation between \(\SFqlep\) and \(\RFq\) suggests now to study the Moore-Penrose inverse of a function \(F\in\SFqlep\).

\bleml{P0848}
 Let \(\lep\in\R\) and \(F\in\SFqlep\). Then \(F^\mpi\) is holomorphic in \(\Cs\).
\elem
\bproof
 In view of formulas \eqref{Nr.NM} and \eqref{L374.1}, we obtain for all \(z\in\Cs\) the identities \(\nul{F(z)}=\nul{F(\iu)}\) and \(\ran{F(z)}=\ran{F(\iu)}\). Thus, the application of~\zitaa{MR3014197}{\cprop{8.4}} completes the proof.
\eproof

Let \(\lep\in\R\) and \(F\in\SFqlep\). Then \rlem{P0848} suggests to look if there are functions closely related to \(F^\mpi\) which belong again to \(\SFqlep\). Against to this background, we are led to the function \(G\colon\Cs\to\Cqq\) defined by \(G(z)\defeq-(z-\lep)^\inv\ek{F(z)}^\mpi\).

\breml{R1620}
 If \(A\in\CEPq\), \ie{}, if \(A\in\Cqq\) fulfills \(\ran{A^\ad}=\ran{A}\), then \(\im(A^\mpi)=-A^\mpi(\im A)(A^\mpi)^\ad\) (see also~\zitaa{MR3014198}{\cpropss{A.5}{A.6}}).
\erem

\bthml{L375}
 Let \(\lep\in\R\) and let \(F\in\SFqlep\). Then \(G\colon\Cs\to\Cqq\) defined by \(G(z)\defeq-(z-\lep)^\inv\ek{F(z)}^\mpi\) belongs to the class \(\SFqlep\) as well.
\ethm
\bproof
 \rlem{P0848} yields that the function \(F^\mpi\) is holomorphic in \(\Cs\). Consequently, the function \(G\) is holomorphic in \(\Cs\). Let \(z\in\ohe\). Using \rpropp{L374}{L374.a}, we have \(\ran{F^\ad(z)}=\ran{F(z)}\). Hence, because of \rrem{R1620}, the equations
 \beql{L375.1}
  \im G(z)
  =G(z)\rk*{\im\ek*{(z-\lep )F(z)}}G^\ad(z)
 \eeq
 and
 \beql{L375.2}
  \im\ek*{(z-\lep )G(z)}
  =\im\ek*{-F^\mpi(z)}
  =F^\mpi(z)\ek*{\im F(z)}\ek*{F^\mpi(z)}^\ad
 \eeq
 hold. Taking into account \(F\in\SFqlep\), the application of \rprop{M1.47N} yields
 \begin{align}\label{L375.3}
  \im F(z)&\in\Cggq&
 &\text{and}&
  \im\ek*{(z-\lep )F(z)}&\in\Cggq.
 \end{align}
 Thus, combining \eqref{L375.1} (resp.\ \eqref{L375.2}) and \eqref{L375.3}, we get \(\im G(z)\in\Cggq\) and \(\im\ek{(z-\lep )G(z)}\in\Cggq\). Now, the application of \rprop{M1.47N} yields \(G\in\SFqlep\).
\eproof

Now we specify the result of \rthm{L375} for functions belonging to \(\SOFqlep\).

\bpropl{L395}
 Let \(\lep\in\R\), let \(F\in\SOFqlep \), and let \(\stm{F}\) be the \tasmo{\lep }{\(F\)}. Then \(G\colon\Cs \to\Cqq\) defined by \(G(z)\defeq-(z-\lep )^\inv\ek{F(z)}^\mpi\) belongs to \(\SFqlep\) and
 \beql{L395.A}
  \gammau{G}
  =\ek*{\stma{F}{\rhl}}^\mpi.
\eeq
In particular, if \(F\) is not the constant function with value \(\Oqq\), then \(G\in\SFqlep\setminus\SOFqlep\).
\eprop
\bproof
 In view of \rthm{L375}, we have \(G\in\SFqlep\). From \rprop{R367Ad} we obtain
 \beql{Nr.GGY}
  \gammau{G}
  =\lim_{y\to+\infty}G(\iu y).
 \eeq
 Since \(F\) belongs to \(\SOFqlep \), we have \(\lim_{y\to+\infty}F(\iu y)=\Oqq\). \rprop{F394} yields \eqref{N14}. Consequently, \(\lim_{y\to+\infty}(\lep -\iu y)F(\iu y)=\stma{F}{\rhl}\). In view of \rprop{F394}, we have \(\ran{F(\iu y)}=\ran{\stma{F}{\rhl}}\) and, in particular, \(\rank[(\lep -\iu y)F(\iu y)]=\rank\stma{F}{\rhl}\) for each \(y\in(0,+\infty)\). Hence, taking into account~\zitaa{MR1105324}{\cthm{10.4.1}}, we obtain
 \beql{Nr.GGZ}
  \lim_{y\to+\infty}\rk*{\ek*{(\lep -\iu y)F(\iu y)}^\mpi}
  =\ek*{\stmA{F}{\rhl }}^\mpi.
 \eeq
 Since \(G(\iu y)=-(\iu y-\lep )^\inv F^\mpi(\iu y)=\ek{(\lep -\iu y)F(\iu y)}^\mpi\) holds true for each \(y\in(0,+\infty)\), from \eqref{Nr.GGY} and \eqref{Nr.GGZ} we get then \eqref{L395.A}. Now assume that \(G\) belongs to \(\SOFqlep\). From the definition of the class \(\SOFqlep\) we obtain then \(\lim_{y\to+\infty}G(\iu y)=\Oqq\), which, in view of \eqref{Nr.GGY} and \eqref{L395.A}, implies \(\stma{F}{\rhl}=\Oqq^\mpi=\Oqq\). \rprop{F394} yields then \(\nul{F(z)}=\Cq\) and hence \(F(z)=\Oqq\) for all \(z\in\Cs\). This proves \(G\notin\SOFqlep\) if \(F\) is not the constant function with value \(\Oqq\).
\eproof

For the special choice \(q=1\) and \(\lep=0\), the following class was introduced by Kats/Krein~\zitaa{KK74}{\cdefn{D1.5.2}}.
\bnotal{D0924}
 Let \(\lep\in\R\). Then let \(\SiFqlep\)\index{s@\(\SiFqlep\)} be the class of all matrix-valued functions \(F\colon\Cs\to\Cqq\) which fulfill the following two conditions:
 \bAeqi{0}
  \il{D0924.I} \(F\) is holomorphic in \(\Cs\) with \(\Rstr_\ohe F\in\RFq\).
  \il{D0924.II} For all \(x\in\lhl\), the matrix \(-F(x)\) is \tnnH{}.
 \eAeqi
\enota

\breml{R1228}
 Let \(\lep\in\R\), then \(F\in\Cs\to\Cqq\) belongs to \(\SiFqlep\) if and only if \(u^\ad Fu\in\SiFuu{1}{\lep}\) for all \(u\in\Cq\).
\erem

\bexal{E0808}
 Let \(\lep\in\R\) and let \(D,E\in\Cggq\). Then \(F\colon\Cs\to\Cqq\) defined by \(F(z)\defeq-D+(z-\lep)E\) belongs to \(\SiFqlep\).
\eexa

\bleml{L1321}
 Let \(\lep\in\R\) and let \(f\colon\Cs\to\C\) be such that there are real numbers \(d\) and \(e\) and a finite signed measure \(\rho\) on \(((\lep,+\infty),\Bsa{(\lep,+\infty)})\) such that \(f(z)=-d +(z-\lep)\ek{e +\int_{(\lep,+\infty)}(1+t-\lep)/(t-z)\rho(\dif t)}\) holds true for all \(z\in\Cs\). Then \(d=-\lim_{x \to+0}f(\lep-x )\) and \(e=-\lim_{x\to+\infty}\ek{f(\lep-x)+d}/x\). Furthermore, \(d\), \(e\), and \(\rho\) are uniquely determined.
\elem
\bproof
 With \(z_x \defeq\lep-x \) we have \((1+t-\lep)/(t-z_x)=(1+t-\lep)/(t-\lep+x)\geq0\) for all \(t\in(\lep,+\infty)\) and all \(x \in(0,+\infty)\), which decreases to \(0\) as \(x \) increases to infinity. Since the signed measure \(\rho\) is finite, its Jordan decomposition \(\rho=\rho_+-\rho_-\) consists of two finite measures. Hence, \(\int_{(\lep,+\infty)}(1+t-\lep)/(t-z_1)\rho_\pm(\dif t)=\rho_\pm((\lep,+\infty))<\infty\) holds true. Thus, Lebesgue's monotone convergence theorem yields \(\lim_{x \to+\infty}\int_{(\lep,+\infty)}(1+t-\lep)/(t-z_x)\rho_\pm(\dif t)=0\), which implies \(\lim_{x\to+\infty}\ek{f(z_x)+d}/(z_x-\lep)=e\). Furthermore, we have \(-(z_x -\lep)(1+t-\lep)/(t-z_x)=(1+t-\lep)/[1+(t-\lep)/x ]\geq0\) for all \(t\in(\lep,+\infty)\) and all \(x \in(0,+\infty)\), which decreases to \(0\) as \(x \) decreases to \(0\). Since \(\int_{(\lep,+\infty)}\ek{-(z_1-\lep)(1+t-\lep)/(t-z_1)}\rho_\pm(\dif t)=\rho_\pm((\lep,+\infty))<\infty\) holds true, Lebesgue's monotone convergence theorem yields \(\lim_{x \to+0}\int_{(\lep,+\infty)}\ek{-(z_x -\lep)(1+t-\lep)/(t-z_x )}\rho_\pm(\dif t)=0\)  showing \(-\lim_{x \to+0}f(z_x )=d\). In particular, \(d\) and \(e\) are uniquely determined. Now let \(\sigma\) be an arbitrary finite signed measure on \(((\lep,+\infty),\Bsa{(\lep,+\infty)})\) such that \(f(z)=-d +(z-\lep)\ek{e +\int_{(\lep,+\infty)}(1+t-\lep)/(t-z)\sigma(\dif t)}\) holds true for all \(z\in\Cs\). Then \(\int_{(\lep,+\infty)}(1+t-\lep)/(t-z)\sigma(\dif t)=\int_{(\lep,+\infty)}(1+t-\lep)/(t-z)\rho(\dif t)\) for all \(z\in\Cs\). Since the signed measure \(\sigma\) is finite, its Jordan decomposition \(\sigma=\sigma_+-\sigma_-\) consists of two finite measures. Hence, we obtain \(\int_{(\lep,+\infty)}(1+t-\lep)/(t-z)(\sigma_++\rho_-)(\dif t)=\int_{(\lep,+\infty)}(1+t-\lep)/(t-z)(\rho_++\sigma_-)(\dif t)\) for all \(z\in\Cs\) with finite measures \(\sigma_++\rho_-\) and \(\rho_++\sigma_-\). Using \rthm{368}, it is readily checked then that \(\sigma_++\rho_-\) and \(\rho_++\sigma_-\) coincide. Consequently, \(\sigma=\rho\) follows.
\eproof

\bleml{L1110}
 Let \(\lep\in\R\) and let \(f\in\SiFuu{1}{\lep}\). Then there are unique \tnn{} real numbers \(d\) and \(e\) and a unique measure \(\rho\in\Mggoa{1}{(\lep,+\infty)}\) such that \(f(z)=-d +(z-\lep)\ek{e +\int_{(\lep,+\infty)}(1+t-\lep)/(t-z)\rho(\dif t)}\) holds true for all \(z\in\Cs\).
\elem
\bproof
 Obviously, the function \(g\colon\C\setminus[0,+\infty)\to\C\) defined by \(g(w)\defeq f(w+\lep)\) belongs to \(\SiFuu{1}{0}\). Hence, by virtue of~\zitaa{KK74}{\cthm{S1.5.2}}, there exist unique numbers \(a\in(-\infty,0]\) and \(b\in[0,+\infty)\) and a unique measure \(\tau\) on \(((0,+\infty),\Bsa{(0,+\infty)})\) with \(\int_{(0,+\infty)}1/(x+x^2)\tau(\dif x)<\infty\) such that \(g(w)=a+bw+\int_{(0,+\infty)}\ek{1/\rk{x-w}-1/x}\tau(\dif x)\)
 for all \(w\in\C\setminus[0,+\infty)\). Then \(\chi\) defined on \(\Bsa{(0,+\infty)}\) by \(\chi(B)\defeq\int_B1/(x+x^2)\tau(\dif x)\) is a finite measure on \(((0,+\infty),\Bsa{(0,+\infty)})\) and the integral \(\int_{(0,+\infty)}(x+x^2)\ek*{1/\rk{x-w}-1/x}\chi(\dif x)\)
 exists for all \(w\in\C\setminus[0,+\infty)\) and equals to \(\int_{(0,+\infty)}\ek{1/(x-w)-1/x}\tau(\dif x)\). We have \((x+x^2)\ek{1/(x-w)-1/x}=w(1+x)/(x-w)\) for all \(w\in\C\setminus[0,+\infty)\) and all \(x\in(0,+\infty)\). In view of \(z-\lep\in\C\setminus[0,+\infty)\) for all \(z\in\Cs\), we obtain thus
 \[\begin{split}
  f(z)
  &=g(z-\lep)
  =a+b(z-\lep)+\int_{(0,+\infty)}\ek*{\frac{1}{x-(z-\lep)}-\frac{1}{x}}\tau(\dif x)\\
  &=a+(z-\lep)\ek*{b+\int_{(0,+\infty)}\frac{1+x}{x-(z-\lep)}\chi(\dif x)}\\
  &=-d+(z-\lep)\ek*{e+\int_{(\lep,+\infty)}\frac{1+t-\lep}{t-z}\rho(\dif t)}
 \end{split}\]
 for all \(z\in\Cs\), where \(d\defeq-a\), \(e\defeq b\), and \(\rho\) is the image measure of \(\chi\) under the translation \(T\colon(0,+\infty)\to(\lep,+\infty)\) defined by \(T(x)\defeq x+\lep\). In particular \(d,e\in[0,+\infty)\) and \(\rho\in\Mggoa{1}{(\lep,+\infty)}\). Hence, the triple \((d,e,\rho)\) is unique by virtue of \rlem{L1321}. 
\eproof

\bthml{T1305}
 Let \(\lep\in\R\) and let  \(F\colon\Cs\to\Cqq\). Then:
 \benui
  \il{T1305.a} If \(F\in\SiFqlep\), then there are unique \tnnH{} complex \tqqa{matrices} \(D\) and \(E\) and a unique \tnnH{} measure \(\rho\in\Mggqa{(\lep,+\infty)}\) such that
  \beql{T1305.A}
   F(z)
   =-D +(z-\lep)\ek*{E +\int_{(\lep,+\infty)}\frac{1+t-\lep}{t-z}\rho(\dif t)}
  \eeq
  for all \(z\in\Cs\). Furthermore, the function \(P\colon\Cs\to\Cqq\) defined by \(P(z)\defeq(z-\lep)^\inv F(z)\) belongs to \(\SFqlep\) with \(D=\muPa{\set{\lep}}\) and \((E,\rho)=(\gammaP,\Rstr_\Bsa{(\lep,+\infty)}\muP)\).
  \il{T1305.b} If \(D\in\Cggq\), \(E \in\Cggq\), and \(\rho\in\Mggqa{(\lep,+\infty)}\) are such that \(F\) can be represented via \eqref{T1305.A} for all \(z\in\Cs\), then \(F\) belongs to \(\SiFqlep\).
 \eenui
\ethm
\bproof
 \eqref{T1305.a} We consider an arbitrary vector \(u\in\Cq\). According to \rrem{R1228}, then \(f_u\defeq u^\ad Fu\) belongs to \(\SiFuu{1}{\lep}\). Hence, \rlem{L1110} yields the existence of a unique triple \((d_u,e_u,\rho_u)\in[0,+\infty)\times[0,+\infty)\times\Mggoa{1}{(\lep,+\infty)}\) such that
 \beql{T1305.1}
  f_u(z)
  =-d_u+(z-\lep)\ek*{e_u+\int_{(\lep,+\infty)}\frac{1+t-\lep}{t-z}\rho_u(\dif t)}
 \eeq
 holds true for all \(z\in\Cs\).
 With the standard basis \((\eu{1},\eu{2},\dotsc,\eu{q})\) of \(\Cq\), let \(d_{jk}\defeq\frac{1}{4}\sum_{\ell=0}^3(-\iu)^\ell d_{\eu{j}+\iu^\ell\eu{k}}\), \(e_{jk}\defeq\frac{1}{4}\sum_{\ell=0}^3(-\iu)^\ell e_{\eu{j}+\iu^\ell\eu{k}}\) and \(\rho_{jk}\defeq\frac{1}{4}\sum_{\ell=0}^3(-\iu)^\ell \rho_{\eu{j}+\iu^\ell\eu{k}}\) for all \(j,k\in\mn{1}{q}\). We have then
 \[
  \begin{split}
   &\eu{j}^\ad\ek*{F(z)}\eu{k}
   =\frac{1}{4}\sum_{\ell=0}^3(-\iu)^\ell f_{\eu{j}+\iu^\ell\eu{k}}(z)\\
   &=\frac{1}{4}\sum_{\ell=0}^3(-\iu)^\ell\rk*{-d_{\eu{j}+\iu^\ell\eu{k}}+(z-\lep)\ek*{e_{\eu{j}+\iu^\ell\eu{k}}+\int_{(\lep,+\infty)}\frac{1+t-\lep}{t-z}\rho_{\eu{j}+\iu^\ell\eu{k}}(\dif t)}}\\
   &=-d_{jk}+(z-\lep)\ek*{e_{jk}+\int_{(\lep,+\infty)}\frac{1+t-\lep}{t-z}\rho_{jk}(\dif t)}
  \end{split}
 \]
for all \(j,k\in\mn{1}{q}\) and all \(z\in\Cs\). Hence, \eqref{T1305.A} follows for all \(z\in\Cs\) with \(D\defeq\matauuo{d_{jk}}{j,k}{1}{q}\), \(E\defeq\matauuo{e_{jk}}{j,k}{1}{q}\), and \(\rho\defeq\matauuo{\rho_{jk}}{j,k}{1}{q}\). For all \(\zeta\in\C\) with \(\abs{\zeta}=1\), we have \(f_{\zeta u}=f_u\) and thus \(d_{\zeta u}=d_u\), \(e_{\zeta u}=e_u\), and \(\rho_{\zeta u}=\rho_u\) by virtue of the uniqueness of the triple \((d_u,e_u,\rho_u)\). A straightforward calculation yields for all \(j,k\in\mn{1}{q}\) then \(d_{kj}=\ko{d_{jk}}\), \(e_{kj}=\ko{e_{jk}}\), and \(\rho_{kj}(B)=\ko{\rho_{jk}(B)}\) for all \(B\in\Bsa{(\lep,+\infty)}\). Thus, the matrices \(D\) and \(E\) are \tH{} and \(\rho\) is an \(\sigma\)\nobreakdash-additive mapping defined on \(\Bsa{(\lep,+\infty)}\) with values in \(\CHq\). From \eqref{T1305.A} we obtain \(f_u(z)=-u^\ad Du+(z-\lep)\rk{u^\ad Eu +\int_{(\lep,+\infty)}\ek{(1+t-\lep)/(t-z)}(u^\ad \rho u)(\dif t)}\) for all \(z\in\Cs\), where \(u^\ad Du\) and \(u^\ad Eu\) belong to \(\R\) and \(u^\ad \rho u\) is a finite signed measure on \(((\lep,+\infty),\Bsa{(\lep,+\infty)})\). In view of \eqref{T1305.1}, \rlem{L1321} yields \(u^\ad Du=d_u\), \(u^\ad Eu=e_u\), and \(u^\ad \rho u=\rho_u\). In particular, \(u^\ad Du\) and \(u^\ad Eu\) belong to \([0,+\infty)\) and \(u^\ad \rho u\in\Mggoa{1}{(\lep,+\infty)}\). Since \(u\in\Cq\) was arbitrarily chosen, hence \(D,E\in\Cggq\), and \(\rho\in\Mggqa{(\lep,+\infty)}\) follow.

 Now let \(D ,E \in\Cggq\), and \(\rho\in\Mggqa{(\lep,+\infty)}\) be such that \eqref{T1305.A} holds true for all \(z\in\Cs\). Denote by \(\Kron{\lep}\) the Dirac measure on \((\rhl,\Bsarhl)\) with unit mass at \(\lep\). Then \(P\) admits the representation \(P(z)=E +\int_{[\lep,+\infty)}(1+t-\lep)/(t-z)\theta(\dif t)\) for all \(z\in\Cs\), where \(\theta\colon\Bsarhl\to\Cggq\) defined by \(\theta(B)\defeq\rho(B\cap(\lep,+\infty))+\ek{\Kron{\lep}(B)}D\) belongs to \(\Mggqrhl\). Hence, \rthmp{368}{368.b} and \rrem{JB} yield \(P\in\SFqlep\) with \(\gammaP=E\) and \(\muP=\theta\). In particular, \(\muPa{\set{\lep}}=D\) and \(\Rstr_\Bsa{(\lep,+\infty)}\muP=\rho\). Hence, the triple \((D,E,\rho)\) is unique.
 
 \eqref{T1305.b} Let \(D ,E \in\Cggq\) and \(\rho\in\Mggqa{(\lep,+\infty)}\) be such that \eqref{T1305.A} holds true for all \(z\in\Cs\). As explained above, then \(P\) belongs to \(\SFqlep\). Since \(F(z)=(z-\lep)P(z)\) for all \(z\in\Cs\), we hence conclude with \rprop{M1.47N} that \(F\) belongs to \(\SiFqlep\).
\eproof

In the following, if \(\lep\in\R\) and \(F\in \SiFqlep \) are given, then we will write \((\DF,\EF,\rhoF)\)\index{d@\(\DF\)}\index{e@\(\EF\)}\index{r@\(\rhoF\)} for the unique triple \((D,E,\rho)\) from \(\Cggq\times\Cggq\times\Mggqa{(\lep,+\infty)}\) which fulfills \eqref{T1305.A} for all \(z\in\Cs\).

\bcorl{C1611}
 Let \(\lep\in\R\). Then \(F\colon\Cs\to\Cqq\) belongs to \(\SiFqlep\) if and only if \(P\colon\Cs\to\Cqq\) defined by \(P(z)=(z-\lep)^\inv F(z)\) belongs to \(\SFqlep\).
\ecor
\bproof
 If \(F\in\SiFqlep\), then \(P\in\SFqlep\), by virtue of \rthmp{T1305}{T1305.a}.
 
 Conversely, now suppose \(P\in\SFqlep\). According to \rthm{368} and \rrem{JB}, then
 \[\begin{split}
  F(z)
  =(z-\lep)P(z)
  &=(z-\lep)\ek*{\gammaP +\int_{\rhl} \frac{1+ t-\alpha}{t-z} \muPa{\dif t}}\\
  &=-D+(z-\lep)\ek*{E +\int_{(\lep,+\infty)}\frac{1+t-\lep}{t-z}\rho(\dif t)}
 \end{split}\]
 for all \(z\in\Cs\), where the matrices \(D\defeq\muPa{\set{\lep}}\) and \(E\defeq\gammaP\) are \tnnH{} and \(\rho\defeq\Rstr_\Bsa{(\lep,+\infty)}\muP\) belongs to \(\Mggqa{(\lep,+\infty)}\). Hence, \rthmp{T1305}{T1305.b} yields \(F\in\SiFqlep\).
\eproof

\bcorl{C0901}
 Let \(\lep\in\R\) and let \(F\in\SiFqlep\). For all \(x_1,x_2\in\lhl\) with \(x_1\leq x_2\), then \(F(x_1)\leq F(x_2)\leq\Oqq\).
\ecor
\bproof
 Using \rthm{T1305}, we obtain
 \[
  F(x_2)-F(x_1)
  =(x_2-x_1)\ek*{\EF+\int_{\rhl}\frac{(1+t-\lep)(t-\lep)}{(t-x_2)(t-x_1)}\rhoFa{\dif t}}
 \]
 for all \(x_1,x_2\in\lhl\) with \(x_1\leq x_2\), by direct calculation. Since \(\EF\in\Cggq\) and \(-F(x)\in\Cggq\) for all \(x\in\lhl\), thus the proof is complete.
\eproof

Now we consider again the situation of \rexa{E0808}.
\bexal{E0914}
 Let \(\lep\in\R\) and let \(D,E\in\Cggq\). Then \(F\colon\Cs\to\Cqq\) defined by \(F(z)\defeq-D+(z-\lep)E\) belongs to \(\SiFqlep\), where \(\DF=D\), \(\EF=E\), and \(\rhoF\) is the constant measure with value \(\Oqq\).
\eexa

\bpropl{P1622}
 Let \(\lep \in\R\) and let \(F\in\SiFqlep \). Then:
 \benui
  \il{P1622.a} If \(z\in\Cs\), then \(\ko z\in\Cs\) and \(\ek{F(z)}^\ad=F(\ko z)\).
  \il{P1622.b} For all \(z\in \Cs \),
  \begin{align}
   \Nul{F(z)}&=\nul{\DF}\cap\nul{\EF}\cap\Nul{\rhoFA{(\lep,+\infty)}},\label{P1622.A}\\
   \Ran{F(z)}&=\ran{\DF}+\ran{\EF}+\Ran{\rhoFA{(\lep,+\infty)}},\label{P1622.B}
  \end{align}
  and, in particular, \(\nul{\ek{F(z)}^\ad}=\nul{F(z)}\) and \(\ran{\ek{F(z)}^\ad}=\ran{F(z)}\).
  \il{P1622.c} Let \(r\in\NO\). Then the following statements are equivalent:
  \baeqii{0}
   \item \(\rank  F(z) = r\) for all \(z\in\Cs \).
   \item There is some \(z_0\in\Cs \) such that \(\rank  F(z_0) = r\).
   \item \(\dim[\ran{\DF}+\ran{\EF}+\ran{\rhoFa{(\lep,+\infty)}}]= r\).
  \eaeqii
 \eenui
\eprop
\bproof
 \eqref{P1622.a} This can be seen from \rthmp{T1305}{T1305.a}.
 
 \eqref{P1622.b} According to \rthmp{T1305}{T1305.a}, the function \(P\colon\Cs\to\Cqq\) defined by \(P(z)\defeq(z-\lep)^\inv F(z)\) belongs to \(\SFqlep\) with \(\DF=\muPa{\set{\lep}}\) and \((\EF,\rhoF)=(\gammaP,\Rstr_\Bsa{(\lep,+\infty)}\muP)\). In particular, \(\muPa{\rhl}=\DF+\rhoFa{(\lep,+\infty)}\). Since the two matrices on the right-hand side of the last equation are both \tnnH{}, we get \(\nul{\muPa{\rhl}}=\nul{\DF}\cap\nul{\rhoFa{(\lep,+\infty)}}\). Now let \(z\in\Cs\). Applying \rpropp{L374}{L374.a} to \(P\), we get \(\nul{P(z)}=\nul{\gammaP}\cap\Nul{\muPa{\rhl}}\). In view of \(\nul{F(z)}=\nul{P(z)}\), then \eqref{P1622.A} follows.
 Thus, \eqref{P1622.A} is proved for all \(z\in\Cs\).
 From~\eqref{P1622.a} and \eqref{P1622.A} we get \(\nul{\ek{F(z)}^\ad}=\nul{F(z)}\) for all \(z\in\Cs\). Taking additionally into account that the matrices \(\DF\), \(\EF\), and \(\rhoFa{(\lep,+\infty)}\) are \tnnH{}, we obtain \eqref{P1622.B} from \eqref{P1622.A} in the same way as in the proof of \rpropp{L374}{L374.a}. Using~\eqref{P1622.a} and \eqref{P1622.B}, we get \(\ran{\ek{F(z)}^\ad}=\ran{F(z)}\) for all \(z\in\Cs\).
 
 \eqref{P1622.c}  This is a consequence of~\eqref{P1622.B}. 
\eproof

\bcorl{C0846}
 Let \(\lep \in\R\), let \(F\in\SiFqlep \), and let \(z_0\in\Cs\). Then \(F(z_0)=\Oqq\) if and only if \(F(z)=\Oqq\) for all \(z\in\Cs\).
\ecor
\bproof
 This is an immediate consequence of \rpropp{P1622}{P1622.c}.
\eproof

\bcorl{C0847}
 Let \(\lep \in\R\), let \(F\in\SiFqlep \), and let \(\lambda\in\R\) be such that the matrix \(\DF+\lambda\Iq\) is \tnnH{}. Then
 \begin{align}\label{C0847.B1}
  \Ran{F(z)-\lambda\Iq}&=\Ran{F(w)-\lambda\Iq}&
  &\text{and}&
  \Nul{F(z)-\lambda\Iq}&=\Nul{F(w)-\lambda\Iq}
 \end{align}
 for every choice of \(z\) and \(w\) in \(\Cs\). In particular, if \(\lambda\geq0\), then \(\lambda\) is an eigenvalue of the matrix \(F(z_0)\) for some \(z_0\in\Cs\) if and only if \(\lambda\) is an eigenvalue of the matrix \(F(z)\) for all \(z\in\Cs\). In this case, the eigenspaces \(\nul{F(z)-\lambda\Iq}\) are independent of \(z\in\Cs\).
\ecor
\bproof
 In view of \rthm{T1305}, we conclude that the function \(G\colon\Cs\to\Cqq\) defined by \(G(z)\defeq F(z)-\lambda\Iq\) belongs to \(\SiFqlep\). The application of \rpropp{P1622}{P1622.b} to the function \(G\) yields then \eqref{C0847.B1}. Since the matrix \(\DF\) is \tnnH{}, we have \(\DF+\lambda\Iq\in\Cggq\) if \(\lambda\geq0\). Thus, the remaining assertions are an immediate consequence of \eqref{C0847.B1}.
\eproof

\bleml{L1004}
 Let \(\lep\in\R\) and \(F\in\SiFqlep\). Then \(F^\mpi\) is holomorphic in \(\Cs\).
\elem
\bproof
 In view of \eqref{P1622.A} and \eqref{P1622.B}, we obtain for all \(z\in\Cs\) the identities \(\nul{F(z)}=\nul{F(\iu)}\) and \(\ran{F(z)}=\ran{F(\iu)}\). Thus, the application of~\zitaa{MR3014197}{\cprop{8.4}} completes the proof.
\eproof

The following proposition is a generalization of a result due to Kats and Krein~\zitaa{KK74}{\clem{D1.5.2}}, who considered the case \(q=1\) and \(\lep=0\).
\bthml{P0931}
 Let \(\lep\in\R\) and let \(F\colon\Cs\to\Cqq\) be a matrix-valued function. Then \(F\) belongs to \(\SiFqlep\) if and only if \(G\defeq-F^\mpi\) belongs to \(\SFqlep\).
\ethm
\bproof
 First suppose \(F\in\SiFqlep\). Then \(G\) is holomorphic in \(\Cs\) by virtue of \rlem{L1004}. In view of \(\Rstr_\ohe F\in\RFq\), we conclude from~\zitaa{MR2988005}{\cprop{3.8}} that \(\Rstr_\ohe G\) belongs to \(\RFq\) as well. In particular, \(\im G(w)\in\Cggq\) for all \(w\in\ohe\). Because of \(-F(x)\in\Cggq\) for each \(x\in(-\infty,\lep)\), we have \(G(x)=\ek{-F(x)}^\mpi\in\Cggq\) for all \(x\in(-\infty,\lep)\) (see, \eg{}~\zitaa{MR1152328}{\clem{1.1.5}}). Hence, \(G\) belongs to \(\SFqlep\).
 
 Now suppose \(G\in\SFqlep\). Then \(F=-G^\mpi\). Hence, \(F\) is holomorphic in \(\Cs\) by virtue of \rlem{P0848}. Since \rprop{M1.47N} yields \(\Rstr_\ohe G\in\RFq\), we conclude from~\zitaa{MR2988005}{\cprop{3.8}} that \(\Rstr_\ohe F\) belongs to \(\RFq\) as well. Because of \(G(x)\in\Cggq\) for all \(x\in(-\infty,\lep)\), we have \(-F(x)=\ek{G(x)}^\mpi\in\Cggq\) for all \(x\in(-\infty,\lep)\) (see, \eg{}~\zitaa{MR1152328}{\clem{1.1.5}}). Hence, \(F\in\SiFqlep\).
\eproof

\section{Integral representations for the class $\TFqrep$}\label{S1159}
The main goal of this section is to derive some integral representations for \trTF{s} of order \(q\). Our strategy is based on using the corresponding results for the class \(\SFuu{q}{-\rep}\). The following observation provides the key to realize our aims.

\breml{R0948}
 Let \(\lep,\rep\in\R\) and let \(T\colon\C\to\C\) be defined by \(T(z)\defeq\lep+\rep-\ko z\). Then \(T(\Cs)=\Ct\), \(T(\Ct)=\Cs\), \(T(\ohe)=\ohe\), \(T(\lhl)=\brhl\) and \(T(\brhl)=\lhl\). Consequently, in view of \rprop{214-1}, one can easily check that, for each \(F\in\SFqlep\), the function \(G\colon\Ct\to\Cqq\) defined by \(G(z)\defeq-\ek{F(\lep+\rep-\ko z)}^\ad\) belongs to \(\TFqrep\) and that, conversely, for any function  \(G\in\TFqrep\), the function \(F\colon\Cs\to\Cqq\) defined by \(F(z)\defeq-\ek{G(\lep+\rep-\ko z)}^\ad\) belongs to \(\SFqlep\).
\erem

\bpropl{P1502}
 Let \(\rep\in\R\) and let \(G\in\TFqrep \). Then the Nevanlinna parametrization \((\Npa,\Npb,\Npnu)\) of \(\Rstr_{\ohe }G\) fulfills \(\Npnua{\brhl}= \Oqq\), \(\Npb= \Oqq\), and \(\Npnu\in\MgguqR{1}\). In particular, for each \(z\in\Ct\), then \(G(z)=\Npa +\int_{\blhl}\rk{1+tz}/\rk{t-z} \Npnua{\dif t}\).
\eprop
\bproof
 According to \rrem{R0948}, the function \(F\colon\C\setminus[-\rep,+\infty)\to\Cqq\) defined by \(F(z)\defeq-\ek{G(-\ko z)}^\ad\) belongs to \(\SFuu{q}{-\rep}\). From \rrem{R1602} we obtain then that the Nevanlinna parametrization of \(\Rstr_{\ohe}F\) is given by \((-\Npa,\Npb,\theta)\), where \(\theta\) is the image measure of \(\Npnu\) under the reflection \(t\mapsto-t\) on \(\R\). Now \rprop{HCT} yields \(\theta((-\infty,-\rep))=\Oqq\), \(\Npb=\Oqq\), \(\theta\in\MgguqR{1}\), and \(F(z)=-\Npa+\int_{[-\rep,+\infty)}(1+tz)/(t-z) \theta(\dif t)\) for all \(z\in\C\setminus[-\rep,+\infty)\). Hence, \(\NpnuA{\brhl}= \Oqq\), \(\Npnu\in\MgguqR{1}\), and
 \[
  G(z)
  =-\ek*{F(-\ko z)}^\ad
  =A+\int_{[-\rep,+\infty)}\frac{1-tz}{-t-z}\theta(\dif t)
  =A+\int_{\blhl} \frac{1+tz}{t-z} \Npnua{\dif t}
 \]
 for all \(z\in\Ct\).
\eproof

\bthml{T1557}
 Let \(\rep\in\R\) and let \(G\colon\Ct\to\Cqq\). Then:
 \benui
  \il{T1557.a} Suppose \(G\in\TFqrep \). Denote by \((\Npa,\Npb,\Npnu)\) the Nevanlinna parametrization of \(\Rstr_{\ohe } G\) and let \(\tilde\Npnu\defeq\Rstr_{\Bsablhl}\Npnu\). Then \(\tilde\Npnu\in\Mgguqblhl{1}\) and there is a unique pair \((C,\eta)\) from \(\Cggq\times\Mgguqblhl{1}\) such that
  \beql{T1557.B1}
   G(z)
   =-C+\int_{\blhl}\frac{1+t^2}{t-z}\eta(\dif t)
  \eeq
  for all \(z\in\Ct\), namely \(C=\int_{\blhl}t\tilde\Npnu(\dif t)-\Npa\) and \(\eta=\tilde\Npnu\). Furthermore, \(C=-\lim_{r\to+\infty}G(\rep+r\ec^{\iu\phi})\) for all \(\phi\in(-\pi/2,\pi/2)\).
  \il{T1557.b} Let \(C\in\Cggq\) and let \(\eta\in\Mgguqblhl{1}\) be such that \eqref{T1557.B1} holds true for all \(z\in\Ct\). Then \(G\) belongs to \(\TFqrep\).
 \eenui
\ethm
\bproof
 \eqref{T1557.a} According to \rrem{R0948}, the function \(F\colon\C\setminus[-\rep,+\infty)\to\Cqq\) defined by \(F(z)\defeq-\ek{G(-\ko z)}^\ad\) belongs to \(\SFuu{q}{-\rep}\). From \rrem{R1602} we obtain then that the Nevanlinna parametrization of \(\Rstr_{\ohe}F\) is given by \((-\Npa,\Npb,\theta)\), where \(\theta\) is the image measure of \(\Npnu\) under the reflection \(t\mapsto-t\) on \(\R\). Now \rthmp{T1542}{T1542.a} yields that \(\tilde\theta\defeq\Rstr_{\Bsa{[-\rep,+\infty)}}\theta\) belongs to \(\Mgguqa{1}{[-\rep,+\infty)}\) and that there is a unique pair \((D,\tau)\in\Cggq\times\Mgguqa{1}{[-\rep,+\infty)}\) such that \(F(z)=D+\int_{[-\rep,+\infty)}(1+t^2)/(t-z)\tau(\dif t)\) for all \(z\in\C\setminus[-\rep,+\infty)\), namely \(D=-\Npa-\int_{[-\rep,+\infty)}t\tilde\theta(\dif t)\) and \(\tau=\tilde\theta\). Since \(\tilde\theta\) is the image measure of \(\tilde\Npnu\) under the transformation \(T\colon\blhl\to[-\rep,+\infty)\) defined by \(T(t)\defeq-t\), we can conclude \(\tilde\Npnu\in\Mgguqblhl{1}\) and
 \[
  \begin{split}
   G(z)
   &=-\ek*{F(-\ko z)}^\ad
   =-D^\ad-\ek*{\int_{[-\rep,+\infty)}\frac{1+t^2}{t+\ko z}\tau(\dif t)}^\ad
   =-D-\int_{[-\rep,+\infty)}\frac{1+t^2}{t+z}\tau(\dif t)\\
   &=-\ek*{\int_{[-\rep,+\infty)}(-t)\tilde\theta(\dif t)-\Npa}+\int_{[-\rep,+\infty)}\frac{1+t^2}{-t-z}\tilde\theta(\dif t)
   =-C+\int_{\blhl}\frac{1+t^2}{t-z}\eta(\dif t)
  \end{split}
 \]
 for all \(z\in\Ct\), where \(C\defeq\int_{\blhl}t\tilde\Npnu(\dif t)-\Npa\) and \(\eta\defeq\tilde\Npnu\). From the above computation we see \(C=D\) and hence \(C\in\Cggq\) follows. Taking additionally into account \rthmp{T1542}{T1542.a}, for all \(\phi\in(-\pi/2,\pi/2)\), we get
 \[
  C
  =D^\ad
  =\ek*{\lim_{r\to+\infty}F(-\rep+r\ec^{\iu(\pi-\phi)})}^\ad
  =\lim_{r\to+\infty}\ek*{F(-\rep-r\ec^{-\iu\phi})}^\ad
  =-\lim_{r\to+\infty}G(\rep+r\ec^{\iu\phi}).
 \]
 
 Now let \(C\in\Cggq\) and \(\eta\in\Mgguqblhl{1}\) be such that \eqref{T1557.B1} holds true for all \(z\in\Ct\). Then \(\chi\colon\BsaR\to\Cggq\) defined by \(\chi(M)\defeq\eta(M\cap\blhl)\) belongs to \(\MgguqR{1}\) and the matrix \(-C+\int_\R t\chi(\dif t)\) is \tH{}. Using \eqref{T1542.0}, we conclude from \eqref{T1557.B1} then that the integral \(\int_\R(1+tz)/(t-z)\chi(\dif t)\) exists and that \(G(z)=-C+\int_\R t\chi(\dif t)+z\cdot\Oqq+\int_\R(1+tz)/(t-z)\chi(\dif t)\) for all \(z\in\ohe\). \rthmp{T1554}{T1554.a} yields then \(-C+\int_\R t\chi(\dif t)=\Npa\) and \(\chi=\Npnu\). Hence \(\eta=\tilde\Npnu\) and \(C=\int_{\blhl}t\tilde\Npnu(\dif t)-\Npa\).
 
 \eqref{T1557.b} Let \(C\in\Cggq\) and \(\eta\in\Mgguqblhl{1}\) be such that \eqref{T1557.B1} holds true for all \(z\in\Ct\). Using a result on holomorphic dependence of an integral on a complex parameter (see, \eg{}~\zitaa{MR2257838}{Ch.~IV, \S5, \cSatz{5.8}}), we conclude then that \(G\) is a matrix-valued function which is holomorphic in \(\Ct\). Furthermore,
 \[
  \im G(z)
  =\int_{\blhl}\im\rk*{\frac{1+t^2}{t-z}}\eta(\dif t)
  =\int_{\blhl}\frac{(1+t^2)\im z}{\abs{t-z}^2}\eta(\dif t)
  \in\Cggq
 \]
 for all \(z\in\ohe\) and
 \[
  -G(x)
  =C+\int_{\blhl}\frac{1+t^2}{x-t}\eta(\dif t)
  \in\Cggq
 \]
 for all \(x\in\brhl\). Thus, \(G\) belongs to \(\TFqrep\).
\eproof

In the following, if \(\rep\in\R\) and \(G\in \TFqrep \) are given, then we will write \((\CG,\etaG)\)\index{c@\(\CG\)}\index{e@\(\etaG\)} for the unique pair \((C,\eta)\) from \(\Cggq\times\Mgguqblhl{1}\) which fulfills \eqref{T1557.B1} for all \(z\in\Ct\).

\breml{R0843}
 Let \(\rep\in\R\) and let \(G\in\TFqrep\). For all \(x_1,x_2\in\brhl\) with \(x_1\leq x_2\), then \(G(x_1)\leq G(x_2)\leq\Oqq\), by virtue of \rthmp{T1557}{T1557.a}.
\erem

\breml{R0941}
 \benui
  \il{R0941.a} Let \(\lep\in\R\) and let \(F\in \SFqlep\). In view of \rthm{T1542}, we have then
 \[
  -\ek*{F(-\ko z)}^\ad
  =-\CF+\int_{\rhl}\frac{1+t^2}{-t-z}\etaF(\dif t)
  =-\CF+\int_{(-\infty,-\lep]}\frac{1+t^2}{t-z}\hat\theta(\dif t)
 \]
 for all \(z\in\C\setminus(-\infty,-\lep]\), where \(\hat\theta\) is the image measure of \(\etaF\) under the transformation \(R\colon\rhl\to(-\infty,-\lep]\) defined by \(R(t)\defeq-t\). Because of \(\CF\in\Cggq\) and \(\hat\theta\in\Mgguqa{1}{(-\infty,-\lep]}\), \rthm{T1557} yields then that \(G\colon\C\setminus(-\infty,-\lep]\to\Cqq\) defined by \(G(z)\defeq-\ek{F(-\ko z)}^\ad\) belongs to \(\TFuu{q}{-\lep}\) and that \((\CG,\etaG)=(\CF,\hat\theta)\).
  \il{R0941.b} Let \(\rep\in\R\) and let \(G\in \TFqrep\). In view of \rthm{T1557} and \rprop{P1100}, we have then
 \[
  -\ek*{G(-\ko z)}^\ad
  =\CG+\int_{\blhl}\frac{1+t^2}{-t-z}\etaG(\dif t)
  =\CG+\int_{[-\rep,+\infty)}\frac{1+t^2}{t-z}\tilde\theta(\dif t)
 \]
 for all \(z\in\C\setminus[-\rep,+\infty)\), where \(\tilde\theta\) is the image measure of \(\etaG\) under the transformation \(T\colon\blhl\to[-\rep,+\infty)\) defined by \(T(t)\defeq-t\). Because of \(\CG\in\Cggq\) and \(\tilde\theta\in\Mgguqa{1}{[-\rep,+\infty)}\), \rthm{T1542} yields then that \(F\colon\C\setminus[-\rep,+\infty)\to\Cqq\) defined by \(F(z)\defeq-\ek{G(-\ko z)}^\ad\) belongs to \(\SFuu{q}{-\rep}\) and that \((\CF,\etaF)=(\CG,\tilde\theta)\).
 \eenui
\erem

Now we get an integral representation for functions which belong to the class \(\TFqrep\).
\bthml{T1332}
 Let \(\rep \in\R\) and let \(G\colon\Ct  \to \Cqq \). Then:
 \benui
 \il{T1332.a} If \(G\in \TFqrep \), then there are a unique matrix \(\gamma \in \Cggq\) and a unique \tnnH{} measure \(\mu \in\Mggqblhl \) such that
 \beql{T1332.B1}
  G(z)
  =-\gamma +\int_{\blhl} \frac{1+\rep-t}{t-z} \mu (\dif t).
 \eeq
 holds for each \(z\in\Ct\). Furthermore, \(\gamma=\CG\) and \(\gamma=-\lim_{y\to+\infty}G\rk{\iu y}\).
 \il{T1332.b} If there are a matrix \(\gamma\in\Cggq\) and a \tnnH{} measure \(\mu \in\Mggqrhl \) such that \(G\) can be represented via \eqref{T1332.B1} for each \(z\in \Ct\), then \(G\) belongs to the class \(\TFqrep \).
 \eenui
\ethm
\bproof
 \eqref{T1332.a} According to \rremp{R0941}{R0941.b}, the function \(F\colon\C\setminus[-\rep,+\infty)\to\Cqq\) defined by \(F(z)\defeq-\ek{G(-\ko z)}^\ad\) belongs to \(\SFuu{q}{-\rep}\) and \(\CF=\CG\). \rthmp{368}{368.a} yields then that there is a unique pair \((\delta,\rho)\) from \(\Cggq\times\Mggqa{[-\rep,+\infty)}\) such that
 \beql{T1332.1}
  F(z)
  =\delta+\int_{[-\rep,+\infty)}\frac{1+t+\rep}{t-z}\rho(\dif t)
 \eeq
 for all \(z\in\C\setminus[-\rep,+\infty)\) and that \(\delta=\CF\). Applying \rprop{P1100}, we now infer
 \[
  \begin{split}
   G(z)
   &=-\ek*{F(-\ko z)}^\ad
   =-\CF-\int_{[-\rep,+\infty)}\frac{1+t+\rep}{t+z}\rho(\dif t)\\
   &=-\CG+\int_{[-\rep,+\infty)}\frac{1+\rep-(-t)}{-t-z}\rho(\dif t)
   =-\gamma +\int_{\blhl} \frac{1+\rep-t}{t-z} \mu (\dif t)
  \end{split}
 \]
 for all \(z\in\Ct\), where \(\gamma\defeq\CG\) and \(\mu\) is the image measure of \(\rho\) under the transformation \(R\colon[-\rep,+\infty)\to\blhl\) defined by \(R(t)\defeq-t\). Since \rprop{R367Ad} yields \(\lim_{y\to+\infty}F(\iu y)=\delta\), we conclude furthermore
 \[
  \gamma
  =\CG
  =\CF
  =\delta
  =\delta^\ad
  =\ek*{\lim_{y\to+\infty}F(\iu y)}^\ad
  =-\lim_{y\to+\infty}G\rk{\iu y}.
 \]
 
 Now let \(\gamma\in\Cggq\) and \(\mu\in\Mggqblhl\) be arbitrary such that \eqref{T1332.B1} holds true for all \(z\in\Ct\). Then using \rprop{P1100} we get
 \beql{T1332.2}
  \begin{split}
   F(z)
   &=-\ek*{G(-\ko z)}^\ad
   =\gamma-\int_{\blhl} \frac{1+\rep-t}{t+z} \mu (\dif t)\\
   &=\gamma+\int_{\blhl} \frac{1-t+\rep}{-t-z} \mu (\dif t)
   =\gamma+\int_{[-\rep,+\infty)} \frac{1+t+\rep}{t-z}\tilde\theta(\dif t)
  \end{split}
 \eeq
 for all \(z\in\C\setminus[-\rep,+\infty)\), where \(\tilde\theta\) is the image measure of \(\mu\) under the transformation \(T\colon\blhl\to[-\rep,+\infty)\) defined by \(T(t)\defeq-t\). Since we know from \rthmp{368}{368.a} that the pair \((\delta,\rho)\in\Cggq\times\Mggqa{[-\rep,+\infty)}\) with \eqref{T1332.1} for all \(z\in\C\setminus[-\rep,+\infty)\) is unique, comparing with \eqref{T1332.2}, we conclude \(\gamma=\delta\) and \(\tilde\theta=\rho\). Hence, \(\gamma=\CF=\CG\) and \(\mu\) is the image measure of \(\rho\) under the transformation \(R\).
 
 \eqref{T1332.b} Let \(\gamma\in\Cggq\) and \(\mu\in\Mggqblhl\) be such that \eqref{T1332.B1} holds true for all \(z\in\Ct\). Then \(F\colon\C\setminus[-\rep,+\infty)\to\Cqq\) defined by \(F(z)\defeq-\ek{G(-\ko z)}^\ad\) fulfills \eqref{T1332.2} for all \(z\in\C\setminus[-\rep,+\infty)\), where \(\tilde\theta\) is the image measure of \(\mu\) under the transformation \(T\). \rthmp{368}{368.b} yields then \(F\in\SFuu{q}{-\rep}\). In view of \(G(z)=-\ek{F(-\ko z)}^\ad\) for all \(z\in\Ct\), hence \(G\) belongs to \(\TFqrep \) by virtue of \rrem{R0948}.
\eproof

In the following, if \(\rep\in\R\) and \(G\in \TFqrep \) are given, then we will write \((\gammaG,\muG)\)\index{g@\(\gammaG\)}\index{m@\(\muG\)} for the unique pair \((\gamma,\mu)\) from \(\Cggq\times\Mggqblhl\) which fulfills \eqref{T1332.B1} for all \(z\in\Ct\).

\breml{R1056}
 \benui
  \il{R1056.a} Let \(\lep\in\R\) and let \(F\in \SFqlep\). In view of \rthm{368} and \rprop{P1100}, we have then
 \[
  -\ek*{F(-\ko z)}^\ad
  =-\gammaF+\int_{\rhl}\frac{1-(-t)-\lep}{-t-z}\muF(\dif t)
  =-\gammaF+\int_{(-\infty,-\lep]}\frac{1-\lep-t}{t-z}\hat\theta(\dif t)
 \]
 for all \(z\in\C\setminus(-\infty,-\lep]\), where \(\hat\theta\) is the image measure of \(\muF\) under the transformation \(R\colon\rhl\to(-\infty,-\lep]\) defined by \(R(t)\defeq-t\). Because of \(\gammaF\in\Cggq\), \rthm{T1332} yields then that \(G\colon\C\setminus(-\infty,-\lep]\to\Cqq\) defined by \(G(z)\defeq-\ek{F(-\ko z)}^\ad\) belongs to \(\TFuu{q}{-\lep}\) and that \((\gammaG,\muG)=(\gammaF,\hat\theta)\).
  \il{R1056.b} Let \(\rep\in\R\) and let \(G\in \TFqrep\). In view of \rthm{T1332}, we have then
 \[
  -\ek*{G(-\ko z)}^\ad
  =\gammaG+\int_{\blhl}\frac{1+\rep-t}{-t-z}\muG(\dif t)
  =\gammaG+\int_{[-\rep,+\infty)}\frac{1+t+\rep}{t-z}\tilde\theta(\dif t)
 \]
 for all \(z\in\C\setminus[-\rep,+\infty)\), where \(\tilde\theta\) is the image measure of \(\muG\) under the transformation \(T\colon\blhl\to[-\rep,+\infty)\) defined by \(T(t)\defeq-t\). Because of \(\gammaG\in\Cggq\), \rthm{368} yields then that \(F\colon\C\setminus[-\rep,+\infty)\to\Cqq\) defined by \(F(z)\defeq-\ek{G(-\ko z)}^\ad\) belongs to \(\SFuu{q}{-\rep}\) and that \((\gammaF,\muF)=(\gammaG,\tilde\theta)\).
 \eenui
\erem

\bpropl{P1049}
 Let \(\rep \in\R\) and let \(G\in\TFqrep \). Then:
 \benui
  \il{P1049.a} Let \(z\in\Ct\). Then \(\ko z\in\Ct\) and \(\ek{G(z)}^\ad=G(\ko z)\).
  \il{P1049.b} For all \(z\in \Ct\),
  \begin{align}\label{P1049.AB}
   \Ran{G(z)}&=\ran{\gammau{G}}+\Ran{\muuA{G}{\blhl }},&
   \Nul{G(z)}&=\nul{\gammaG}\cap\Nul{\muuA{G}{\blhl }},
  \end{align}
  and, in particular, \(\ran{\ek{G(z)}^\ad}=\ran{G(z)}\) and \(\nul{\ek{G(z)}^\ad}=\nul{G(z)}\).
  \il{P1049.c} Let \(r\in\NO\). Then the following statements are equivalent:
  \baeqii{0}
   \item \(\rank G(z) = r\) for all \(z\in\Ct\).
   \item There is some \(z_0\in\Ct\) such that \(\rank G(z_0) = r\).
   \item \(\dim[ \ran{\gammaG} + \ran{\muu{G} ( \blhl  }]= r\).
  \eaeqii
 \eenui
\eprop
\bproof
 \eqref{P1049.a} This can be concluded from the representation \eqref{T1332.B1} in \rthmp{T1332}{T1332.a}.
 
 \eqref{P1049.b} According to \rremp{R1056}{R1056.b}, the function \(F\colon\C\setminus[-\rep,+\infty)\to\Cqq\) defined by \(F(z)\defeq-\ek{G(-\ko z)}^\ad\) belongs to \(\SFuu{q}{-\rep}\) and \((\gammaF,\muF)=(\gammaG,\tilde\theta)\), where \(\tilde\theta\) is the image measure of \(\muG\) under the transformation \(T\colon\blhl\to[-\rep,+\infty)\) defined by \(T(t)\defeq-t\). In particular, \(\muua{F}{[-\rep,+\infty)}=\muua{G}{\blhl}\). \rpropp{L374}{L374.a} yields \(\ran{\ek{F(w)}^\ad}=\ran{F(w)}=\ran{\gammau{F}}+\ran{\muua{F}{[-\rep,+\infty)}}\) and \(\nul{\ek{F(w)}^\ad}=\nul{F(w)}=\nul{\gammau{F}}\cap\nul{\muua{F}{[-\rep,+\infty)}}\) for all \(w\in\C\setminus[-\rep,+\infty)\). We now infer
 \[
  \Ran{G(z)}
  =\Ran{\ek*{F(-\ko z)}^\ad}
  =\ran{\gammaF}+\Ran{\muuA{F}{[-\rep,+\infty)}}
  =\ran{\gammaG}+\Ran{\muuA{G}{\blhl}}
 \]
 and
 \[
  \Nul{G(z)}
  =\Nul{\ek*{F(-\ko z)}^\ad}
  =\nul{\gammaF}\cap\Nul{\muuA{F}{[-\rep,+\infty)}}
  =\nul{\gammaG}\cap\Nul{\muuA{G}{\blhl }}
 \]
 for all \(z\in\Ct\).  From~\eqref{P1049.a} and \eqref{P1049.AB} we get \(\ran{\ek{G(z)}^\ad}=\ran{G(z)}\) and \(\nul{\ek{G(z)}^\ad}=\nul{G(z)}\).
 
 \eqref{P1049.c} This is an immediate consequence of~\eqref{P1049.b}.
\eproof

\bcorl{C1344}
 Let \(\rep \in\R\), let \(G\in\TFqrep \), and let \(z_0\in\Ct\). Then \(G(z_0)=\Oqq\) if and only if \(G\) is identically \(\Oqq\).
\ecor
\bproof
 This is an immediate consequence of \rpropp{P1049}{P1049.c}.
\eproof

\bcorl{C1227}
 Let \(\rep \in\R\), let \(G\in\TFqrep \), and let \(\lambda\in\R\) be such that the matrix \(\gammaG+\lambda\Iq\) is \tnnH{}. Then
 \begin{align}\label{C1227.B1}
  \Ran{G(z)-\lambda\Iq}&=\Ran{G(w)-\lambda\Iq}&
  &\text{and}&
  \Nul{G(z)-\lambda\Iq}&=\Nul{G(w)-\lambda\Iq}
 \end{align}
 for all \(z,w\in\Cs\). In particular, if \(\lambda\geq0\), then \(\lambda\) is an eigenvalue of the matrix \(G(z_0)\) for some \(z_0\in\Cs\) if and only if \(\lambda\) is an eigenvalue of the matrix \(G(z)\) for all \(z\in\Cs\). In this case, the eigenspaces \(\nul{G(z)-\lambda\Iq}\) are independent of \(z\in\Cs\).
\ecor
\bproof
 In view of \rthm{T1332}, we conclude that the function \(F\colon\Cs\to\Cqq\) defined by \(F(z)\defeq G(z)-\lambda\Iq\) belongs to \(\TFqrep\). The application of \rpropp{P1049}{P1049.b} to the function \(F\) yields then \eqref{C1227.B1}. Since the matrix \(\gammaG\) is \tnnH{}, we have \(\gammaG+\lambda\Iq\in\Cggq\) if \(\lambda\geq0\). Thus, the remaining assertions are an immediate consequence of \eqref{C1227.B1}.
\eproof

Now we apply the preceding results to the subclass \(\TdFqrep\) of \(\TFqrep\).

\breml{R1604}
 Let \(\rep\in\R\). From \rthmp{T1332}{T1332.a} we see that the class \(\TdFqrep\) consists of all \(G\in\TFqrep\) with \(\gammaG=\Oqq\).
\erem

\breml{R1609}
 Let \(\rep\in\R\) and let \(G\in\TdFqrep \). In view of \rrem{R1604} and \rpropp{P1049}{P1049.b}, then \(\nul{G(z)}=\nul{\muua{G}{\blhl }}\) and \(\ran{G(z)}=\ran{\muua{G}{\blhl }}\) for all \(z\in\Ct\).
\erem

\section{Characterizations of the class $\TFqrep$}\label{S1356}
While we have discussed the class \(\SFqlep\) in \rsec{S1038}, we give here now the corresponding results for the class \(\TFqrep\).
\bpropl{P1356}
 Let \(\rep\in\R\), let \(G\colon\Ct\to\Cqq\) be a matrix-valued function, and let \(G^\mulz\colon\Ct\to\Cqq\) be defined by \(G^\mulz(z)\defeq(\rep-z)G(z)\). Then \(G\) belongs to \(\TFqrep\) if and only if the following two conditions hold true:
 \bAeqi{0}
  \il{P1356.I} \(G\) is holomorphic in \(\Ct\).
  \il{P1356.II} The matrix-valued functions \(\Rstr_\ohe G\) and \(\Rstr_\ohe G^\mulz\) both belong to \(\RFq\).
 \eAeqi
\eprop
\bproof
 According to \rrem{R0948}, the function \(G\) belongs to \(\TFqrep\) if and only if the function \(F\colon\C\setminus[-\rep,+\infty)\to\Cqq\) defined by \(F(z)\defeq-\ek{G(-\ko z)}^\ad\) belongs to \(\SFuu{q}{-\rep}\). Furthermore, \(G\) is holomorphic in \(\Ct\) if and only if \(F\) is holomorphic in \(\C\setminus[-\rep,+\infty)\). By virtue of \rrem{R1602}, the function \(\Rstr_\ohe G\) belongs to \(\RFq\) if and only if \(\Rstr_\ohe F\) belongs to \(\RFq\). Let \(F^\mulz\colon\C\setminus[-\rep,+\infty)\to\Cqq\) be defined by \(F^\mulz(z)\defeq(z+\rep)F(z)\). Then
 \[
  -\ek*{G^\mulz(-\ko z)}^\ad
  =-\ek*{(\rep+\ko z)G(-\ko z)}^\ad
  =(z+\rep)\rk*{-\ek*{G(-\ko z)}^\ad}
  =F^\mulz(z)
 \]
 for all \(z\in\C\setminus[-\rep,+\infty)\). According to \rrem{R1602}, hence \(\Rstr_\ohe G^\mulz\) belong to \(\RFq\) if and only if \(\Rstr_\ohe F^\mulz\) belong to \(\RFq\). The application of \rprop{M1.47N} completes the proof.
\eproof

For each \(\rep\in\R\), let \(\rherep\defeq\setaa{z\in\C}{\re z>\rep}\)\index{c@\(\rherep\)}.

\bpropl{P1530}
 Let \(\rep \in\R\) and let \(G\colon\Ct\to\Cqq\) be a matrix-valued function. Then \(G\) belongs to \(\TFqrep\) if and only if the following four conditions are fulfilled:
 \bAeqi{0}
  \il{P1530.I} \(G\) is holomorphic in \(\Ct\).
  \il{P1530.II} \(\im G(z)\) is \tnnH{} for all \(z\in\ohe\).
  \il{P1530.III} \(-\im G(z)\) is \tnnH{} for all \(z\in\uhe\).
  \il{P1530.IV} \(-\re G(z)\) is \tnnH{} for all \(z\in\rherep\).
 \eAeqi
\eprop
\bproof
 According to \rrem{R0948}, the function \(G\) belongs to \(\TFqrep\) if and only if the function \(F\colon\C\setminus[-\rep,+\infty)\to\Cqq\) defined by \(F(z)\defeq-\ek{G(-\ko z)}^\ad\) belongs to \(\SFuu{q}{-\rep}\). Furthermore, \(G\) is holomorphic in \(\Ct\) if and only if \(F\) is holomorphic in \(\C\setminus[-\rep,+\infty)\). For all \(z\in\C\setminus[-\rep,+\infty)\), we have \(\re F(z)=-\re G(-\ko z)\) and \(\im F(z)=\im G(-\ko z)\). Hence,
 \begin{align*}
  \im F(\ohe)&=\im G(\ohe),&
  \im F(\uhe)&=\im G(\uhe)&
  &\text{and}&
  \re F(\lhe{-\rep})&=-\re G(\rherep).
 \end{align*}
 The application of \rprop{214-1} completes the proof.
\eproof

\section{The class $\TOFqrep$}\label{S1600}
In \rsec{S1242}, we have studied the class \(\SOFqlep\). The aim of this section is to derive corresponding results for the dual class \(\TOFqrep\).
The following observation establishes the bridge to \rsec{S1242}.

\breml{R1627}
 \benui
  \il{R1627.a} Let \(\lep\in\R\), let \(F\in \SOFqlep\), and let \(G\colon\C\setminus(-\infty,-\lep]\to\Cqq\) be defined by \(G(z)\defeq-\ek{F(-\ko z)}^\ad\). We have then \(\normE{G(\iu y)}=\normE{F(\iu y)}\) for all \(y\in[1,+\infty)\). Taking additionally into account \rrem{R0948}, one can see that \(G\) belongs to \(\TOFuu{q}{-\lep}\).
  \il{R1627.b} Let \(\rep\in\R\), let \(G\in \TOFqrep\), and let \(F\colon\C\setminus[-\rep,+\infty)\to\Cqq\) be defined by \(F(z)\defeq-\ek{G(-\ko z)}^\ad\). We have then \(\normE{F(\iu y)}=\normE{G(\iu y)}\) for all \(y\in[1,+\infty)\). Taking additionally into account \rrem{R0948}, one can see that \(F\) belongs to \(\SOFuu{q}{-\rep}\).
 \eenui
\erem

\bthml{T1610}
 Let \(\rep\in\R\) and let \(G\colon\Ct\to \Cqq \). Then:
 \benui
  \il{T1610.a} If \(G\in \TOFqrep  \), then there is a unique measure \(\sigma\in\Mggqblhl \) such that
  \beql{T1610.B1}
   G(z)
   =\int_{\blhl} \frac{1}{t-z} \sigma(\dif t)
  \eeq
  for each \(z\in\Ct\). Furthermore, \(\sigma(\blhl)=-\iu\lim_{y\to+\infty}yG(\iu y)\).
  \il{T1610.b} If there is a measure \(\sigma\in\Mggqblhl \) such that \(G\) can be represented via \eqref{T1610.B1} for each \(z\in \Ct\), then \(G\) belongs to the class \(\TOFqrep\).
 \eenui
\ethm
\bproof
 \eqref{T1610.a} According to \rremp{R1627}{R1627.b}, the function \(F\colon\C\setminus[-\rep,+\infty)\to\Cqq\) defined by \(F(z)\defeq-\ek{G(-\ko z)}^\ad\) belongs to \(\SOFuu{q}{-\rep}\). \rthmp{T1434}{T1434.a} yields then the existence of a unique \tnnH{} measure \(\tau\in\Mggqa{[-\rep,+\infty)}\) such that
 \beql{T1610.1}
  F(z)
  =\int_{[-\rep,+\infty)}\frac{1}{t-z}\tau(\dif t)
 \eeq
 for all \(z\in\C\setminus[-\rep,+\infty)\). Using \rprop{P1100}, we now infer
 \[
   G(z)
   =-\ek*{F(-\ko z)}^\ad
   =\int_{[-\rep,+\infty)}\frac{1}{-t-z}\tau(\dif t)
   =\int_{\blhl} \frac{1}{t-z} \sigma(\dif t)
 \]
 for all \(z\in\Ct\), where \(\sigma\) is the image measure of \(\tau\) under the transformation \(R\colon[-\rep,+\infty)\to\blhl\) defined by \(R(t)\defeq-t\). Since \rprop{F394} yields \(\tau([-\rep,+\infty))=-\iu \lim_{y\to+\infty}yF(\iu y)\), we conclude furthermore
 \[
  \sigma\rk*{\blhl}
  =\ek*{\sigma\rk*{\blhl}}^\ad
  =\ek*{\tau\rk*{[-\rep,+\infty)}}^\ad
  =\iu\lim_{y\to+\infty}y\ek*{F(\iu y)}^\ad
  =-\iu\lim_{y\to+\infty}yG(\iu y).
 \]
 
 Now let \(\sigma\in\Mggqblhl\) be such that \eqref{T1610.B1} holds true for all \(z\in\Ct\). Using \rprop{P1100}, we get then
 \beql{T1610.2}
   F(z)
   =-\ek*{G(-\ko z)}^\ad
   =\int_{\blhl} \frac{1}{-t-z} \sigma(\dif t)
   =\int_{[-\rep,+\infty)} \frac{1}{t-z}\tilde\theta(\dif t)
 \eeq
 for all \(z\in\C\setminus[-\rep,+\infty)\), where \(\tilde\theta\) is the image measure of \(\sigma\) under the transformation \(T\colon\blhl\to[-\rep,+\infty)\) defined by \(T(t)\defeq-t\). Since we know from \rthmp{T1434}{T1434.a} that the measure \(\tau\in\Mggqa{[-\rep,+\infty)}\) with \eqref{T1610.1} for all \(z\in\C\setminus[-\rep,+\infty)\) is unique, we obtain \(\tilde\theta=\tau\). Hence, \(\sigma\) is the image measure of \(\tau\) under the transformation \(R\).
 
 \eqref{T1610.b} Let \(\sigma\in\Mggqblhl\) be such that \eqref{T1610.B1} holds true for all \(z\in\Ct\). Thus, in view of \rprop{P1100}, the function \(F\colon\C\setminus[-\rep,+\infty)\to\Cqq\) defined by \(F(z)\defeq-\ek{G(-\ko z)}^\ad\) fulfills \eqref{T1610.2} for all \(z\in\C\setminus[-\rep,+\infty)\), where \(\tilde\theta\) is the image measure of \(\sigma\) under the transformation \(T\). \rthmp{T1434}{T1434.b} yields then \(F\in\SOFuu{q}{-\rep}\). Because of \(G(z)=-\ek{F(-\ko z)}^\ad\) for all \(z\in\Ct\), the function \(G\) belongs to \(\TOFqrep \) by virtue of \rremp{R1627}{R1627.a}.
\eproof

If \(\rep\in\R\) and \(\sigma\) is a measure belonging to \(\Mggqblhl \), then we will call the matrix-valued function \(G\colon\Cs \to\Cqq\) which is, for each \(z\in\Ct\), given by \eqref{T1610.B1} the \emph{\tbsto{\rep}{\(\sigma\)}}. If \(G\in\TOFqrep \), then the unique measure \(\sigma\in\Mggqblhl \) which fulfills \eqref{T1610.B1} for each \(z\in\Ct\) is said to be the \emph{\tbsmo{\rep}{\(G\)}} and will be denoted by \(\stm{G}\)\index{s@\(\stm{G}\)}.

\breml{R1023}
 \benui
  \il{R1023.a} Let \(\lep\in\R\) and let \(F\in \SOFqlep\). In view of \rthm{T1434} and \rprop{P1100}, we have then
  \[
   -\ek*{F(-\ko z)}^\ad
   =\int_{\rhl}\frac{1}{-t-z}\stma{F}{\dif t}
   =\int_{(-\infty,-\lep]} \frac{1}{t-z}\hat\theta(\dif t)
  \]
  for all \(z\in\C\setminus(-\infty,-\lep]\), where \(\hat\theta\) is the image measure of \(\stm{F}\) under the transformation \(R\colon\rhl\to(-\infty,-\lep]\) defined by \(R(t)\defeq-t\). \rthm{T1610} yields then that \(G\colon\C\setminus(-\infty,-\lep]\to\Cqq\) defined by \(G(z)\defeq-\ek{F(-\ko z)}^\ad\) belongs to \(\TOFuu{q}{-\lep}\) and that \(\stm{G}=\hat\theta\).
  \il{R1023.b} Let \(\rep\in\R\) and let \(G\in \TOFqrep\). In view of \rthm{T1610} and \rprop{P1100}, we have then
  \[
  -\ek*{G(-\ko z)}^\ad
   =\int_{\blhl}\frac{1}{-t-z}\stma{G}{\dif t}
   =\int_{[-\rep,+\infty)} \frac{1}{t-z}\tilde\theta(\dif t)
  \]
  for all \(z\in\C\setminus[-\rep,+\infty)\), where \(\tilde\theta\) is the image measure of \(\stm{G}\) under the transformation \(T\colon\blhl\to[-\rep,+\infty)\) defined by \(T(t)\defeq-t\). \rthm{T1434} yields then that \(F\colon\C\setminus[-\rep,+\infty)\to\Cqq\) defined by \(F(z)\defeq-\ek{G(-\ko z)}^\ad\) belongs to \(\SOFuu{q}{-\rep}\) and that \(\stm{F}=\tilde\theta\).
 \eenui
\erem

\bpropl{P0929}
 Let \(\rep\in\R\) and let \(G\in\TOFqrep \). Then \(\ran{G(z)}=\ran{\stma{G}{\blhl}}\) and \(\nul{G(z)}=\nul{\stma{G}{\blhl}}\) for all \(z\in\Ct\).
\eprop
\bproof
 According to \rremp{R1023}{R1023.b}, the function \(F\colon\C\setminus[-\rep,+\infty)\to\Cqq\) defined by \(F(z)\defeq-\ek{G(-\ko z)}^\ad\) belongs to \(\SOFuu{q}{-\rep}\) and \(\stm{F}\) is the image measure of \(\stm{G}\) under the transformation \(T\colon\blhl\to[-\rep,+\infty)\) defined by \(T(t)\defeq-t\). In particular, \(\stma{F}{[-\rep,+\infty)}=\stma{G}{\blhl}\). \rprop{F394} yields \(\ran{F(w)}=\ran{\stma{F}{[-\rep,+\infty)}}\) and \(\nul{F(w)}=\nul{\stma{F}{[-\rep,+\infty)}}\) for all \(w\in\C\setminus[-\rep,+\infty)\). Furthermore, \(\ran{\ek{F(w)}^\ad}=\ran{F(w)}\) and \(\nul{\ek{F(w)}^\ad}=\nul{F(w)}\) follow for all \(w\in\C\setminus[-\rep,+\infty)\) from \rpropp{L374}{L374.a}. Finally for all \(z\in\Ct\), we infer
 \[
  \Ran{G(z)}
  =\Ran{\ek*{F(-\ko z)}^\ad}
  =\Ran{F(-\ko z)}
  =\Ran{\stmA{F}{[-\rep,+\infty)}}
  =\Ran{\stmA{G}{\blhl}}
 \]
 and
 \[
  \Nul{G(z)}
  =\Nul{\ek*{F(-\ko z)}^\ad}
  =\Nul{F(-\ko z)}
  =\Nul{\stmA{F}{[-\rep,+\infty)}}
  =\Nul{\stmA{G}{\blhl }}.\qedhere
 \]
\eproof

\section{Moore-Penrose inverses of functions belonging to the class $\TFqrep$}\label{S0950}
This section is the dual counterpart to \rsec{S1044}.

\bpropl{P0832}
 Let \(\rep\in\R\) and let \(G\in\TFqrep\). Then \(G^\mpi\) is holomorphic in \(\Ct\).
\eprop
\bproof
 In view of \rpropp{P1049}{P1049.b}, we obtain the identities \(\nul{G(z)}=\nul{G(\iu)}\) and \(\ran{G(z)}=\ran{G(\iu)}\) for all \(z\in\Ct\). Thus, the application of~\zitaa{MR3014197}{\cprop{8.4}} completes the proof.
\eproof

Let \(\rep\in\R\) and \(G\in\TFqrep\). Then \rlem{P0832} suggests to look if there are functions closely related to \(G^\mpi\) which belong again to \(\TFqrep\). Against to this background, we are led to the function \(F\colon\Ct\to\Cqq\) defined by \(F(z)\defeq-(\rep-z)^\inv\ek{G(z)}^\mpi\).

\bthml{T1507}
 Let \(\rep\in\R\) and let \(G\in\TFqrep\). Then \(F\colon\Ct\to\Cqq\) defined by \(F(z)\defeq-(\rep-z)^\inv\ek{G(z)}^\mpi\) belongs to \(\TFqrep\).
\ethm
\bproof
 According to \rrem{R0948}, the function \(P\colon\C\setminus[-\rep,+\infty)\to\Cqq\) defined by \(P(z)\defeq-\ek{G(-\ko z)}^\ad\) belongs to \(\SFuu{q}{-\rep}\). \rthm{L375} yields then that the function \(Q\colon\C\setminus[-\rep,+\infty)\to\Cqq\) defined by \(Q(z)\defeq-(z+\rep)^\inv\ek{P(z)}^\mpi\) belongs to \(\SFuu{q}{-\rep}\) as well. We now infer
 \beql{T1507.1}
  -\ek*{Q(-\ko z)}^\ad
  =-(-z+\rep)^\inv\rk*{-\ek*{P(-\ko z)}^\mpi}^\ad
  =-(\rep-z)^\inv\rk*{-\ek*{P(-\ko z)}^\ad}^\mpi
  =F(z)
 \eeq
 for all \(z\in\Ct\). Hence, \(F\in\TFqrep\) by virtue of \rrem{R0948}.
\eproof

Now we specify the result of \rthm{T1507} for functions belonging to \(\TOFqrep\).

\bpropl{P1056}
 Let \(\rep\in\R\) and let \(G\in\TOFqrep \). Then \(F\colon\Ct\to\Cqq\) defined by \(F(z)\defeq-(\rep-z )^\inv\ek{G(z)}^\mpi\) belongs to \(\TFqrep\) and \(\gammaF=\ek{\stma{G}{\blhl}}^\mpi\). If \(G\) is not the constant function with value \(\Oqq\), then \(F\in\TFqrep\setminus\TOFqrep\).
\eprop
\bproof
 According to \rremp{R1023}{R1023.b}, the function \(P\colon\C\setminus[-\rep,+\infty)\to\Cqq\) defined by \(P(z)\defeq-\ek{G(-\ko z)}^\ad\) belongs to \(\SOFuu{q}{-\rep}\) and \(\stm{P}\) is the image measure of \(\stm{G}\) under the transformation \(T\colon\blhl\to[-\rep,+\infty)\) defined by \(T(t)\defeq-t\). In particular, \(\stma{P}{[-\rep,+\infty)}=\stma{G}{\blhl}\). \rprop{L395} yields that the function \(Q\colon\C\setminus[-\rep,+\infty)\to\Cqq\) defined by \(Q(z)\defeq-(z+\rep)^\inv\ek{P(z)}^\mpi\) belongs to \(\SFuu{q}{-\rep}\), that \(\gammau{Q}=\ek{\stma{P}{[-\rep,+\infty)}}^\mpi\), and that \(Q\notin\SOFuu{q}{-\rep}\) if \(P\) is not the constant function with value \(\Oqq\). Furthermore, we have \eqref{T1507.1} for all \(z\in\Ct\). Hence, \(F\in\TFqrep\) and \(\gammaF=\gammau{Q}\) by virtue of \rremp{R1056}{R1056.a}. We now infer \(\gammaF=\ek{\stma{G}{\blhl}}^\mpi\). From \eqref{T1507.1} we conclude \(Q(z)=-\ek{F(-\ko z)}^\ad\) for all \(z\in\C\setminus[-\rep,+\infty)\). Since \(G\not\equiv\Oqq\) implies \(P\not\equiv\Oqq\), the proof is complete in view of \rremp{R1627}{R1627.b}.
\eproof

Now we introduce the dual counterpart of the class introduced in \rnota{D0924}.

\bnotal{D1106}
 Let \(\rep\in\R\). Then let \(\TiFqrep\)\index{s@\(\TiFqrep\)} be the class of all matrix-valued functions \(G\colon\Ct\to\Cqq\) which fulfill the following two conditions:
 \bAeqi{0}
  \il{D1106.I} \(G\) is holomorphic in \(\Ct\) with \(\Rstr_\ohe G\in\RFq\).
  \il{D1106.II} For all \(x\in\brhl\), the matrix \(G(x)\) is \tnnH{}.
 \eAeqi
\enota

For the special case \(q=1\) and \(\rep=0\) the class introduced in \rnota{D1106} was studied by Katsnelson~\zita{MR2805421}.
\bexal{E1107}
 Let \(\rep\in\R\) and let \(D,E\in\Cggq\). Then \(G\colon\Ct\to\Cqq\) defined by \(G(z)\defeq D-(\rep-z)E\) belongs to \(\TiFqrep\).
\eexa

\breml{R1154}
 Let \(\lep,\rep\in\R\). Then:
 \benui
  \il{R1154.a} If \(F\in\SiFqlep\), then \(G\colon\Ct\to\Cqq\) defined by \(G(z)\defeq-\ek{F(\lep+\rep-\ko z)}^\ad\) belongs to \(\TiFqrep\).
  \il{R1154.b} If \(G\in\TiFqrep\), then \(F\colon\Cs\to\Cqq\) defined by \(F(z)\defeq-\ek{G(\lep+\rep-\ko z)}^\ad\) belongs to \(\SiFqlep\).
 \eenui
\erem

\bthml{T1209}
 Let \(\rep\in\R\) and let  \(G\colon\Ct\to\Cqq\). Then:
 \benui
  \il{T1209.a} If \(G\in\TiFqrep\), then there are unique \tnnH{} complex \tqqa{matrices} \(D\) and \(E\) and a unique \tnnH{} measure \(\rho\in\Mggqa{(-\infty,\rep)}\) such that
  \beql{T1209.A}
   G(z)
   =D+(\rep-z)\ek*{-E+\int_{(-\infty,\rep)}\frac{1+\rep-t}{t-z}\rho(\dif t)}
  \eeq
  for all \(z\in\Ct\). Furthermore, the function \(Q\colon\Ct\to\Cqq\) defined by \(Q(z)\defeq(\rep-z)^\inv G(z)\) belongs to \(\TFqrep\) with \(D=\muQa{\set{\rep}}\) and \((E,\rho)=(\gammaQ,\Rstr_\Bsa{(-\infty,\rep)}\muQ)\).
  \il{T1209.b} If \(D\in\Cggq\), \(E \in\Cggq\) and \(\rho\in\Mggqa{(-\infty,\rep)}\) are such that \(G\) can be represented via \eqref{T1209.A} for all \(z\in\Ct\), then \(G\) belongs to \(\TiFqrep\).
 \eenui
\ethm
\bproof
 \eqref{T1209.a} Let \(G\in\TiFqrep\). According to \rremp{R1154}{R1154.b}, the function \(F\colon\C\setminus[-\rep,+\infty)\to\Cqq\) defined by \(F(z)\defeq-\ek{G(-\ko z)}^\ad\) belongs to \(\SiFuu{q}{-\rep}\). \rthmp{T1305}{T1305.a} yields then that there is a unique triple \((M,N,\omega)\) from \(\Cggq\times\Cggq\times\Mggqa{(-\rep,+\infty)}\) such that \(F(z)=-M +(z+\rep)\ek{N +\int_{(-\rep,+\infty)}(1+t+\rep)/(t-z)\omega(\dif t)}\)
 for all \(z\in\C\setminus[-\rep,+\infty)\) and that the function \(P\colon\C\setminus[-\rep,+\infty)\to\Cqq\) defined by \(P(z)  \defeq(z+\rep)^\inv F(z)\) belongs to \(\SFuu{q}{-\rep}\) with \(M=\muPa{\set{-\rep}}\) and \((N,\omega)=(\gammaP,\Rstr_\Bsa{(-\rep,+\infty)}\muP)\). Applying \rprop{P1100}, we now infer
 \[
  \begin{split}
   G(z)
   &=-\ek*{F(-\ko z)}^\ad
   =M-(-z+\rep)\ek*{N +\int_{(-\rep,+\infty)}\frac{1+t+\rep}{t+z}\omega(\dif t)}\\
   &=M+(\rep-z)\ek*{-N+\int_{(-\rep,+\infty)}\frac{1+\rep-(-t)}{-t-z}\omega(\dif t)}\\
   &=D+(\rep-z)\ek*{-E+\int_{(-\infty,\rep)}\frac{1+\rep-t}{t-z}\rho(\dif t)}
  \end{split}
 \]
 for all \(z\in\Ct\), where \(D\defeq M\), where \(E\defeq N\), and where \(\rho\) is the image measure of \(\omega\) under the transformation \(R_0\colon(-\rep,+\infty)\to(-\infty,\rep)\) defined by \(R_0(t)\defeq-t\). The function \(Q\) fulfills \(Q(z)=-\ek{P(-\ko z)}^\ad\) for all \(z\in\Ct\). \rremp{R1056}{R1056.a} yields then \(Q\in\TFqrep\) with \(\gammaQ=\gammaP\) and \(\muQ\) being the image measure of \(\muP\) under the transformation \(R\colon[-\rep,+\infty)\to(-\infty,\rep]\) defined by \(R(t)\defeq-t\). Hence, \(D=M=\muPa{\set{-\rep}}=\muQa{\set{\rep}}\) and \(E=N=\gammaP=\gammaQ\) follow. Furthermore, \(\Rstr_\Bsa{(-\infty,\rep)}\muQ\) is the image measure of \(\Rstr_\Bsa{(-\rep,+\infty)}\muP\) under the transformation \(R_0\), implying \(\rho=\Rstr_\Bsa{(-\infty,\rep)}\muQ\). In particular, the triple \((D,E,\rho)\) is unique.
 
 \eqref{T1209.b} Let \(D ,E \in\Cggq\) and let \(\rho\in\Mggqa{(-\infty,\rep)}\) be such that \eqref{T1209.A} holds true for all \(z\in\Ct\). Let \(F\colon\C\setminus[-\rep,+\infty)\to\Cqq\) be defined by \(F(z)\defeq-\ek{G(-\ko z)}^\ad\). Then, using \rprop{P1100}, we get
 \[
  \begin{split}
   F(z)
   &=-\ek*{G(-\ko z)}^\ad
   =-D-(\rep+z)\ek*{-E+\int_{(-\infty,\rep)}\frac{1+\rep-t}{t+z}\rho(\dif t)}\\
   &=-D+(z+\rep)\ek*{E+\int_{(-\infty,\rep)}\frac{1-t+\rep}{-t-z}\rho(\dif t)}\\
   &=-D+(z+\rep)\ek*{E+\int_{(-\rep,\infty)}\frac{1+t+\rep}{t-z}\tilde\theta(\dif t)}
  \end{split}
 \]
 for all \(z\in\C\setminus[-\rep,+\infty)\), where \(\tilde\theta\) is the image measure of \(\rho\) under the transformation \(T_0\colon(-\infty,\rep)\to(-\rep,+\infty)\) defined by \(T_0(t)\defeq-t\). \rthmp{T1305}{T1305.b} yields then \(F\in\SiFuu{q}{-\rep}\). In view of \(G(z)=-\ek{F(-\ko z)}^\ad\) for all \(z\in\Ct\), hence \(G\) belongs to \(\TiFqrep \) by virtue of \rremp{R1154}{R1154.a}.
\eproof

In the following, if \(\rep\in\R\) and \(G\in \TiFqrep \) are given, then we will write \((\DG,\EG,\rhoG)\)\index{d@\(\DG\)}\index{e@\(\EG\)}\index{r@\(\rhoG\)} for the unique triple \((D,E,\rho)\in\Cggq\times\Cggq\times\Mggqa{(-\infty,\rep)}\) which fulfills \eqref{T1209.A} for all \(z\in\Ct\).

\breml{R0911}
 \benui
  \il{R0911.a} Let \(\lep\in\R\) and let \(F\in \SiFqlep\). In view of \rthmp{T1305}{T1305.a} and \rprop{P1100}, we have then
  \[\begin{split}
   -\ek*{F(-\ko z)}^\ad
   &=\DF-(-z-\lep)\ek*{\EF+\int_{(\lep,+\infty)}\frac{1+t-\lep}{t+z}\rhoF(\dif t)}\\
   &=\DF+(-\lep-z)\ek*{-\EF+\int_{(\lep,+\infty)}\frac{1-\lep-\rk{-t}}{-t-z}\rhoF(\dif t)}\\
   &=\DF+(-\lep-z)\ek*{-\EF+\int_{(-\infty,-\lep)}\frac{1-\lep-t}{t-z}\hat\theta(\dif t)}
  \end{split}\]
  for all \(z\in\C\setminus(-\infty,-\lep]\), where \(\hat\theta\) is the image measure of \(\rhoF\) under the transformation \(R_0\colon(\lep,+\infty)\to(-\infty,-\lep)\) defined by \(R_0(t)\defeq-t\). Because of \(\DF\in\Cggq\) and \(\EF\in\Cggq\), \rthm{T1209} yields then that \(G\colon\C\setminus(-\infty,-\lep]\to\Cqq\) defined by \(G(z)\defeq-\ek{F(-\ko z)}^\ad\) belongs to \(\TiFuu{q}{-\lep}\) and that \((\DG,\EG,\rhoG)=(\DF,\EF,\hat\theta)\).
  \il{R0911.b} Let \(\rep\in\R\) and let \(G\in \TiFqrep\). In view of \rthmp{T1209}{T1209.a} and \rprop{P1100}, we have then
  \[\begin{split}
   -\ek*{G(-\ko z)}^\ad
   &=-\DG-(\rep+z)\ek*{-\EG+\int_{(-\infty,\rep)}\frac{1+\rep-t}{t+z}\rhoG(\dif t)}\\
   &=-\DG+(z+\rep)\ek*{\EG+\int_{(-\infty,\rep)}\frac{1-t+\rep}{-t-z}\rhoG(\dif t)}\\
   &=-\DG+(z+\rep)\ek*{\EG+\int_{(-\rep,+\infty)}\frac{1+t+\rep}{t-z}\tilde\theta(\dif t)}
  \end{split}\]
  for all \(z\in\C\setminus[-\rep,+\infty)\), where \(\tilde\theta\) is the image measure of \(\rhoG\) under the transformation \(T_0\colon(-\infty,\rep)\to(-\rep,+\infty)\) defined by \(T_0(t)\defeq-t\). Because of \(\DG\in\Cggq\) and \(\EG\in\Cggq\), \rthm{T1305} yields then that \(F\colon\C\setminus[-\rep,+\infty)\to\Cqq\) defined by \(F(z)\defeq-\ek{G(-\ko z)}^\ad\) belongs to \(\SiFuu{q}{-\rep}\) and that \((\DF,\EF,\rhoF)=(\DG,\EG,\tilde\theta)\).
 \eenui
\erem

\bcorl{C0938}
 Let \(\rep\in\R\). Then \(G\colon\Ct\to\Cqq\) belongs to \(\TiFqrep\) if and only if \(Q\colon\Ct\to\Cqq\) defined by \(Q(z)=(\rep-z)^\inv G(z)\) belongs to \(\TFqrep\).
\ecor
\bproof
 From \rrem{R0911} we conclude that \(G\) belongs to \(\TiFqrep\) if and only if \(F\colon\C\setminus[-\rep,+\infty)\to\Cqq\) defined by \(F(z)\defeq-\ek{G(-\ko z)}^\ad\) belongs to \(\SiFuu{q}{-\rep}\). \rcor{C1611} yields that \(F\) belongs to \(\SiFuu{q}{-\rep}\) if and only if \(P\colon\C\setminus[-\rep,+\infty)\to\Cqq\) defined by \(P(z)=(z+\rep)^\inv F(z)\) belongs to \(\SFuu{q}{-\rep}\). Since \(Q(z)=-\ek{P(-\ko z)}^\ad\) holds true for all \(z\in\Ct\), we see by virtue of \rrem{R1056} that \(P\) belongs to \(\SFuu{q}{-\rep}\) if and only if \(Q\) belongs to \(\TFqrep\), which completes the proof.
\eproof

\bcorl{C0851}
 Let \(\rep\in\R\) and let \(G\in\TiFqrep\). For all \(x_1,x_2\in\brhl\) with \(x_1\leq x_2\), then \(\Oqq\leq G(x_1)\leq G(x_2)\).
\ecor
\bproof
 Using \rthm{T1209}, we obtain \(G(x_2)-G(x_1)=(x_2-x_1)\rk{\EG+\int_{\blhl}\ek{(1+\rep-t)(\rep-t)}/\ek{(t-x_2)(t-x_1)}\rhoGa{\dif t}}\) for all \(x_1,x_2\in\brhl\) with \(x_1\leq x_2\), by direct calculation. Since \(\EG\in\Cggq\) and \(G(x)\in\Cggq\) for all \(x\in\brhl\), thus the proof is complete.
\eproof

\bpropl{P0955}
 Let \(\rep \in\R\) and let \(G\in\TiFqrep \). Then:
 \benui
  \il{P0955.a} If \(z\in\Ct\), then \(\ek{G(z)}^\ad=G(\ko z)\).
  \il{P0955.b} For all \(z\in \Ct \),
  \begin{align}
   \Ran{G(z)}&=\ran{\DG}+\ran{\EG}+\Ran{\rhoGA{(-\infty,\rep)}},\label{P0955.A}\\
   \Nul{G(z)}&=\nul{\DG}\cap\nul{\EG}\cap\Nul{\rhoGA{(-\infty,\rep)}},\label{P0955.B}
  \end{align}
  and, in particular, \(\ran{\ek{G(z)}^\ad}=\ran{G(z)}\) and \(\nul{\ek{G(z)}^\ad}=\nul{G(z)}\).
  \il{P0955.c} Let \(r\in\NO\). Then the following statements are equivalent:
  \baeqii{0}
   \item \(\rank  G(z) = r\) for all \(z\in\Ct \).
   \item There is some \(z_0\in\Ct \) such that \(\rank  G(z_0) = r\).
   \item \(\dim[\ran{\DG}+\ran{\EG}+\ran{\rhoGa{(-\infty,\rep)}}]= r\).
  \eaeqii
 \eenui
\eprop
\bproof
 \eqref{P0955.a} This can be concluded from the representation \eqref{T1209.A} in \rthmp{T1209}{T1209.a}.
 
 \eqref{P0955.b} According to \rremp{R0911}{R0911.b}, the function \(F\colon\C\setminus[-\rep,+\infty)\to\Cqq\) defined by \(F(z)\defeq-\ek{G(-\ko z)}^\ad\) belongs to \(\SiFuu{q}{-\rep}\) and \((\DF,\EF,\rhoF)=(\DG,\EG,\tilde\theta)\), where \(\tilde\theta\) is the image measure of \(\rhoG\) under the transformation \(T_0\colon(-\infty,\rep)\to(-\rep,+\infty)\) defined by \(T_0(t)\defeq-t\). In particular, \(\rhoFa{(-\rep,+\infty)}=\rhoGa{(-\infty,\rep)}\). \rpropp{P1622}{P1622.b} yields \(\ran{\ek{F(w)}^\ad}=\ran{F(w)}=\ran{\DF}+\ran{\EF}+\ran{\rhoFa{(-\rep,+\infty)}}\) and \(\nul{\ek{F(w)}^\ad}=\nul{F(w)}=\nul{\DF}\cap\nul{\EF}\cap\nul{\rhoFa{(-\rep,+\infty)}}\) for all \(w\in\C\setminus[-\rep,+\infty)\). We infer
 \[\begin{split}
  \Ran{G(z)}
  =\ran{\ek*{F(-\ko z)}^\ad}
  &=\ran{\DF}+\ran{\EF}+\Ran{\rhoFA{(-\rep,+\infty)}}\\
  &=\ran{\DG}+\ran{\EG}+\Ran{\rhoGA{(-\infty,\rep)}}
 \end{split}\]
 and, analogously, \(\nul{G(z)}=\nul{\DG}\cap\nul{\EG}\cap\nul{\rhoGa{(-\infty,\rep)}}\)
 for all \(z\in\Ct\). 
 In view of~\eqref{P0955.a}, \rpart{P0955.b} is proved.
 
 \eqref{P0955.c} This is an immediate consequence of~\eqref{P0955.b}.
\eproof

\bcorl{C1029}
 Let \(\rep \in\R\), let \(G\in\TiFqrep \), and let \(z_0\in\Ct\). Then \(G(z_0)=\Oqq\) if and only if \(G(z)=\Oqq\) for all \(z\in\Ct\).
\ecor
\bproof
 This is an immediate consequence of \rpropp{P0955}{P0955.c}.
\eproof

\bcorl{C1034}
 Let \(\rep \in\R\), let \(G\in\TiFqrep \), and let \(\lambda\in\R\) be such that the matrix \(\DG-\lambda\Iq\) is \tnnH{}. Then
 \begin{align}\label{C1034.B1}
  \Ran{G(z)-\lambda\Iq}&=\Ran{G(w)-\lambda\Iq}&
  &\text{and}&
  \Nul{G(z)-\lambda\Iq}&=\Nul{G(w)-\lambda\Iq}
 \end{align}
 for every choice of \(z\) and \(w\) in \(\Ct\). In particular, if \(\lambda\leq0\), then \(\lambda\) is an eigenvalue of the matrix \(G(z_0)\) for some \(z_0\in\Ct\) if and only if \(\lambda\) is an eigenvalue of the matrix \(G(z)\) for all \(z\in\Ct\). In this case, the eigenspaces \(\nul{G(z)-\lambda\Iq}\) are independent of \(z\in\Ct\).
\ecor
\bproof
 In view of \rthm{T1209}, we conclude that the function \(F\colon\Ct\to\Cqq\) defined by \(F(z)\defeq G(z)-\lambda\Iq\) belongs to \(\TiFqrep\). The application of \rpropp{P0955}{P0955.b} to the function \(F\) yields then \eqref{C1034.B1}. Since the matrix \(\DG\) is \tnnH{}, we have \(\DG-\lambda\Iq\in\Cggq\) if \(\lambda\leq0\). Thus, the remaining assertions are an immediate consequence of \eqref{C1034.B1}.
\eproof

\bleml{L1425}
 Let \(\rep\in\R\) and \(G\in\TiFqrep\). Then \(G^\mpi\) is holomorphic in \(\Ct\).
\elem
\bproof
 In view of formulas \eqref{P0955.A} and \eqref{P0955.B}, we obtain \(\nul{G(z)}=\nul{G(\iu)}\) and \(\ran{G(z)}=\ran{G(\iu)}\) for all \(z\in\Ct\). Thus, the application of~\zitaa{MR3014197}{\cprop{8.4}} completes the proof.
\eproof

The following result is an analogue of \rthm{P0931}.
\bthml{P1101}
 Let \(\rep\in\R\) and \(G\colon\Ct\to\Cqq\) be a matrix-valued function. Then \(G\) belongs to \(\TiFqrep\) if and only if \(F\defeq-G^\mpi\) belongs to \(\TFqrep\).
\ethm
\bproof
 From \rrem{R0911} we get that \(G\) belongs to \(\TiFqrep\) if and only if \(Q\colon\C\setminus[-\rep,+\infty)\to\Cqq\) defined by \(Q(z)\defeq-\ek{G(-\ko z)}^\ad\) belongs to \(\SiFuu{q}{-\rep}\). \rthm{P0931} yields that \(Q\) belongs to \(\SiFuu{q}{-\rep}\) if and only if \(P\defeq-Q^\mpi\) belongs to \(\SFuu{q}{-\rep}\). Since \(F(z)=-\ek{P(-\ko z)}^\ad\) is true for all \(z\in\Ct\), we see from \rrem{R1056} that \(P\) belongs to \(\SFuu{q}{-\rep}\) if and only if \(F\) belongs to \(\TFqrep\), which completes the proof.
\eproof

\appendix
\section{Some considerations on \tnnH{} measures}\label{A1608}
In this appendix, we summarize some facts on integration with respect to \tnnH{} measures. For each \tnnH{} \tqqa{measure} \(\mu = (\mu_{jk})^q_{j,k=1}\) on a measurable space \((\Omega,\mathfrak{A})\), we denote by \(\Loaaaa{1}{\Omega}{\gA}{\mu}{\C}\)\index{l@\(\Loaaaa{1}{\Omega}{\gA}{\mu}{\C}\)} the set of all \(\gA\)\nobreakdash-\(\BsaC\)\nobreakdash-measurable functions \(f\colon\Omega\to\C\) such that the integral \(\int_\Omega f\dif\mu\) exists, i.\,e.\ that \(\int_\Omega\abs{f}\dif\tilde\mu_{jk}<\infty\) for every choice of \(j,k\in\mn{1}{q}\), where \(\tilde\mu_{jk}\) is the variation of the complex measure \(\mu_{jk}\).

For each \(A\in\Cqq\), let \(\tr A\)\index{t@\(\tr A\)} be the trace of \(A\).

\breml{R1643}
 Let \(\mu\) be a \tnnH{} measure on a measurable space \((\Omega,\mathfrak{A})\), let \(\tau\defeq\tr\mu\) be the trace measure of \(\mu\), and let \(f\colon\Omega\to\C\) be an \(\mathfrak{A}\)\nobreakdash-\(\BsaC\)\nobreakdash-measurable function. Then \(f\) belongs to \(\LaaaC{\Omega}{\mathfrak{A}}{\mu}\) if and only if \(f\) belongs to \(\LaaaC{\Omega}{\mathfrak{A}}{\tau}\).
\erem

\breml{R1141}
 Let \(\mu\) be a \tnnH{} \tqqa{measure} on a measurable space \((\Omega,\mathfrak{A})\) and let \(u\in\Cq\). Then \(\nu\defeq u^\ad\mu u\) is a finite measure on \((\Omega,\mathfrak{A})\) which is absolutely continuous with respect to the trace measure of \(\mu\). If \(f\) belongs to \(\Loaaaa{1}{\Omega}{\gA}{\mu}{\C}\), then \(\int_\Omega\abs{f}\dif\nu<\infty\) and \(\int_\Omega f\dif\nu=u^\ad\rk{\int_\Omega f\dif\mu}u\).
\erem

\breml{R1558}
 Let \(\mu\) be a \tnnH{} \tqqa{measure} on a measurable space \((\Omega,\mathfrak{A})\). An \(\mathfrak{A}\)\nobreakdash-\(\BsaC\)\nobreakdash-measurable function \(f\colon\Omega\to\C\) belongs to \(\Loaaaa{1}{\Omega}{\gA}{\mu}{\C}\) if and only if \(\int_\Omega\abs{f}\dif(u^\ad\mu u)<\infty\) for all \(u\in\Cq\).
\erem

\blemnl{cf.~\zitaa{MR2988005}{\clem{B.2}}}{62B}
 Let \((\Omega,\gA )\) be a measurable space, let \(\sigma\) be a \tnnH{} \tqqa{measure} on \((\Omega,\gA )\), and let \(\tau\) be the trace measure of \(\sigma\). Then:
 \benui
  \il{62B.c} If \(f\in\Loaaaa{1}{\Omega}{\gA}{\sigma}{\C}\), then \(\ran{\int_{\Omega} f\dif\sigma}\subseteq\ran{\sigma (\Omega)}\) and \(\nul{\sigma (\Omega)}\subseteq\nul{\int_{\Omega} f\dif\sigma}\).
  \il{62B.d} If \(f\in\Loaaaa{1}{\Omega}{\gA}{\sigma}{\C}\) fulfills \(\tau\rk{\{f\notin(0,+\infty)\}}=0\), then \(\ran{\int_{\Omega} f\dif\sigma}=\ran{\sigma (\Omega)}\) and \(\nul{\sigma (\Omega)}=\nul{\int_{\Omega} f\dif\sigma}\).
 \eenui
\elem

\bpropnl{\zitaa{MR3133464}{\cprop{B.1}}}{P1100}
 Let \((\Omega,\mathfrak{A})\) and \((\tilde\Omega,\tilde{\mathfrak{A}})\) be measurable spaces and let \(\mu\) be a \tnnH{} \tqqa{measure} on \((\Omega,\mathfrak{A})\). Further, let \(T\colon\Omega\to\tilde\Omega\) be an \(\mathfrak{A}\)\nobreakdash-\(\tilde{\mathfrak{A}}\)\nobreakdash-measurable mapping. Then \(T(\mu)\colon\tilde{\mathfrak{A}}\to\Cqq\) defined by \([T(\mu)](\tilde A)\defeq\mu(T^\inv(\tilde A))\) is a \tnnH{} \tqqa{measure} on \((\tilde\Omega,\tilde{\mathfrak{A}})\). Furthermore, if \(\tilde f\colon\tilde\Omega\to\C\) is an \(\tilde{\mathfrak{A}}\)\nobreakdash-\(\Bori{\C}\)\nobreakdash-measurable mapping, then \(\tilde f\in\LaaaC{\tilde\Omega}{\tilde{\mathfrak{A}}}{T(\mu)}\) if and only if \(\tilde f\circ T\in\LaaaC{\Omega}{\mathfrak{A}}{\mu}\). If \(\tilde f\) belongs to \(\LaaaC{\tilde\Omega}{\tilde{\mathfrak{A}}}{T(\mu)}\), then \(\int_{\tilde A}\tilde f\dif[T(\mu)]=\int_{T^\inv(\tilde A)}(\tilde f\circ T)\dif\mu\) for all \(\tilde A\in\tilde{\mathfrak{A}}\).
\eprop

\bpropnl{Lebesgue's dominated convergence for \tnnH{} measures}{P1102}
 Let \(\mu\) be a \tnnH{} \tqqa{measure} on a measurable space \((\Omega,\mathfrak{A})\) with trace measure \(\tau\). For all \(n\in\N\), let \(f_n\colon\Omega\to\C\) be an \(\mathfrak{A}\)\nobreakdash-\(\BsaC\)\nobreakdash-measurable function. Let \(f\colon\Omega\to\C\) be an \(\mathfrak{A}\)\nobreakdash-\(\BsaC\)\nobreakdash-measurable function and let \(g\in\LaaaC{\Omega}{\mathfrak{A}}{\mu}\) be such that \(\lim_{n\to+\infty}f_n(\omega) = f(\omega)\) for  \(\tau\)\nobreakdash-a.\,a.\ \(\omega\in\Omega\) and that \(\abs{f_n(\omega)}\leq\abs{g(\omega)}\) for all \(n\in\N\) and \(\tau\)\nobreakdash-a.\,a.\ \(\omega\in\Omega\). Then \(f\in\LaaaC{\Omega}{\mathfrak{A}}{\mu}\), \(f_n\in\LaaaC{\Omega}{\mathfrak{A}}{\mu}\) for all \(n\in\N\), and \( \lim_{n\to+\infty}\int_\Omega f_n\dif\mu=\int_\Omega f\dif\mu\).
\eprop
\bproof
 Let \(u\in\Cq\) and let \(\nu\defeq u^\ad\mu u\). According to \rrem{R1141}, then \(\nu\) is a finite measure on \((\Omega,\mathfrak{A})\) which is absolutely continuous with respect to \(\tau\) and \(\int_\Omega\abs{g}\dif\nu<\infty\). In particular, \(\lim_{n\to+\infty}f_n(\omega) = f(\omega)\) for  \(\nu\)\nobreakdash-a.\,a.\ \(\omega\in\Omega\) and \(\abs{f_n(\omega)}\leq\abs{g(\omega)}\) for all \(n\in\N\) and \(\nu\)\nobreakdash-a.\,a.\ \(\omega\in\Omega\). Thus, Lebesgue's dominated convergence theorem provides us \(\int_\Omega\abs{f}\dif\nu<\infty\), \(\int_\Omega\abs{f_n}\dif\nu<\infty\) for all \(n\in\N\), and \(\lim_{n\to+\infty}\int_\Omega f_n\dif\nu=\int_\Omega f\dif\nu\). Since \(u\in\Cq\) was arbitrarily chosen, \(f\in\LaaaC{\Omega}{\mathfrak{A}}{\mu}\) and \(f_n\in\LaaaC{\Omega}{\mathfrak{A}}{\mu}\) for all \(n\in\N\) follow by virtue of \rrem{R1558}. Thus, using \rrem{R1141}, we get \( \lim_{n\to+\infty}u^\ad\rk{\int_\Omega f_n\dif\mu}u=u^\ad\rk{\int_\Omega f\dif\mu}u\) for all \(u\in\Cq\). Hence, \(\lim_{n\to+\infty}\int_\Omega f_n\dif\mu=\int_\Omega f\dif\mu\) holds true.
\eproof

\breml{62C}
 If \(A\in\Cqq \) is such that \(\im  A\in\Cggq\), then \(\nul{A}\subseteq\nul{\im  A}\) (see, \eg{}~\zitaa{MR3014198}{\clem{A.10}}).
\erem

\bleml{63A}
 Let \(\Omega\) be a \tne{} closed subset of \(\R \), let \(\sigma\in\Mggqa{\Omega}\), and let \(\tau\) be the trace measure of \(\sigma\). Then:
 \benui
  \il{63A.a} For each \(w\in\C\setminus \Omega\), the function \(g_w \colon\Omega\to\C\)\index{g@\(g_w\)} defined by \(g_w (t)\defeq1/(t-w)\) belongs to \(\Loaaaa{1}{\Omega}{\BsaO}{\sigma}{\C}\).
  \il{63A.b} The matrix-valued function \(S\colon\C\setminus \Omega\to\Cqq \) given by \(S(w)\defeq\int_\Omega g_w \dif\sigma\) satisfies \(\ran{S(z)}=\ran{\sigma (\Omega)}\) and, in particular, \(\rank S(z) =\rank\sigma (\Omega)\) for each  \(z\in\C\setminus [\inf \Omega, \sup \Omega]\).
 \il{63A.c} \(-\iu\lim_{y\to+\infty}yS(\iu y)=\sigma(\Omega)\).
 \eenui
\elem
\bproof
 \rPart{63A.a} is readily checked. In particular, the function \(S\) is well defined. Let \(z\in\C\setminus [\inf\Omega, \sup\Omega]\). From \rlemp{62B}{62B.c} we get then
 \[
  \Ran{S(z)}
  =\Ran{\int_{\Omega} g_z\dif\sigma} 
  \subseteq \Ran{\sigma (\Omega)}.
 \]
 Thus, in order to prove \rpart{63A.b}, it remains to verify that \(\ran{\sigma (\Omega)} \subseteq\ran{S(z)}\), \ie{}, that
 \beql{Nr.DD}
  \Nul{S^\ad (z)}
  \subseteq \Nul{\sigma (\Omega)}
 \eeq
 is true. \rPart{63A.a} implies
 \beql{Nr.D4}
  \im  g_z
  \in \Loaaaa{1}{\Omega}{\BsaO}{\sigma}{\C}.
 \eeq
 For each \(t\in\Omega\), one can easily check the equations
 \begin{align}\label{Nr.6-18}
  g_z (t)&= \frac{t-\ko z }{\abs{t-z}^2},&
  \ko{g_z} (t)&=\frac{t-z}{\abs{t-z}^2},&
  &\text{and}&
  \im  g_z (t)&= \frac{\im  z}{\abs{t-z}^2}.%
 \end{align}

 First we consider the case that \(z\) belongs to \(\ohe \). In view of \(\im\rk{-S^\ad (z)}=\im S(z)=\int_\Omega\im g_z\dif\sigma\) and \eqref{Nr.6-18}, we see that the matrix \(\im\rk{-S^\ad (z)}\) is \tnnH{}. \rrem{62C} and \eqref{Nr.6-18} yield then
 \beql{Nr.D3}
  \Nul{S^\ad (z)}
  =\Nul{-S^\ad (z)}
  \subseteq \Nul{\im\rk*{-S^\ad (z)}}
  = \Nul{ \int_\Omega \im  g_z \dif\sigma}.
 \eeq
 From \eqref{Nr.6-18} we know that  \(\tau( \set{ \im  g_z\in (-\infty,0] }) = \tau (\emptyset) = 0\). Taking into account \eqref{Nr.D4} and \rlemp{62B}{62B.d}, then
 \beql{Nr.D10}
  \Nul{ \int_\Omega \im  g_z \dif\sigma}
  = \Nul{\sigma (\Omega)}
 \eeq
 follows. Combining \eqref{Nr.D3} and \eqref{Nr.D10} we obtain \eqref{Nr.DD}.

 Now we study the case that \(z\) belongs to \(\uhe\). By virtue of
 \beql{Nr.FGZ}
  \im  S^\ad (z)
  = -\im  S(z)
  = \int_\Omega (-\im  g_z) \dif\sigma
 \eeq
 and \eqref{Nr.6-18}, the matrix \(\im  S^\ad (z)\) belongs to \(\Cggq\). Hence, \rrem{62C} and \eqref{Nr.FGZ} yield
 \beql{Nr.FN}
  \Nul{S^\ad (z)}
  \subseteq\Nul{\im  S^\ad (z)}
  = \Nul{ \int_\Omega (-\im  g_z) \dif\sigma}.
 \eeq
 Because of \eqref{Nr.6-18}, we have \(\tau( \set{ - \im  g_z\in (-\infty, 0]}) = \tau (\emptyset) = 0\). Thus, using \eqref{Nr.D4} and \rlemp{62B}{62B.d}, we have \(\nul{\int_\Omega (-\im  g_z) \dif\sigma}=\nul{\sigma (\Omega)}\). Taking into account \eqref{Nr.FN}, this implies \eqref{Nr.DD}.

 Now we discuss the case that \(\inf \Omega > -\infty\) and that \(z\in (-\infty, \inf \Omega)\). In view of \(t\in\Omega\) and \eqref{Nr.6-18}, we have then \(\tau(\set{ \ko{g_z} \in (-\infty, 0]}) = \tau (\emptyset) = 0\). Since one can see from \rpart{63A.a} that \(\ko{g_z}\) belongs to \(\Loaaaa{1}{\Omega}{\BsaO}{\sigma}{\C}\), application of \rlemp{62B}{62B.d} provides us then 
 \[
  \Nul{S^\ad (z)}
  = \Nul{ \int_\Omega \ko{g_z} \dif\sigma}
  = \Nul{\sigma (\Omega)}.
 \]
 Hence, \eqref{Nr.DD} is valid. Similarly, in the case that \(\sup \Omega < +\infty\) and \(z\in (\sup\Omega, +\infty)\) hold, we get that \eqref{Nr.DD} is fulfilled. Thus \eqref{Nr.DD} is  verified for each \(z\in\C\setminus [\inf \Omega,\sup\Omega]\).

 \eqref{63A.c} For each \(y\in(0,+\infty)\) and each \(t\in\Omega\), we have \(\re\rk{(-\iu y)/(t-\iu y)}=y^2/(t^2+y^2)\), \(\im\rk{(-\iu y)/(t-\iu y)}=-yt/(t^2+y^2)\), and \(\abs{(-\iu y)/(t-\iu y)}\leq1\). From \rprop{P1102} we obtain then
 \[
  \sigma(\Omega)
  =\int_\Omega1\dif\sigma
  =\lim_{y\to+\infty}\int_\Omega\frac{-\iu y}{t-\iu y}\sigma(\dif t)
  =-\iu\lim_{y\to+\infty}yS(\iu y).\qedhere
 \]
\end{proof}

\bibliography{142arxiv}

\def\cprime{$'$}
\begin{thebibliography}{10}

\bibitem{MR2828331}
Y.~Arlinskii, S.~Belyi, and E.~Tsekanovskii.
\newblock {\em Conservative realizations of {H}erglotz-{N}evanlinna functions},
  volume 217 of {\em Operator Theory: Advances and Applications}.
\newblock Birkh\"auser/Springer Basel AG, Basel, 2011.

\bibitem{MR0168769}
N.~Aronszajn and W.~F. Donoghue.
\newblock A supplement to the paper on exponential representations of analytic
  functions in the upper half-plane with positive imaginary part.
\newblock {\em J. Analyse Math.}, 12:113--127, 1964.

\bibitem{MR1105324}
S.~L. Campbell and C.~D. Meyer, Jr.
\newblock {\em Generalized inverses of linear transformations}.
\newblock Dover Publications Inc., New York, 1991.
\newblock Corrected reprint of the 1979 original.

\bibitem{MR2222521}
A.~E. Choque~Rivero, {\relax Yu}.~M. Dyukarev, B.~Fritzsche, and B.~Kirstein.
\newblock A truncated matricial moment problem on a finite interval.
\newblock In {\em Interpolation, {S}chur functions and moment problems}, volume
  165 of {\em Oper. Theory Adv. Appl.}, pages 121--173. Birkh\"auser, Basel,
  2006.

\bibitem{MR1152328}
V.~K. Dubovoj, B.~Fritzsche, and B.~Kirstein.
\newblock {\em Matricial version of the classical {S}chur problem}, volume 129
  of {\em Teubner-Texte zur Mathematik [Teubner Texts in Mathematics]}.
\newblock B. G. Teubner Verlagsgesellschaft mbH, Stuttgart, 1992.
\newblock With German, French and Russian summaries.

\bibitem{MR686076}
{\relax Yu}.~M. Dyukarev.
\newblock Multiplicative and additive {S}tieltjes classes of analytic
  matrix-valued functions and interpolation problems connected with them. {II}.
\newblock {\em Teor. Funktsi\u\i\ Funktsional. Anal. i Prilozhen.},
  (38):40--48, 127, 1982.

\bibitem{MR645305}
{\relax Yu}.~M. Dyukarev and V.~{\`E}. Katsnel{\cprime}son.
\newblock Multiplicative and additive {S}tieltjes classes of analytic
  matrix-valued functions and interpolation problems connected with them. {I}.
\newblock {\em Teor. Funktsi\u\i\ Funktsional. Anal. i Prilozhen.},
  (36):13--27, 126, 1981.

\bibitem{MR752057}
{\relax Yu}.~M. Dyukarev and V.~{\`E}. Katsnel{\cprime}son.
\newblock Multiplicative and additive {S}tieltjes classes of analytic
  matrix-valued functions, and interpolation problems connected with them.
  {III}.
\newblock {\em Teor. Funktsi\u\i\ Funktsional. Anal. i Prilozhen.},
  (41):64--70, 1984.

\bibitem{MR2257838}
J.~Elstrodt.
\newblock {\em Ma\ss- und {I}ntegrationstheorie}.
\newblock Springer-Lehrbuch. [Springer Textbook]. Springer-Verlag, Berlin,
  fourth edition, 2005.
\newblock Grundwissen Mathematik. [Basic Knowledge in Mathematics].

\bibitem{MR3014198}
B.~Fritzsche, B.~Kirstein, A.~Lasarow, and A.~Rahn.
\newblock On reciprocal sequences of matricial {C}arath\'eodory sequences and
  associated matrix functions.
\newblock In {\em Interpolation, {S}chur functions and moment problems. {II}},
  volume 226 of {\em Oper. Theory Adv. Appl.}, pages 57--115.
  Birkh\"auser/Springer Basel AG, Basel, 2012.

\bibitem{MR2988005}
B.~Fritzsche, B.~Kirstein, and C.~M{\"a}dler.
\newblock On matrix-valued {H}erglotz-{N}evanlinna functions with an emphasis
  on particular subclasses.
\newblock {\em Math. Nachr.}, 285(14-15):1770--1790, 2012.

\bibitem{MR3133464}
B.~Fritzsche, B.~Kirstein, and C.~M{\"a}dler.
\newblock Transformations of matricial {$\alpha$}-{S}tieltjes non-negative
  definite sequences.
\newblock {\em Linear Algebra Appl.}, 439(12):3893--3933, 2013.

\bibitem{1151}
B.~Fritzsche, B.~Kirstein, and C.~M{\"a}dler.
\newblock On a simultaneous approach to the even and odd truncated matricial
  {H}amburger moment problems.
\newblock In {\em Recent Advances in Inverse Scattering, {S}chur Analysis and
  Stochastic Processes}, volume 244 of {\em Oper. Theory Adv. Appl.}, pages
  181--285. Birkh\"auser/Springer Basel AG, Basel, 2015.

\bibitem{MR3014199}
B.~Fritzsche, B.~Kirstein, C.~M{\"a}dler, and T.~Schwarz.
\newblock On a {S}chur-type algorithm for sequences of complex {$p\times
  q$}-matrices and its interrelations with the canonical {H}ankel
  parametrization.
\newblock In {\em Interpolation, {S}chur functions and moment problems. {II}},
  volume 226 of {\em Oper. Theory Adv. Appl.}, pages 117--192.
  Birkh\"auser/Springer Basel AG, Basel, 2012.

\bibitem{MR3014197}
B.~Fritzsche, B.~Kirstein, C.~M{\"a}dler, and T.~Schwarz.
\newblock On the concept of invertibility for sequences of complex {$p\times
  q$}-matrices and its application to holomorphic {$p\times q$}-matrix-valued
  functions.
\newblock In {\em Interpolation, {S}chur functions and moment problems. {II}},
  volume 226 of {\em Oper. Theory Adv. Appl.}, pages 9--56.
  Birkh\"auser/Springer Basel AG, Basel, 2012.

\bibitem{MR1784638}
F.~Gesztesy and E.~Tsekanovskii.
\newblock On matrix-valued {H}erglotz functions.
\newblock {\em Math. Nachr.}, 218:61--138, 2000.

\bibitem{KK74}
I.~S. Kats and M.~G. Kre{\u\i}n.
\newblock {$R$}\nobreakdash-functions---analytic functions mapping the upper
  halfplane into itself ({R}ussian).
\newblock Appendix~I in F.~V.~Atkinson.\newblock {\em Diskretnye i nepreryvnye
  granichnye zadachi}.\newblock Translated from the English by I. S. Iohvidov
  and G. A. Karal\cprime nik. Edited and supplemented by I. S. Kats and M. G.
  Kre\u\i n. Izdat. ``Mir'', Moscow, 1968.
\newblock English translation in {\em American {M}athematical {S}ociety
  {T}ranslations, {S}eries 2. {V}ol. 103: {N}ine papers in analysis}.\newblock
  American Mathematical Society, Providence, R.I., 1974, pages~1--18.

\bibitem{MR2805421}
V.~Katsnelson.
\newblock Stieltjes functions and {H}urwitz stable entire functions.
\newblock {\em Complex Anal. Oper. Theory}, 5(2):611--630, 2011.

\bibitem{MR0458081}
M.~G. Kre{\u\i}n and A.~A. Nudel{\cprime}man.
\newblock {\em The {M}arkov moment problem and extremal problems}.
\newblock American Mathematical Society, Providence, R.I., 1977.
\newblock Ideas and problems of P. L. {\v{C}}eby{\v{s}}ev and A. A. Markov and
  their further development, Translated from the Russian by D. Louvish,
  Translations of Mathematical Monographs, Vol. 50.

\bibitem{MR0038923}
H.~Schwerdtfeger.
\newblock {\em Introduction to {L}inear {A}lgebra and the {T}heory of
  {M}atrices}.
\newblock P. Noordhoff, Groningen, 1950.

\end{thebibliography}
\bibliographystyle{abbrv}%

\vfill\noindent
\begin{minipage}{0.5\textwidth}
 Universit\"at Leipzig\\
Fakult\"at f\"ur Mathematik und Informatik\\
PF~10~09~20\\
D-04009~Leipzig
\end{minipage}
\begin{minipage}{0.49\textwidth}
 \begin{flushright}
  \texttt{
   fritzsche@math.uni-leipzig.de\\
   kirstein@math.uni-leipzig.de\\
   maedler@math.uni-leipzig.de
  } 
 \end{flushright}
\end{minipage}

\end{document}